\documentclass[11pt]{article}
\usepackage{amsfonts,latexsym,amsmath,amscd,geometry}
\usepackage{amssymb}
\usepackage{latexsym}

\newcommand \nc{\newcommand}
\newtheorem{theorem}{Theorem}[section]
\newtheorem{lemma}[theorem]{Lemma}

\nc{\ba}{\begin{array}}\nc{\ea}{\end{array}}
\nc{\be}{\begin{eqnarray}}\nc{\ee}{\end{eqnarray}}
\nc{\beq}{\begin{equation}}\nc{\eeq}{\end{equation}}
\nc{\bex}{\begin{eqnarray*}}\nc{\eex}{\end{eqnarray*}}
\nc{\btm}{\begin{theorem}} \nc{\etm}{\end{theorem}}
\nc{\blm}{\begin{lemma}} \nc{\elm}{\end{lemma}}
\nc{\R}{\mathbb{R}}  \nc{\ld}{\lambda}
\nc{\va}{\varphi}
\nc{\ve}{\varepsilon}

\def\pf{\noindent{\bf Proof.\quad}}

\newcommand \qed {\hfill $\Box$}

\begin{document}
\title{Isotropic-Nematic Phase Transition and Liquid Crystal Droplets}
\author{Fanghua Lin, \ Changyou Wang\\
\\
Courant Institute of Mathematical Sciences, New York University, NY 10012\\
Department of Mathematics,
Purdue University
West Lafayette, IN 47907}
\date{}
\maketitle

\begin{abstract}
Liquid crystal droplets are of great interest from physics and applications. Rigorous mathematical analysis is challenging as the problem involves harmonic maps (and in general the Oseen-Frank model), free interfaces and topological defects which could be either inside the droplet or on its surface along with some intriguing boundary anchoring conditions for the orientation configurations. In this paper, through a study of the phase transition between the isotropic and nematic states of liquid crystal based on the Ericksen model, we can show, when the size of droplet is much larger in comparison with the ratio of the Frank constants to the surface tension, a $\Gamma$-convergence theorem for minimizers. This $\Gamma$-limit is in fact the sharp interface limit for the phase transition between the isotropic and nematic regions when the small parameter $\varepsilon$, corresponding to the transition layer width, goes to zero. This limiting process not only provides a geometric description of the shape of the droplet as one would expect, and surprisingly it also gives  the anchoring conditions for the orientations of liquid crystals on the surface of the droplet depending on material constants. In particular, homeotropic, tangential, and even free boundary conditions as assumed in earlier phenomenological modelings arise naturally provided that the surface tension, Frank and Ericksen constants are in suitable ranges.

\end{abstract}

\section{Introduction and statement of results}
\setcounter{equation}{0}
\setcounter{theorem}{0}
\subsection{Isotropic and nematic phase transitions}

Liquid crystal is a state of matter between liquid and solid, in which the molecules tend to align locally in
a common direction and form an anisotropic structure. This orientational order produces an anisotropic complex
fluid with remarkable optical features, which have profound applications in optical and display devices. There
are many phases in liquid crystals including isotropic, nematic, and smectic phases.
Perhaps the most common one is the nematic phase, where the molecules exhibit an orientational order in the absence of translational order. Under the influences of
either external electric-magnetic fields, or thermal changes, or compositional changes, liquid crystals often undergo phase transitions. In the process of phase transitions, there form regions of different phases and thin fast transitional layers around sharp surfaces, across which the nematic order parameter becomes discontinuous.

For simple fluids, the phase transitions between two mixed fluids are usually driven by the interface tensions so that the geometric shape of a sharp free interface takes the form of either an area minimizing surface or surface of constant mean curvature that minimizes areas with volume constraint of enclosed regions. Mathematically,
this phase transition problem has been extensively studied by many authors, including Modica-Mortola \cite{Modica-Mortola}, Modica \cite{Modica}, Sternberg \cite{S1},
Kohn-Sternberg \cite{KS}, through the technique of De Giorgi's Gamma-convergence in the framework of (scalar-valued) Allen-Cahn energy functional with double-well potentials. While Fonseca-Tartar \cite{FT}, Sternberg \cite{S2}, and Andre-Shafrir \cite{AS} studied the gradient theory of phase transitions involving
Allen-Cahn type energy functionals with potential wells of points or curves in $\R^2$. More recently, 
partly motivated by the Keller-Rubinstein-Sternberg problem \cite{RSK}, Lin-Pan-Wang \cite{LPW} have made a systematic study of the vectorial singular perturbation problem of general high dimensional wells,
see also \cite{LW}.

In contrast with simple fluids, the anisotropic structure of liquid crystals implies that both elastic constants of liquid crystal materials in the nematic region and anchoring angles of nematic liquid crystal director fields on the transitional interface will play important roles in determining the shape, possible defect structures and the stability of the interfaces. These are mathematically very challenging problems.
There have been numerous works in the literature, including modeling and experiments, modular simulations and numerics, on phase transitions in liquid crystals by physicists and engineers, see \cite {LT1, LT2}.
The study of the isotropic-nematic interface based on the Ginzburg-Landau-de Gennes (LGdG) theory was initiated in a paper by de Gennes \cite{Gennes}, where the structure of the infinite, flat isotropic-nematic interface in a uniform uniaxial ansatz solution was analyzed. In general, nematic ordering is strongly influenced by confining surfaces, which can impose and favour a particular anchoring condition  (e.g. homeotropic, planar, or oblique anchoring) on the nematic state. It turns out that the relative sizes of different elastic coefficients also play important roles on anchoring conditions on the interface, see Kamil-Bhattacharjee-Adhikari-Menon \cite{KBAM1, KBAM2}. There are proposed forms of the surface energy by Chanderashka \cite{Chan}, Ericksen \cite{E1} based on a phenomenological theory. Besides some earlier works by Virga \cite{V} and Lin-Poon \cite{LP1}, there is an obvious lack of mathematical understanding of these problems.
In \cite{DK}, Dio and Kuzuu studied the structure of the interface between the coexisting isotropic and nematic phases of a lyotropic liquid crystal,
and found an explicit formula of the interfacial tensions in terms of anchoring angles and the length and diameter of the liquid crystal molecule, which favour the planar anchoring. It is known that the nematic structure in the interfacial region can differ substantially from the nematic structure in the bulk. For example, it was shown by
Popa-Nita-Sluckin-Wheeler \cite{PSW} a region proximate to the interface can exhibit biaxiality within the LGdG theory, even if the stable nematic phase is purely uniaxial, provided planar anchoring is enforced on the interface. Such a biaxiality is absent if the anchoring is homeotropic \cite{Gennes}, see also \cite{KBAM1, KBAM2}.

The Landau-Ginzburg-de Gennes model is certainly more flexible and may also be more consistent from both
mathematical and physical point of views. For example, it can be derived rigorously from microscopic (molecular/kinetic) models. It can be used to describe more complex defect structures, both uniaxial and biaxial. In fact, purely uniaxial solutions are very rare in the de Gennes-Ginzburg-Landau model though in many situations they are well approximated by uniaxial ones, see Henno-Majumdar \cite{HM} and Majumdar-Zaenescu \cite{MZ}. The de Gennes-Ginzburg-Landau model in principle may also lead to anisotropic surface energies. This would result in different shapes of droplets and defect patterns within them. But for the analysis for Landau-Ginzburg-de Gennes model, the complexity is formidable. If the energy density functions are quadratic in gradient of Q with coefficients that are quadratic polynomials in Q, then there are $22$ invariants (and $13$ surface terms) along with $4$ null-Lagrangians.  If one would consider additional chiral effects one may need 2-4 additional terms, see \cite{GP} and \cite{LMT}.  On the other hand, with much simplified energy functionals as considered by various authors recently, we believe that the analysis in this current paper can be applied without essential difficulties. There are very few mathematical works on the phase transitions of nematic liquid crystals in high dimensions within the LGdG theory. Let us mention in this direction the work by Park-Wang-Zhang-Zhang \cite{PWZZ} in dimension one, and Golovaty-Novack-Sternberg-Venkatraman \cite{GNSV1, GNSV2} in dimension two.

One can formally derive for these simplified models that in the isotropic-nematic sharp interface transitions, the biaxial property of solutions and their defects contribute only lower orders to the total energy of the system. Thus it does not affect too much the shape of droplet, rather the detailed structure of defects that could be biaxial in nature near defects. This is one of the important reasons that in this article, we shall adopt the reduced Landau-Ginzburg-de Gennes model, or the so-called Ericksen's model for the uniaxial nematic liquid crystals of variable degrees of orientations. This model is natural and relatively simple, and mathematically it is self-contained and consistent. It can accommodate point defects, disinclinations and domain walls in liquid crystals for which rigorous analysis are possible. Moreover, it can keep the classical Oseen-Frank model, favored by many physicists, intact. Since the surface tension in such highly viscous fluids is quite large, the ratio between Frank elastic constants and surface tension is often very small compared to the size of typical droplets. It leads us to the study of the phase transition problems formulated in a form of singular perturbations for the Ericksen energy functional and consequently a sharp interface limiting problem. We will need some $\Gamma$-convergence techniques that the authors have developed for vector-valued variational problems in \cite{LPW, LW}. Here the coexisting isotropic and nematic states are separated by an interface in which order parameters rise from zero on the isotropic side of the interface to saturated, non-zero values on the nematic side. The nematic regions are what we have referred as droplets, and in this way we treat nematic droplets (positioned in an isotropic liquid), and their boundaries are the isotropic-nematic interfaces within the same framework of Ericksen's model. Naturally, one can also study droplets containing isotropic liquid immersed in a volume of nematic liquid crystals.  Of particular interest is that the anchoring conditions for nematic liquid crystal configurations at the boundary of droplets are intrinsically determined by the material constants, and can be derived from this sharp interface limit instead of that it needs to be assumed in phenomenological models.

\subsection{Descriptions of main theorems}
\setcounter{equation}{0}
\setcounter{theorem}{0}
In the framework of Ericksen's theory \cite{E2} (see also \cite{Lin0, Lin}), a nematic liquid crystal configuration is described by a
pair of parameters $(s,n):\Omega\subset\mathbb R^3\mapsto [-\frac12,1]\times \mathbb S^2$, where
$s(x)$ denotes the degree of orientation and $n(x)$ denotes the orientation field at a point $x\in\Omega$.
In particular, at a point $x\in\Omega$ molecules are perfectly aligned in the direction $n(x)$ when $s(x)=1$,
while molecules are perpendicular to $n(x)$ when $s(x)=-\frac12$.

The Ericksen energy density function is assumed to take the form \cite{E2}
$$\mathcal{W}(s, n, \nabla s, \nabla n)=\mathcal{W}_2(s, n, \nabla s, \nabla n)+\mathcal{W}_0(s),$$
where
\begin{eqnarray}\label{w2-1}
&&\mathcal{W}_2(s, n, \nabla s, \nabla n)\nonumber\\
&&= s^2\big[ k_1({\rm{div}} n)^2+k_2(n\cdot {\rm{curl}} n)^2+k_3 |n\wedge {\rm{curl}} n|^2
+(k_2+k_4)({\rm{tr}}(\nabla n)^2-({\rm{div}} n)^2)\big]+\alpha s^2 |\nabla n|^2\nonumber\\
&&\ \ \ +\beta |\nabla s|^2+L_1(\nabla s\cdot n)^2+L_2|\nabla s\wedge n|^2+L_3(\nabla s\cdot n)(s{\rm{div}} n)+L_4 s\nabla s\cdot (\nabla n)n.
\end{eqnarray}
Here $\alpha, \beta>0$, $k_1\ge 0, k_2\ge 0, k_3\ge 0$ are Frank elasticity constants,
and $k_4$ is another constant such that
$$k_2\ge |k_4|\ \ {\rm{and}}\ \ 2k_1\ge k_2+k_4,$$
and $L_1\ge 0, L_2\ge 0, L_3, L_4$ are also constants.
$\mathcal{W}_0(s)$ is a Landau type bulk potential function that dictates isotropic and nematic phases.

It follows from the fact that $|n|=1$ and direct calculations that
$$|n\wedge {\rm{curl}} n|^2=|(\nabla n) n|^2,\
|\nabla s|^2=(\nabla s\cdot n)^2+|\nabla s\wedge n|^2.
$$
Hence, as in Ericksen \cite{E2} and Lin-Poon \cite{LP, LP1}, we can reorganize the expression of $\mathcal{W}_2$ into
the form
\begin{eqnarray}\label{w2-2}
\mathcal{W}_2(s, n, \nabla s, \nabla n)
&=& s^2\big[ \bar{k}_1|{\rm{div}}n|^2+k_2(n\cdot {\rm{curl}}n)^2+\bar{k}_3 |n\wedge {\rm{curl}}n|^2\nonumber\\
&&+(k_2+k_4)({\rm{tr}}(\nabla n)^2-({\rm{div}}n)^2)\big]+\alpha s^2 |\nabla n|^2 +(\beta-L_2) |\nabla s|^2\nonumber\\
&&+k_5 |\nabla s-(\nabla s\cdot n) n-\nu s(\nabla n) n|^2+k_6|\nabla s\cdot n-\sigma s{\rm{div}} n|^2,
\end{eqnarray}
where
\begin{equation*}
\begin{cases}
\sigma=-\frac{L_3}{2(\beta+L_1)},\ \ \nu=-\frac{L_4}{2(\beta+L_2)},\\
\bar{k}_1=k_1-\sigma^2 k_6 =k_1-\frac{L_3^2}{4(\beta+L_1)},\\
\bar{k}_3=k_3-\nu^2 k_5 =k_3-\frac{L_4^2}{4(\beta+L_2)},\\
k_5=\beta+L_2, \ \ k_6=\beta+L_1.
\end{cases}
\end{equation*}
Based on the physical hypothesis of positivity of the energy density \cite{E2}, it is usually assumed
\begin{equation}\label{positivity}
\bar{k}_1\ge 0, \ k_2\ge 0, \ \bar{k}_3\ge 0, \ k_2\ge |k_4|, \ k_5\ge 0, \ k_6\ge 0, \ \beta>L_2.
\end{equation}
Hence $\mathcal{W}_2(s,n,\nabla s, \nabla n)$ enjoys
the coercivity in $sn$ and of quadratic growth of $\nabla(sn)$:
\begin{equation}\label{coer}
\lambda (|\nabla s|^2+s^2|\nabla n|^2) \le\mathcal{W}_2(s, n,\nabla s, \nabla n)
\le \Lambda (|\nabla s|^2+s^2|\nabla n|^2)
\end{equation}
for two positive constants $\lambda<\Lambda<\infty$ depending only on the coefficients in \eqref{w2-1}.

A sharp interface forms when the size of a liquid crystal droplet is much larger than the ratio of the Frank constants to the surface tension.
In order to study the sharp interface formation between the isotropic phase,
corresponding to $\{s=0\}$, and the nematic
phase, corresponding to $\{s=s_+\}$ for some $s_+\in (0,1)$, we will assume the bulk potential function
takes the form:
$$
\mathcal{W}_0(s)=\frac{1}{\varepsilon^2}W(s),
$$
where $\varepsilon>0$ is a small parameter representing the width of
the interfacial transition region, and the potential function $W\in C^\infty((-\frac12, 1))$
is assumed to be nonnegative, and there exists a unique $s_+\in (0,1)$ such that
\begin{equation}\label{W-cond}
\begin{cases}
W(0)=W(s_+)=0,\ W(s)=W(s_+-s) \ \forall 0\le s\le s_+,\\
\forall \delta_1>0, \exists\delta_2>0 \  {\rm{such\ that}}\ |s|\ge \delta_1 \ {\rm{and}}\ |s-s_+|\ge \delta_1
\Longrightarrow W(s)\ge\delta_2,\\
\displaystyle\lim_{s\rightarrow (-\frac12)^+}{W}(s)=\lim_{s\rightarrow 1^-}{W}(s)=+\infty,
\end{cases}
\end{equation}
In particular, $W$ has two minimal wells of depth zero at $0$ and $s_+$, and
$$W(s)\approx s^2 \ \ {\rm{for}}\ \ |s|<<1, \ \ {\rm{and}}\ \  W(s)\approx (s-s_+)^2
\ \ {\rm{for}}\ \ |s-s_+|<<1.
$$
For $s\not\in [-\frac12, 1]$, we can simply let $W(s)=+\infty$ so that $W$ is defined for all $s\in \R$.

By direct calculations we have the identity
\begin{eqnarray}\label{null-lag1}
{\rm{div}}\big(s^2((\nabla n)n-({\rm{div}} n)n)\big)
=s^2\big({\rm{tr}}(\nabla n)^2-({\rm{div}} n)^2\big)
+2s\nabla s\big((\nabla n)n-({\rm{div}} n)n\big).
\end{eqnarray}
Also recall the null Lagrangian property of
${\rm{div}}\big(s^2((\nabla n)n-({\rm{div}} n)n)\big)$ (see Hardt-Kinderlehrer-Lin \cite{HKL}), that is,
\begin{equation}\label{null-lag2}
\int_\Omega {\rm{div}}\big(s^2((\nabla n)n-({\rm{div}} n)n)\big)\,dx, \ \ \ sn\in H^1(\Omega,\R^3),
\end{equation}
depends only on the value of $(s, n)$ on $\partial\Omega$.

It turns out that both \eqref{null-lag1} and \eqref{null-lag2} will
play a crucial role in our study of phase transitions between the isotropic and the nematic phases.
From now on, we set
\begin{equation}\label{w2-3}
\mathcal{W}_\epsilon(s,n,\nabla s, \nabla n)=\mathcal{W}_2(s,n,\nabla s,\nabla n)
+\frac{1}{\epsilon^2}W(s).
\end{equation}

The problem of sharp interface formations between the isotropic and nematic phases depends
on the relations between the Frank constants $L_1$ and $L_2$. From
\begin{equation}\label{s-decomp}
|\nabla s|^2=|\nabla s\cdot n|^2+|\nabla s\wedge n|^2,
\end{equation}
we see that
$$\beta|\nabla s|^2+L_1(\nabla s\cdot n)^2+L_2|\nabla s\wedge n|^2
=\begin{cases}(\beta+L_2)|\nabla s|^2+(L_1-L_2)(\nabla s\cdot n)^2, & L_1>L_2,\\
(\beta+L_1)|\nabla s|^2+(L_2-L_1) |\nabla s\wedge n|^2, & L_1<L_2,\\
(\beta+L_1)|\nabla s|^2, & L_1=L_2.
\end{cases}
$$
Hence we can reduce the case $L_1>L_2$ into the case ({\bf A});
the case $L_1<L_2$ into the case ({\bf B}); and the case $L_1=L_2$ into
the case ({\bf C}). More precisely, we have
\begin{enumerate}
\item[({\bf A})] \fbox{$L_1>0$ and $L_2=0$} By adding
a null-Lagrangian
$$-\frac12{L_4} {\rm{div}}\big(s^2((\nabla n)n -({\rm{div}}n) n)\big)$$
to $\mathcal{W}_\epsilon$ and applying \eqref{null-lag1} and \eqref{null-lag2},
we can convert $\mathcal{W}_\epsilon$ into
\begin{eqnarray*}
\widetilde{\mathcal{W}}_\epsilon(s,n,\nabla s, \nabla n)
&=&s^2\widetilde{\mathcal{W}}_2(s,n,\nabla s, \nabla n)+\beta|\nabla s|^2+L_1|\nabla s\cdot n|^2\\
&&+(L_3-L_4) (\nabla s\cdot n ) (s {\rm{div}}n) +\frac{1}{\epsilon^2}W(s),
\end{eqnarray*}
where
$$\widetilde{\mathcal{W}}_2(s,n,\nabla s, \nabla n)=\mathcal{W}_2(s,n,\nabla s, \nabla n)
-\frac12{L_4s^2}\big({\rm{tr}}(\nabla n)^2-({\rm{div}} n)^2\big).
$$
Thus, without loss of generality, we will further assume \fbox{$L_4=0$}.

\item[({\bf B})] \fbox{$L_1=0$ and $L_2>0$} By adding
a null-Lagrangian
$$\frac12{L_3} {\rm{div}}\big(s^2((\nabla n)n -({\rm{div}}n) n)\big)$$
to $\mathcal{W}_\epsilon$, we can convert $\mathcal{W}_\epsilon$ into
\begin{eqnarray*}
\widetilde{\mathcal{W}}_\epsilon(s,n,\nabla s, \nabla n)&=&
s^2\widetilde{\mathcal{W}}_2(s,n,\nabla s, \nabla n)+\beta|\nabla s|^2+L_2|\nabla s\wedge n|^2\\
&&+
(L_3+L_4) (s\nabla s) (\nabla n)n+\frac 1{\epsilon^2}W(s),
\end{eqnarray*}
where
$$\widetilde{\mathcal{W}}_2(s,n,\nabla s, \nabla n)=\mathcal{W}_2(s,n,\nabla s, \nabla n)+\frac12{L_3s^2}
\big({\rm{tr}}(\nabla n)^2-({\rm{div}} n)^2\big).
$$
Thus, without loss of generality, we will further assume \fbox{$L_3=0$}.

\item[({\bf C})] \fbox{$L_1=L_2=0$} In this case, we will only consider the following two subcases:
\begin{enumerate}
\item [({\bf C}1)] \fbox{$L_3=L_4=0$}. Hence $\mathcal{W}_2$ can be rewritten as
\begin{eqnarray} \label{w21}
\widetilde{\mathcal{W}}_2(s,n,\nabla s,\nabla n)&=&
s^2\big[k_1({\rm{div}} n)^2+k_2(n\cdot {\rm{curl}} n)^2+k_3 |n\wedge {\rm{curl}} n|^2\nonumber\\
&+&(k_2+k_4)({\rm{tr}}(\nabla n)^2-({\rm{div}} n)^2)\big]+\alpha s^2|\nabla n|^2+\beta|\nabla s|^2.
\end{eqnarray}
\item [({\bf C}2)] \fbox{$L_4=-L_3$}. Hence after adding the null Lagrangian
$$L_3 {\rm{div}}\big(s^2(({\rm{div}} n) n-(\nabla n) n)\big)$$
to $\mathcal{W}_2$, we can convert $\mathcal{W}_2$ to $\widetilde{\mathcal{W}}_2$
\begin{eqnarray} \label{w22}
\widetilde{\mathcal{W}}_2(s,n,\nabla s,\nabla n)
&=&s^2\big[k_1({\rm{div}} n)^2+k_2(n\cdot {\rm{curl}} n)^2+k_3 |n\wedge {\rm{curl}} n|^2\\
&+&(k_2+k_4+L_3)({\rm{tr}}(\nabla n)^2-({\rm{div}} n)^2\big]+\alpha s^2|\nabla n|^2
+\beta|\nabla s|^2.\nonumber
\end{eqnarray}
Thus, after replacing $(k_2+k_4+L_3)$ by $(k_2+k_4)$, ({\bf C}2) can be reduced to ({\bf C}1).
\end{enumerate}
\end{enumerate}

Through this paper, we denote by $\mathcal{H}^2$ the two dimensional Hausdorff measure in $\mathbb R^3$.
Define
the $1$-dimensional minimal connecting energy by
\begin{equation}\label{1d-energy}
{\alpha}_0=\inf\Big\{\int_{-\infty}^\infty \big(\beta \dot s(t)^2+W(s(t))\big)\,dt
\ \big|\ s\in C^{0,1}((-\infty,\infty), \R),
\ s(-\infty)=0, \ s(\infty)=s_+\Big\}.
\end{equation}
It is well-known that $\alpha_0$ is attained by an $\xi\in C^\infty((-\infty,\infty), \R)$,
which satisfies
\begin{equation}\label{1d_orbit1}
\sqrt{\beta}\xi'(t)=\sqrt{W(\xi(t))} \ \ {\rm{in}}\ \ (-\infty,\infty); \ \ \xi(-\infty)=0, \ \xi(\infty)=s_+,
\end{equation}
and
%\begin{equation}\label{1d_orbit2}
%\sqrt{\beta}\xi'(t)=\sqrt{W(\xi(t))}, \ t\in (-1,1).
%\end{equation}
%In particular,
\begin{eqnarray}\label{1d-energy1}
{\alpha}_0&=& 2\sqrt{\beta}\int_0^{s_+} \sqrt{W(t)}\,dt\nonumber\\
&=&\inf\Big\{2\int_{-\infty}^\infty \sqrt{\beta W(s(t))}|\dot{s}(t)|\,dt
\ \big|\ s\in C^{0,1}((-\infty,\infty), \R),
\ s(-\infty)=0, \ s(\infty)=s_+\Big\}\nonumber\\
&=&\inf\Big\{2\int_{-a}^a \sqrt{\beta W(s(t))}|\dot{s}(t)|\,dt
\ \big|\ s\in C^{0,1}([-a,a], \R),
\ s(-a)=0, \ s(a)=s_+\Big\},
\end{eqnarray}
for any $a>0$.

Now we state our first theorem. It concerns the $\Gamma$-convergence of minimizers of the Ericksen energy functional
\begin{equation}\label{w2-4}
\mathcal{E}(s_\epsilon, n_\epsilon):=\int_\Omega\widetilde{\mathcal{W}}_\epsilon(s,n,\nabla s,\nabla n)\,dx,
\ \ {\rm{with}}\ \ \widetilde{\mathcal{W}}_\epsilon(s,n,\nabla s,\nabla n)=\widetilde{\mathcal{W}}_2(s,n,\nabla s,\nabla n)
+\frac{1}{\epsilon^2}W(s),
\end{equation}
either under well prepared Dirichlet boundary values $(t_\epsilon, g_\epsilon)$ when $\Omega\subset\R^3$ is a bounded smooth domain, or under the volume constraint for nematic region when $\Omega=\R^3$, as $\epsilon\rightarrow 0$. Notice that for any fixed $\epsilon>0$,
the existence and regularity of minimizer $(s_\epsilon, n_\epsilon)$ to \eqref{w2-4} have been studied by Lin \cite{Lin0, Lin},
Lin-Poon \cite{LP} and Ambrosio \cite{Ambrosio1, Ambrosio2}.

For a bounded smooth $\Omega\subset\R^3$, we prescribe $(t_\epsilon, g_\epsilon):\partial\Omega\to \R\times\mathbb{S}^2$ as follows. Let $\Sigma^\pm\subset\Omega$ be two disjoint,
connected open subset of $\partial\Omega$ such that
\begin{itemize}
\item [i)]
$\partial\Sigma^\pm=\Sigma^0$ is a smooth, closed curve of $\partial\Omega$, and
$\partial\Omega=\Sigma^+\cup\Sigma^-\cup\Sigma^0$.
% and $\Sigma^0$ spans an area minimizing surface $\Gamma\subset\Omega$.
\item [ii)]
there exists $L>0$ such that
$t_\epsilon\in H^1\big(\partial\Omega\big)$
satisfies
%$s_\epsilon:\partial\Omega\to \R$ satisfies
\begin{equation}\label{b_cond1.1}
\big\|t_\epsilon\big\|_{L^2(\widehat{\Sigma}^{-})}\rightarrow 0
 \ \ {\rm{and}} \ \
 \big\|t_\epsilon-s_+\big\|_{L^2(\widehat{\Sigma}^+)}\rightarrow 0,
 \ \ {\rm{as}}\ \ \epsilon\rightarrow 0,
\end{equation}
\begin{equation}\label{b_cond1.2}
\sup_{0<\epsilon<1}\int_{\partial\Omega} \big(\epsilon (|\nabla_{\rm{tan}} t_\epsilon|^2
+\frac{1}{\epsilon}W(t_\epsilon)\big)\,d\mathcal{H}^2 \le L,
\end{equation}
and there exists a map  $\widehat{n}_\epsilon\in H^1(\Omega,\mathbb S^2)$ such that
$\widehat{n}_\epsilon=g_\epsilon$ on $\partial\Omega$, and
\begin{equation}\label{b_cond1.20}
\begin{split}
\begin{cases}
\widehat{n}_\epsilon\cdot\nu_\Gamma=0 \ {\rm{on}}\ \Gamma & \mbox{under the condition ({\bf A})}: L_1>L_2=L_4=0,\\
\widehat{n}_\epsilon\wedge\nu_{\Gamma}=0 \ {\rm{on}}\ \Gamma & \mbox{under the condition ({\bf B})}: L_2>L_1=L_3=0,
\end{cases}\\
\lim_{\epsilon\rightarrow 0}\epsilon\int_\Omega |\nabla \widehat{n}_\epsilon|^2\,dx=0,
\end{split}
\end{equation}
where $\Gamma\subset\Omega$ is an area minimizing surface such that $\partial\Gamma=\Sigma^0$.
\end{itemize}

%where $\Sigma^\pm_\epsilon=\big\{x\in\Sigma^\pm: \ {\rm{d}}(x,\Sigma)>\epsilon\big\}$,

\begin{theorem}\label{existence} Assume either the condition $({\bf A})$ $L_1>L_2=L_4=0$
and $3L_3^2\le 4 L_1\alpha$, or the condition $({\bf B})$ $L_2>L_1=L_3=0$ and $L_4^2\le 4L_2\alpha$,
or the condition $({\bf C})$ $L_1=L_2=L_3=L_4=0$ holds. Then we have
\begin{itemize}
\item [(i)] If $\Omega\subset\R^3$ is a bounded smooth domain and $(t_\epsilon, g_\epsilon):\partial\Omega\to\R\times\mathbb S^2$
satisfies \eqref{b_cond1.1}, \eqref{b_cond1.2} and \eqref{b_cond1.20}, then
\begin{eqnarray}\label{minimality0}
&&\lim_{\epsilon\rightarrow 0}\ \inf\Big\{ \int_\Omega\epsilon \widetilde{\mathcal{W}}_\epsilon(s,n,\nabla s,\nabla n)\,dx\ \big|\ (s,n):\Omega\to\R\times \mathbb S^2, \nonumber\\
&&\qquad\qquad\ sn\in H^1(\Omega, \R^3), \ (s,n)\big|_{\partial\Omega}= (t_\epsilon, g_\epsilon)\Big\}=\alpha_0\mathcal{H}^2(\Gamma),
\end{eqnarray}
where $\Gamma\subset\Omega$ is an area minimizing surface with $\partial\Gamma=\Sigma^0$.
\item [(ii)] If $\Omega=\mathbb R^3$,
then
\begin{eqnarray}\label{minimality00}
&&\lim_{\epsilon\rightarrow 0}\ \inf\Big\{ \int_{\R^3} \epsilon\widetilde{\mathcal{W}}_\epsilon(s,n,\nabla s,\nabla n)\,dx
\ \big|\ (s,n): \R^3\to\R\times \mathbb S^2,\ sn\in H^1(\R^3, \R^3), \nonumber\\
&&\qquad\qquad \big|\{x\in\R^3: s(x)\ge {s_+}\}\big|=|B_1|\Big\}=\alpha_0\mathcal{H}^2(\partial B_1),
\end{eqnarray}
where $B_1$ is a ball of radius $1$.\footnote{the volume constraint $|\{x\in\R^3, s(x)>{s_+}\}|=|B_1|$ can be replaced by $|\{x\in\R^3, s(x)>{s_+}\}|=\lambda$ for any
$\lambda>0$. For convenience, we choose $\lambda=|B_1|$.}
\end{itemize}
\end{theorem}

Theorem \ref{existence} can be proved in the framework of $\Gamma$-convergence:
\begin{itemize}
\item [(1)] First, under the conditions on the coefficients $L_i$'s and
$\alpha$, we can show the energy is bounded below by
$$\int_\Omega \big(\epsilon\beta|\nabla s_\epsilon|^2+\frac{1}{\epsilon}W(s_\epsilon)\big)\,dx,$$
which becomes a scalar-valued Allen-Cahn functional so the technique in the BV function space,
as in \cite{Modica-Mortola}, or the isoperimetric inequality in $\R^3$
can be employed to show it is bounded by $\alpha_0\mathcal{H}^2(\Gamma)$.
\item [(2)] Secondly, we construct a comparison map $(s_\epsilon, n_\epsilon)$ by letting
$n_\epsilon=\widehat{n}_\epsilon$
and by placing an almost optimal $1$-dimensional orbit $s_\epsilon(x)=\xi_L(\frac{d_\Gamma(x)}{\epsilon})$ in the transversal direction to $\Gamma$ within $L\epsilon$-width (with $L>>1$), and away from this region, $s_\epsilon$ is made to have very small spatial variations. It turns out that the contribution of anchoring energy
$\displaystyle\int_{\Gamma_{L\epsilon}}|\nabla s_\epsilon\cdot\widehat{n}_\epsilon|^2\,dx$ or
$\displaystyle\int_{\Gamma_{L\epsilon}}|\nabla s_\epsilon\wedge\widehat{n}_\epsilon|^2\,dx$ can be made
arbitrarily small.
\end{itemize}

It follows from Theorem \ref{existence} that the leading order term of $\mathcal{E}(s_\epsilon, n_\epsilon)$ for
minimizers $(s_\epsilon, n_\epsilon)$ is
$\displaystyle\frac{\alpha_0}{\epsilon}\mathcal{H}^2(\Gamma)$, so that
$$\mathcal{E}(s_\epsilon, n_\epsilon)=\displaystyle\frac{\alpha_0}{\epsilon}\mathcal{H}^2(\Gamma)+\mathcal{D}(\epsilon).$$

An important question is to ask for the asymptotic behavior of $\mathcal{D}(\epsilon)$, which is our focus
in this work. For this purpose, we will need to assume that
the boundary value $(t_\epsilon, g_\epsilon)$ provides an almost optimal transition in the fast transition area on $\partial\Omega$
across the interfacial curve $\Sigma^0$.

\medskip
To describe our results, we need to introduce some notations. Set the Oseen-Frank energy density for $n$, with $|n|=1$, by
\begin{equation}\label{OF}
W_{OF}(n, \nabla n)=k_1({\rm{div}} n)^2+k_2(n\cdot {\rm{curl}} n)^2+k_3 |n\wedge {\rm{curl}} n|^2
+(k_2+k_4)\big({\rm{tr}}(\nabla n)^2-({\rm{div}} n)^2\big).
\end{equation}
For $\delta>0$, define the $\delta$-neighborhood of $\partial\Omega$ and $\Gamma$ by
$$(\partial\Omega)_\delta=\big\{x\in\overline\Omega: \ {\rm{d}}(x,\partial\Omega)<\delta\big\},
$$
and
$$\Gamma_{\delta}=\big\{x\in\overline\Omega: \ {\rm{d}}(x,\Gamma)<\delta\big\}=\cup_{-\delta<\lambda<\delta}
\Gamma(\lambda),$$
where $\Gamma(\lambda)=\big\{x\in\overline\Omega: \ d_{\Gamma}(x)=\lambda\big\}$, and
$d_{\Gamma}(x)$ is the signed distance function
of $x$ to $\Gamma$:
$$
d_{\Gamma}(x)=\begin{cases} -d(x,\Gamma) & \ x\in\Omega^-,\\
\ \ d(x,\Gamma) & \ x\in\Omega^+.
\end{cases}
$$

Let $\Omega^\pm\subset\Omega$ be the connected components
such that $\partial\Omega^\pm=\Sigma^\pm\cup\Gamma$.
Set
$$\Gamma_\delta^\pm=\Gamma_\delta\cap\Omega^\pm, \ \ U^\pm_\delta=\Omega^\pm\setminus\Gamma_\delta,
\ \ \Sigma^\pm_\delta=\Sigma^\pm\setminus\Gamma_\delta,
\ \ Q_\delta^\pm=U_\delta^\pm\cap (\partial\Omega)_\delta,
$$
$$
\Omega_\delta=\Omega\setminus (\partial\Omega)_\delta,
\ \ \Omega_\delta^\pm=\Omega_\delta\cap\Omega^\pm, \ \
V_\delta^\pm=\Omega_\delta^\pm\setminus\Gamma_\delta,\ \
W_\delta^\pm=\Omega^\pm_\delta\cap \Gamma_\delta, \ \ O_\delta=\Gamma_\delta
\setminus (W_\delta^+\cup W_\delta^-).
$$

It follows from the condition \eqref{W-cond} and Proposition A.4 of \cite{LPW} that
for any $\gamma\in (\frac12,1)$, there exist an almost minimal orbit
$\xi_{\epsilon,\gamma}\in C^\infty([-\epsilon^{\gamma-1}, \epsilon^{\gamma-1}],\R)$ and
$C_1, C_2, C_3>0$ independent of $\epsilon$ such that
\begin{equation}\label{xi_cond}
\begin{cases}
\xi_{\epsilon, \gamma}(-\epsilon^{\gamma-1})=0,\ \xi_{\epsilon, \gamma}(\epsilon^{\gamma-1})=s_+,\\
\xi_{\epsilon, \gamma} (t) =s_+-\xi_{\epsilon, \gamma} (-t), \ \forall t\in (-\epsilon^{\gamma-1}, \epsilon^{\gamma-1}),\\
\displaystyle\max_{|t|\le\epsilon^{\gamma-1}}\Big|\beta|{\xi_{\epsilon,\gamma}}'|^2+W(\xi_{\epsilon,\gamma})-2\sqrt{\beta W(\xi_{\epsilon,\gamma})}|\xi'_{\epsilon,\gamma}|\Big|\le C_2e^{-C_1\epsilon^{\gamma-1}},\\
\displaystyle\int_{-\epsilon^{\gamma-1}}^{\epsilon^{\gamma-1}}\big(\beta|{\xi_{\epsilon,\gamma}}'|^2+W(\xi_{\epsilon,\gamma})\big)\,d\tau\le
\alpha_0+C_2e^{-C_1\epsilon^{\gamma-1}},\\
\big|{\xi_{\epsilon,\gamma}}'(t)\big|\le C_2e^{-C_1t}, \ \forall |t|\ge C_3.
\end{cases}
\end{equation}
We assume that
\begin{equation}\label{b_cond1.3}
\begin{cases}
\displaystyle t_\epsilon(x)=\xi_{\epsilon,\gamma}(\frac{d_\Gamma(x)}{\epsilon}), \ x\in\Sigma^+_{\epsilon^{\gamma}}\cup \Sigma^{-}_{\epsilon^{\gamma}}
=\big\{x\in\partial\Omega: -\epsilon^{\gamma}\le d_\Gamma(x)\le \epsilon^{\gamma}\big\},\\
\displaystyle\big\|t_\epsilon-s_+\big\|_{L^2(\Sigma^+\setminus \Gamma_{\epsilon^\gamma})}=O(\epsilon)
\ \ {\rm{and}}\ \ \big\|t_\epsilon\big\|_{L^2(\Sigma^-\setminus \Gamma_{\epsilon^\gamma})}=O(\epsilon),\\
\displaystyle\lim_{\epsilon\to 0}\epsilon^\gamma\int_{\partial\Omega\setminus\Gamma_{\epsilon^\gamma}}|\nabla_{\rm{tan}}t_\epsilon|^2\,d\mathcal{H}^2=0.
\end{cases}
\end{equation}
%Here $\xi\in C^\infty([-c, c], \R)$, with $\xi(-c)=0$ and $\xi(c)=s_+$, is a minimizer of $\alpha_0$.

To simplify the technical presentation,
we further assume that there exists an $g\in H^1(\partial\Omega)$ such that $g_\epsilon\rightarrow g$
in $H^1(\partial\Omega)$, and
\begin{equation}\label{b_cond1.4}
g_\epsilon=g\ \ {\rm{on}}\ \ \partial\Omega\cap\Gamma_{\epsilon^\gamma},\ \ {\rm{and}}\ \
\big\|g_\epsilon-g\big\|_{L^2(\partial\Omega\setminus\Gamma_{\epsilon^\gamma})}=o_{\epsilon}(1)\epsilon^\gamma.
\end{equation}

%Hence the condition \eqref{b_cond1.2} reduces to
%\begin{equation}\label{b_cond1.5}
%\int_{\partial\Omega} \big(|\nabla_{\rm{tan}} s_\epsilon|^2+\frac{1}{\epsilon^2} W(s_\epsilon)\big)\,d\mathcal{H}^2\le \frac{L}{\epsilon}, \ \forall 0<\epsilon<1.
%\end{equation}

Recall that a minimal surface $S$ is called strictly stable, if, in addition,
\begin{equation}\label{stable}
\frac{d^2}{dt^2}\big|_{t=0} \mathcal{H}^2\big(\{x+t\phi(x)\nu_S(x): \ x\in S\}\big)>0,
\ \forall 0\not\equiv\phi\in C^\infty_0(S),
\end{equation}
where $\nu_S$ is a unit normal vector field of $S$.

The main contributions of our paper concern the characterization of the $O(1)$-term, $\mathcal{D}(\epsilon)$,
in the energy expansion of $\mathcal{E}(s_\epsilon, n_\epsilon)$. We divide our results into two separate theorems.
The first one deals the case that $\Omega$ is a bounded domain in $\R^3$.

\begin{theorem} \label{sharp} Let $\Omega\subset\mathbb R^3$
be a bounded smooth domain. Assume that $\Gamma$ is a unique, strictly stable, area minimal surface
spanned by $\Sigma^0$, and the boundary values $(t_\epsilon, g_\epsilon)$ satisfy conditions
\eqref{b_cond1.2}, \eqref{xi_cond}, \eqref{b_cond1.3}, and \eqref{b_cond1.4}.
Let $(s_\epsilon, n_\epsilon)$, $0<\epsilon<1$, be minimizers of
$\int_\Omega \widetilde{\mathcal{W}}_\epsilon(s,n,\nabla s, \nabla n)\,dx$, subject to  the boundary condition $(s_\epsilon, n_\epsilon)=(t_\epsilon, g_\epsilon)$ on ${\partial\Omega}$. Then the following statements
hold:
\begin{enumerate}
\item[{\rm{(}}A{\rm{)}}] If $L_1>0$, $L_2=L_4=0$, and
\begin{equation}\label{cond1}
3L_3^2< 4L_1\alpha,
%\ \big( or\  |L_3(\nabla s\cdot n)(s{\rm{div}}n)|\le \frac12{L_1}|\nabla s\cdot n|^2+\frac12{\alpha} s^2|\nabla n|^2\big),
\end{equation}
then
\begin{equation}\label{asymp1}
\mathcal{E}(s_\epsilon, n_\epsilon)
=\frac{\alpha_0}{\epsilon}\mathcal{H}^2(\Gamma)+\mathcal{D}_A+o_\epsilon(1).
\end{equation}
Here $\mathcal{D}_A$ is given by
\begin{equation}\label{OFM1}
\mathcal{D}_A=\inf\Big\{E(n; \Omega^+)=s_+^2\int_{\Omega^+} \big(W_{OF}(n,\nabla n)+\alpha|\nabla n|^2\big)\,dx\Big\}
\end{equation}
among all maps $n\in H^1(\Omega^+,\mathbb S^2)$ satisfying the planar anchoring condition on $\Gamma$:
\begin{equation}\label{bdry1}
n=g \ \ {\rm{on}}\ \ \Sigma^+;
\ \ n\cdot\nu_\Gamma=0 \ \ {\rm{on}}\ \ \Gamma,
\end{equation}
where $\nu_\Gamma$ is the outward unit normal of the nematic region
$\Omega^+\equiv\big\{x\in\Omega: \ s(x)=s_+\big\}$.
\item[{\rm{(}}B{\rm{)}}] If $L_2>0$, $L_1=L_3=0$, and
\begin{equation}\label{cond2}
L_4^2< 4L_2\alpha,
%\ (or\ |L_4(s\nabla s)(\nabla n)n |\le \frac12{L_2}|\nabla s\wedge n|^2+\frac12{\alpha} s^2|\nabla n|^2),
\end{equation}
then
\begin{equation}\label{asymp2}
\mathcal{E}(s_\epsilon, n_\epsilon)
=\frac{\alpha_0}{\epsilon}\mathcal{H}^2(\Gamma)+\mathcal{D}_B+o_\epsilon(1).
\end{equation}
Here $D_B$ is given by
\begin{equation}\label{OFM2}
\mathcal{D}_B=\inf\Big\{E(n;\Omega^+)=s_+^2\int_{\Omega^+} \big(W_{OF}(n,\nabla n)+\alpha|\nabla n|^2\big)\,dx\Big\}
\end{equation}
among all maps $n\in H^1(\Omega^+,\mathbb S^2)$ satisfying
the homeotropic anchoring condition on $\Gamma$:
\begin{equation}\label{bdry2}
n=g \ \ {\rm{on}}\ \ \Sigma^+; \ \ \
n\wedge\nu_\Gamma=0  \ \ {\rm{on}}\ \ \Gamma.
\end{equation}
\item [{\rm{(}}C{\rm{)}}] If $L_1=L_2=L_3=L_4=0$,
then
\begin{equation}\label{asymp3}
\mathcal{E}(s_\epsilon, n_\epsilon)
=\frac{\alpha_0}{\epsilon}\mathcal{H}^2(\Gamma)+\mathcal{D}_C+o_\epsilon(1).
\end{equation}
Here $\mathcal{D}_C$ is given by
\begin{equation}\label{OFM3}
\mathcal{D}_C=\inf\Big\{E(n; \Omega^+)=s_+^2\int_{\Omega^+} \big(W_{OF}(n,\nabla n)+\alpha|\nabla n|^2\big)\,dx\Big\}
\end{equation}
among all maps $n\in H^1(\Omega^+,\mathbb S^2)$ satisfying the free boundary condition on $\Gamma$:
\begin{equation}\label{bdry3}
 n=g \ \ \ {\rm{on}}\ \ \ \Sigma^+.
\end{equation}
\end{enumerate}
\end{theorem}

We would like to remark that the interior regularity and boundary regularity near $\Sigma^+$
of minimizing harmonic maps $n\in H^1(\Omega^+,\mathbb S^2)$ achieving $\mathcal{D}_A$,
or $\mathcal{D}_B$, or $\mathcal{D}_C$ has been studied by Hardt-Kinderlehrer-Lin \cite{HKL}.
For the boundary regularity of $n$ near the interface $\Gamma$ when the isotropic Oseen-Frank energy is considered,
we refer to Hardt-Lin \cite{HL2} and Duzaar-Steffen \cite{DS1,DS2}
for partially constrained or free boundary conditions, and Day-Zarnescu \cite{DZ} under the planar anchoring condition.

\medskip
The second one considers the entire space $\Omega=\R^3$.

\begin{theorem} \label{sharp1} Let $(s_\epsilon, n_\epsilon):\R^3\to\R\times\mathbb{S}^2$, $0<\epsilon<1$, be minimizers of
$\int_{\R^3} \widetilde{\mathcal{W}}_\epsilon(s,n,\nabla s, \nabla n)\,dx$, subject to  the
constraint:
$$\big|\{x\in\R^3:\ s_\epsilon\ge {s_+}\}\big|=|B_1|.$$
Then the following statements hold:
\begin{enumerate}
\item[{\rm{(}}A1{\rm{)}}] If $L_1>0$, $L_2=L_4=0$, and
\begin{equation}\label{cond10}
3L_3^2< 4L_1\alpha,
%\ \big( or\  |L_3(\nabla s\cdot n)(s{\rm{div}}n)|\le \frac12{L_1}|\nabla s\cdot n|^2+\frac12{\alpha} s^2|\nabla n|^2\big),
\end{equation}
then
\begin{equation}\label{asymp10}
\mathcal{E}(s_\epsilon, n_\epsilon)
=\frac{\alpha_0}{\epsilon}\mathcal{H}^2(\partial B_1)+\mathcal{D}_A+o_\epsilon(1),
\end{equation}
where
 $\mathcal{D}_A$ is given by
\begin{equation}\label{OFM10}
\mathcal{D}_A=\inf\Big\{E(n; B_1)=s_+^2\int_{B_1} \big(W_{OF}(n,\nabla n)+\alpha|\nabla n|^2\big)\,dx
\ \big| \ n\in H^1\big(B_1, \mathbb S^2\big)\Big\},
\end{equation}
subject to the planar anchoring condition:
\begin{equation}\label{bdry10}
n(x)\cdot x=0  \ \ {\rm{on}}\ \ \partial B_1.
\end{equation}
\item[{\rm{(}}B1{\rm{)}}] If $L_2>0$, $L_1=L_3=0$, and
\begin{equation}\label{cond20}
L_4^2< 4L_2\alpha,
%\ (or\ |L_4(s\nabla s)(\nabla n)n |\le \frac12{L_2}|\nabla s\wedge n|^2+\frac12{\alpha} s^2|\nabla n|^2),
\end{equation}
then
\begin{equation}\label{asymp20}
\mathcal{E}(s_\epsilon, n_\epsilon)
=\frac{\alpha_0}{\epsilon}\mathcal{H}^2(\partial B_1)+\mathcal{D}_B+o_\epsilon(1),
\end{equation}
where  $D_B$ is given by
\begin{equation}\label{OFM20}
\mathcal{D}_B=\inf\Big\{E(n; B_1)=s_+^2\int_{B_1} \big(W_{OF}(n,\nabla n)+\alpha|\nabla n|^2\big)\,dx
\ \big|
\ n\in H^1(B_1,\mathbb S^2)\Big\},
\end{equation}
subject to the homeotropic anchoring condition:
\begin{equation}\label{bdry20}
n(x)\wedge x=0  \ \ {\rm{on}}\ \ \partial B_1.
\end{equation}
\item [{\rm{(}}C1{\rm{)}}] If $L_1=L_2=L_3=L_4=0$,
then
\begin{equation}\label{asymp30}
\mathcal{E}(s_\epsilon, n_\epsilon)
=\frac{\alpha_0}{\epsilon}\mathcal{H}^2(\partial B_1)+o_\epsilon(1).
\end{equation}
\end{enumerate}
\end{theorem}

We would like to point out that for a bounded domain $\Omega\subset\mathbb R^n$, while the boundary conditions imposed on $(t_\epsilon, g_\epsilon)$
in Theorems 1.1, 1.2, and 1.3 are physically natural,  mathematically they are rather technical to describe. On the other hand, if we consider the
same type problems on a compact manifold $M$ without boundary or a torus $\mathbb T^n$, then the natural condition would be the
volume constraint on approximate nematic regions $\{s_\epsilon\ge s_+\}$. Hence the problem can be significantly simplified because we will
have the compactness of the space (in contrast with $\mathbb R^n$) and avoid the technical issues arising from both 
the physical boundary and the boundary values.

We would like to remark that the regularity of minimizing harmonic maps in the case (B1) was studied
by \cite{HKL} and \cite{HL1, HL2}. See \cite{DZ} for some work related to the boundary regularity of minimizing harmonic maps
in the case (A1).

While the approach to prove Theorem \ref{sharp} and Theorem \ref{sharp1}
is based on the technique of $\Gamma$-convergence,
it is very delicate to obtain the exact characterization of $O(1)$-term in the expansion of $\mathcal{E}(s_\epsilon, n_\epsilon)$ especially when we deal with a bounded domain $\Omega$ with physical boundary data.
\begin{itemize}
\item [(1)] For the construction of sharp upper bound, we need to place an almost minimal $1$-dimensional
orbit $\xi_{\epsilon, \gamma}(\frac{d_\Gamma(x)}{\epsilon})$ in the transversal direction of $\Gamma$
within the width of $O(\epsilon^\gamma)$, which guarantees the GL energy is of $\alpha_0+o_\epsilon(1)\epsilon$,
while we have to utilize the decay property of $\xi_{\epsilon,\gamma}'(t)$ ensuring
$\displaystyle\int_{\Gamma_{\epsilon^\gamma}}|\nabla s_\epsilon\cdot n_\epsilon|^2\,dx$
(or $\displaystyle\int_{\Gamma_{\epsilon^\gamma}}|\nabla s_\epsilon\wedge n_\epsilon|^2\,dx$) is of order $o_\epsilon(1)$.
\item [(2)] To achieve a sharp lower bound, we need to extract a sequence of sets of finite perimeters
$E_\epsilon=\big\{x\in\Omega: s_\epsilon(x)\ge \delta_\epsilon\big\}$ with uniformly bounded perimeters
such that $\mathcal{H}^2(\partial^*E_{\epsilon}\lfloor\Omega)\approx \mathcal{H}^2(\Gamma)$,
$\displaystyle\int_{E_\epsilon}|\nabla n_\epsilon|^2\,dx\le C$, and
$\displaystyle\int_{\partial^*E_\epsilon\lfloor\Omega} \big|\frac{\nabla s_\epsilon}{|\nabla s_\epsilon|}\cdot n_\epsilon\big|^2\,d\mathcal{H}^2 \ (or \displaystyle\int_{\partial^*E_\epsilon\lfloor\Omega} \big|\frac{\nabla s_\epsilon}{|\nabla s_\epsilon|}\wedge n_\epsilon\big|^2\,d\mathcal{H}^2) \le C$. Then we adapt the techniques from \cite{LPW} and
some measure theoretic arguments to show that $E_{\epsilon}\rightarrow \Omega^+$, and
$n_\epsilon\rightarrow n$ in SBV($\Omega$) for some $n\in H^1(\Omega^+,\mathbb S^2)$ with
$n\cdot\nu_\Gamma=0$ (or $n\wedge\nu_\Gamma=0$).
\item [(3)] Utilize the strict stability of $\Gamma$ to show that the leading order coefficients of $\frac{1}{\epsilon}$
in both lower and upper bound estimates match up to order $o_\epsilon(1)\epsilon$.
\end{itemize}

When dealing with the entire space $\Omega=\R^3$, we observe that the approximate nematic
region $E_\epsilon$ constitutes a minimizing sequence of sets that approach the isoperimetric
inequality so that we can apply the quantitative stability theorem by Fusco-Maggi-Pratelli \cite{FMP} (see also
Maggi \cite{Maggi1}) to show, after suitable translations, $E_\epsilon$ converges to $B_1$ in $L^1$.

Theorem 1.3 is also related to the optimal shape problem of variational problems on liquid crystal droplets previously studied by Lin-Poon \cite{LP1}.
More precisely, Lin and Poon \cite{LP1} considered  the following minimization problem
\begin{equation}\label{LD}
\inf\big\{\int_\Omega |\nabla n|^2\,dx+\mu\mathcal{H}^{n-1}(\partial\Omega)\ | \ n\in H^1(\Omega, \mathbb S^2), n(x)=\nu_{\partial\Omega}, \ |\Omega|=|B_1|\big\}.
\end{equation}
Among the class of convex domains $\Omega$, it was shown by \cite{LP1} that $(\Omega, u)=(B_1, \frac{x}{|x|})$ is a unique minimizer of  \eqref{LD}.
Very recently, this result was extended by Li-Wang \cite{LiWang} to the class of star-shaped mean convex domains  in $\mathbb R^3$.

The paper is organized as follows. In section 2, we will establish both lower and upper bounds of
$\epsilon \mathcal{E}_\epsilon(s_\epsilon, n_\epsilon)$ and prove Theorem \ref{existence}.
In section 3, we will study the bounded domain case and
establish both refined lower and upper bounds for $\mathcal{E}_\epsilon(s_\epsilon, n_\epsilon)$
for all three cases and then prove Theorem \ref{sharp}. In section 4, we will study the case
that $\Omega$ entire space $\mathbb R^3$ and prove Theorem \ref{sharp1}.

\section{Proof of Theorem \ref{existence}}
\setcounter{equation}{0}
\setcounter{theorem}{1}

In this section, we will provide a proof of Theorem \ref{existence}. It involves (a) a concrete
construction of comparison map $(s_\epsilon, n_\epsilon)$ in which $s_\epsilon$ exhibits
a fast transition near $\Gamma$ with energy order $\frac{\alpha_0}{\epsilon}\mathcal{H}^2(\Gamma)$;
and (b) obtain the lower bound  by typical arguments of singular perturbations
of functions of bounded variations.

\subsection{Lower bound estimates}

For either a bounded $\Omega\subset\R^3$ or $\Omega=\R^3$ itself, we assume that
\begin{equation}\label{LB1}
\Lambda=\liminf_{\epsilon\to 0}\epsilon\mathcal{E}(s_\epsilon, n_\epsilon)=
\liminf_{\epsilon\to 0} \int_\Omega \epsilon\widetilde{\mathcal{W}}_\epsilon(s_\epsilon,
n_\epsilon, \nabla s_\epsilon, \nabla n_\epsilon)\,dx<\infty.
\end{equation}
It follows from the condition \eqref{positivity} that
$$W_{OF}(n_\epsilon, \nabla n_\epsilon)\ge 0.$$
Observe that by Cauchy-Schwarz inequality, the following properties hold:
\begin{itemize}
\item[(i)] If $3L_3^2\le 4 L_1\alpha$, then
\begin{equation}\label{CS1}
|L_3(\nabla s_\epsilon\cdot n_\epsilon)s_\epsilon {\rm{div}} n_\epsilon|
\le \sqrt{3}|L_3| |\nabla s_\epsilon\cdot n_\epsilon| |s_\epsilon| |\nabla n_\epsilon|
\le L_1 (\nabla s_\epsilon\cdot n_\epsilon)^2+\alpha s_\epsilon^2 |\nabla n_\epsilon|^2,
\end{equation}
where we have used the inequality $|{\rm{div}} n_\epsilon|\le \sqrt{3}|\nabla n_\epsilon|$.
\item[(ii)] If $L_4^2\le 4L_2\alpha$, then
\begin{eqnarray}\label{CS2}
|L_4 s_\epsilon\nabla s_\epsilon (\nabla n_\epsilon) n_\epsilon|
&=& |L_4| \big|s_\epsilon\big(\nabla s_\epsilon-(\nabla s_\epsilon\cdot n_\epsilon)n_\epsilon\big) (
\nabla n_\epsilon) n_\epsilon\big|\nonumber\\
&\le& L_2 \big|\nabla s_\epsilon-(\nabla s_\epsilon\cdot n_\epsilon)n_\epsilon\big|^2
+\alpha s_\epsilon^2|(\nabla n_\epsilon) n_\epsilon|^2\nonumber\\
&\le & L_2|\nabla s_\epsilon\wedge n_\epsilon|^2+\alpha s_\epsilon^2 |\nabla n_\epsilon|^2,
\end{eqnarray}
where we have used the fact that $n_\epsilon\cdot ((\nabla n_\epsilon)n_\epsilon)=0$, \eqref{s-decomp},
and $|(\nabla n_\epsilon)n_\epsilon|\le |\nabla n_\epsilon|$.
\end{itemize}
Since
\begin{eqnarray*}
&&\widetilde{\mathcal{W}}_\epsilon(s_\epsilon, n_\epsilon, \nabla s_\epsilon,\nabla n_\epsilon)
=s_\epsilon^2 W_{OF}(n_\epsilon, \nabla n_\epsilon)
+\alpha s_\epsilon^2|\nabla n_\epsilon|^2+\beta |\nabla s_\epsilon|^2
+L_1 (\nabla s_\epsilon\cdot n_\epsilon)^2\\
&&\qquad\qquad+L_2|\nabla s_\epsilon\wedge n_\epsilon|^2
+L_3(\nabla s_\epsilon\cdot n_\epsilon)s_\epsilon {\rm{div}} n_\epsilon
+L_4s_\epsilon\nabla s_\epsilon (\nabla n_\epsilon) n_\epsilon+\frac{1}{\epsilon^2} W(s_\epsilon),
\end{eqnarray*}
we obtain that \\
\begin{itemize}
\item [(i)] if $L_1>0, L_2=L_4=0$, and $4L_1\alpha\ge 3L_3^2$, then
\begin{eqnarray*}
&&\widetilde{\mathcal{W}}_\epsilon(s_\epsilon, n_\epsilon, \nabla s_\epsilon,\nabla n_\epsilon)\\
&&\ge\alpha s_\epsilon^2|\nabla n_\epsilon|^2+\beta |\nabla s_\epsilon|^2+L_1 (\nabla s_\epsilon\cdot n_\epsilon)^2
+L_3  (\nabla s_\epsilon\cdot n_\epsilon) s_\epsilon {\rm{div}} n_\epsilon
+\frac{1}{\epsilon^2} W(s_\epsilon)\\
&&\ge \beta |\nabla s_\epsilon|^2+\frac{1}{\epsilon^2} W(s_\epsilon),
\end{eqnarray*}
\item [(ii)] if $L_2>0$, $L_1=L_3=0$, and $4L_2\alpha\ge L_4^2$, then
\begin{eqnarray*}
&&\widetilde{\mathcal{W}}_\epsilon(s_\epsilon, n_\epsilon, \nabla s_\epsilon,\nabla n_\epsilon)\\
&&\ge\alpha s_\epsilon^2|\nabla n_\epsilon|^2+\beta |\nabla s_\epsilon|^2+L_2 |\nabla s_\epsilon\wedge n_\epsilon|^2
+L_4  s_\epsilon\nabla s_\epsilon (\nabla n_\epsilon) n_\epsilon
+\frac{1}{\epsilon^2} W(s_\epsilon)\\
&&\ge \beta |\nabla s_\epsilon|^2+\frac{1}{\epsilon^2} W(s_\epsilon),
\end{eqnarray*}
\item [(iii)] if $L_1=L_2=L_3=L_4=0$, then
\begin{eqnarray*}
\widetilde{\mathcal{W}}_\epsilon(s_\epsilon, n_\epsilon, \nabla s_\epsilon,\nabla n_\epsilon)
&\ge& \alpha s_\epsilon^2|\nabla n_\epsilon|^2+\beta |\nabla s_\epsilon|^2+\frac{1}{\epsilon^2} W(s_\epsilon)\\
&\ge &\beta |\nabla s_\epsilon|^2+\frac{1}{\epsilon^2} W(s_\epsilon).
\end{eqnarray*}
\end{itemize}

Now we proceed by dividing the discussion into two separate cases:

\subsubsection{$\Omega\subset \R^3$ is a bounded domain}

For any $\delta>0$, define
$$\Omega^+_{\epsilon,\delta}=\Big\{x\in\Omega: \ |s_\epsilon-s_+|\le\delta\Big\};
\ \ \Omega^-_{\epsilon,\delta}=\Big\{x\in\Omega: \ |s_\epsilon| \le \delta\Big\},
$$
and
$$
{E}_{\epsilon,\delta}=\Omega\setminus\big(\Omega_{\epsilon,\delta}^+\cup \Omega_{\epsilon,\delta}^-\big)
=\Big\{x\in\Omega: \ |s_\epsilon| >\delta, \ \ |s_\epsilon-s_+| >\delta\Big\}.
$$

By the condition \eqref{W-cond} and Federer's co-area formula, we have that for any $0<\delta<\frac{s_+}2$,
\begin{eqnarray} \label{good_slice0}
\epsilon\mathcal{E}(s_\epsilon, n_\epsilon)&=& \int_\Omega \epsilon \widetilde{\mathcal{W}}_\epsilon(s_\epsilon, n_\epsilon, \nabla s_\epsilon,\nabla n_\epsilon)\,dx\nonumber\\
&\ge& \int_\Omega\big(\epsilon\beta |\nabla s_\epsilon|^2+\frac{1}{\epsilon} W(s_\epsilon)\big)\,dx\nonumber\\
&\ge & 2\sqrt{\beta}\int_\Omega \sqrt{W(s_\epsilon)}|\nabla s_\epsilon|\,dx\nonumber\\
&\ge & 2\sqrt{\beta}\Big(\int_{\frac{\delta}2}^\delta \sqrt{W(\tau)} \mathcal{H}^2(\partial \Omega^-_{\epsilon,\tau}\cap\Omega)\,d\tau+\int_{s_+-\delta}^{s_+-\frac{\delta}2}\sqrt{W(\tau)}\mathcal{H}^2(\partial \Omega^+_{\epsilon,\tau}\cap\Omega)\,d\tau \Big)\nonumber\\
&\ge& 2\sqrt{\beta}C_\delta\Big(\int_{\frac{\delta}2}^\delta \mathcal{H}^2(\partial \Omega^-_{\epsilon,\tau}\cap\Omega)\,d\tau+\int_{s_+-\delta}^{s_+-\frac{\delta}2}\mathcal{H}^2(\partial \Omega^+_{\epsilon,\tau}\cap\Omega)\,d\tau \Big).
\end{eqnarray}
Therefore, by Fubini's theorem there exists $\delta_*\in (\frac{\delta}2, \delta)$
such that
\begin{equation}\label{good_slice1}
\mathcal{H}^2(\partial \Omega^-_{\epsilon,\delta_*}\cap\Omega)+ \mathcal{H}^2(\partial \Omega^+_{\epsilon,\delta_*}\cap\Omega)\le C(\beta, \Lambda, \delta).
\end{equation}
From \eqref{W-cond}, we know that there exists $C_\delta>0$ such that
\begin{equation}\label{jump_set}
\big|E_{\epsilon,\delta}\big|\le \frac{1}{C_\delta}\int_\Omega W(s_\epsilon)\,dx\le \frac{\Lambda\epsilon}{C_\delta}\rightarrow 0.
\end{equation}

From \eqref{good_slice1},  there exist two subsets $E^\pm\subset\Omega$ with finite perimeters
in $\Omega$ such that, after passing to a subsequence,
$$\chi_{\Omega_{\epsilon,\delta_*}^\pm}\rightharpoonup \chi_{E^\pm}
\ \ {\rm{in}}\ \ {\rm{BV}}(\R^3)\ \ {\rm{and}}\ \  \chi_{\Omega_{\epsilon,\delta_*}^\pm}\rightarrow \chi_{E^\pm}
\ \ {\rm{in}}\ \ L^1(\R^3).
$$
This and \eqref{jump_set} imply that
$$\big|\Omega\setminus(E^+\cup E^-)\big|=\big|E^+\cap E^-\big|=0$$
so that $\Omega=E^+\cup E^-$ (modular a set of zero Lebesgue measure).

Define an auxiliary function $\phi: (-\frac12, 1)\to \R$ by letting
$$\phi(t)=2\sqrt{\beta}\int_0^t \sqrt{W(\tau)}\,d\tau, \ t\in (-\frac12, 1).$$
Notice that the $1$-dimensional minimal connecting energy $\alpha_0=\phi(s_+)$.
It follows from \eqref{good_slice0} that
\begin{eqnarray}\label{BV_bound1}
\int_\Omega |\nabla (\phi(s_\epsilon))|\,dx
&\le& 2\sqrt{\beta}\int_\Omega \sqrt{W(s_\epsilon)}|\nabla s_\epsilon|\,dx\nonumber\\
&\le&\int_\Omega \epsilon \widetilde{\mathcal{W}}_\epsilon(s_\epsilon, n_\epsilon, \nabla s_\epsilon, \nabla n_\epsilon)\,dx\le \Lambda+o(1).
\end{eqnarray}
From the boundary condition \eqref{b_cond1.1}, we know that
\begin{equation}\label{bdry_value1}
\phi(s_\epsilon)\rightarrow 0\ \ {\rm{in}}\ \ L^2(\Sigma^-)
\ \ {\rm{and}}\ \ \phi(s_\epsilon)\rightarrow \phi(s_+)\ \ {\rm{in}}\ \ L^2(\Sigma^+),
\ \ {\rm{as}}\ \ \epsilon\rightarrow 0.
\end{equation}
In particular, we may assume that
$$\sup_{0<\epsilon<1}\int_{\partial\Omega} |\phi(s_\epsilon)|^2\,d\mathcal{H}^2\le
C<\infty.$$
Hence by the Poincar\'e inequality we have that
\begin{eqnarray}\label{L1_bound1}
\int_\Omega |\phi(s_\epsilon)|\,dx \le C\big(\int_\Omega |\nabla (\phi(s_\epsilon))|+
\int_{\partial\Omega} |\phi(s_\epsilon)|\,d\mathcal{H}^2\big)
\le C(1+\Lambda).
\end{eqnarray}
It follows from \eqref{BV_bound1} and \eqref{L1_bound1} that there exists
$\psi\in BV(\Omega)$ such that after passing to a subsequence,
$\phi(s_\epsilon)\rightarrow\psi$ weakly in ${\rm{BV}}(\Omega)$ and strongly in
$L^1(\Omega)\cap L^1(\partial\Omega)$. By the lower semicontinuity, we have that
\begin{equation}\label{BV_bound2}
\begin{split}
|D\psi|(\Omega)&\le\liminf_{\epsilon\to 0}\int_\Omega|\nabla(\phi(s_\epsilon))|\,dx
\le \liminf_{\epsilon\rightarrow 0} \int_{\Omega} \big(\epsilon|\nabla s_\epsilon|^2+\frac{1}{\epsilon}W(s_\epsilon)\big)\,dx\\
&\le\liminf_{\epsilon\rightarrow 0}\int_\Omega \epsilon\widetilde{\mathcal{W}}_\epsilon
(s_\epsilon, n_\epsilon, \nabla s_\epsilon, \nabla n_\epsilon)\,dx\le \Lambda.
\end{split}
\end{equation}
It follows from \eqref{bdry_value1} that
$\psi=0$ on $\Sigma^-$, and $\psi=\phi(s_+)$ on $\Sigma^+$. We claim that
\begin{equation}\label{interior_value}
\psi(x)=\begin{cases} 0 & x\in E^-,\\ \phi(s_+) & x\in E^+. \end{cases}
\end{equation}
To see this, observe by Fatou's lemma that
\begin{eqnarray*}
&&\int_{E^-} |\psi|^2\,dx=\int_{\R^3}|\psi|^2\chi_{E^-}\,dx\\
&&\le \liminf_{\epsilon\rightarrow 0} \int_{\R^3} |\phi(s_\epsilon)|^2 \chi_{\Omega^{-}_{\epsilon,\delta_*}}\,dx\\
&&= \liminf_{\epsilon\rightarrow 0}\int_{\Omega^-_{\epsilon,\delta_*}}|\phi(s_\epsilon)|^2\,dx\\
&&\le C\big\|\phi'\big\|_{L^\infty([-\delta_*, \delta_*])}^2\liminf_{\epsilon\rightarrow 0}\int_{\Omega^-_{\epsilon,\delta_*}}|s_\epsilon|^2\,dx\\
&&\le C\liminf_{\epsilon\rightarrow 0}\int_{\Omega^-_{\epsilon,\delta_*}}W(s_\epsilon)\,dx
\le C\Lambda \epsilon\rightarrow 0,
\end{eqnarray*}
where we have used the fact that $W(\tau)\approx \tau^2$ for $|\tau|\le\delta_*$.
Thus $\psi=0$ a.e. in $E^-$.

Similarly, by using the fact that $W(\tau)\approx |\tau-s_+|^2$ for $|\tau-s_+|\le\delta_*$, we can estimate
\begin{eqnarray*}
&&\int_{E^+} |\psi-\phi(s_+)|^2\,dx=\int_{\R^3}|\psi-\phi(s_+)|^2\chi_{E^+}\,dx\\
&&\le \liminf_{\epsilon\rightarrow 0} \int_{\R^3} |\phi(s_\epsilon)-\phi(s_+)|^2 \chi_{\Omega^{+}_{\epsilon,\delta_*}}\,dx\\
&&= \liminf_{\epsilon\rightarrow 0}\int_{\Omega^+_{\epsilon,\delta_*}}|\phi(s_\epsilon)-\phi(s_+)|^2\,dx\\
&&\le C\|\phi'\|_{L^\infty([s_+-\delta_*, s_+ +\delta_*])}\liminf_{\epsilon\rightarrow 0}
\int_{\Omega^+_{\epsilon,\delta_*}}|s_\epsilon-s_+|^2\,dx\\
&&\le C\liminf_{\epsilon\rightarrow 0}\int_{\Omega^+_{\epsilon,\delta_*}}W(s_\epsilon)\,dx
\le C\Lambda \epsilon\rightarrow 0,
\end{eqnarray*}
this yields that $\psi=\phi(s_+)$ a.e. in $E^+$.

It follows from \eqref{interior_value} that
\begin{equation}\label{nematic_region1}
E^+=\big\{x\in\Omega: \ \psi(x)\ge t\big\}, \ \forall t\in (0, \alpha_0).
\end{equation}
In what following, for a subset $E\subset\R^3$
we denote by $[[E]]$  the corresponding $3$-dimensional current (through integration),
and denote by $\partial[[E]]$ the boundary current of $[[E]]$.
%Since \footnote{need to specify more on the boundary value $s_\epsilon$}

Since $s_\epsilon\in H^1(\partial\Omega, (-\frac12, 1))$ satisfies \eqref{b_cond1.1},
$$
\partial\big{[\big[}\partial[[\Omega_{\epsilon,\delta_*}^+]]\big\lfloor~\Omega\big{]\big]}
=[[\{x\in\partial\Omega: \ s_\epsilon(x)=s_+-\delta_*\}]]\rightarrow [[\Sigma^0]], \ {\rm{as}}\ \epsilon\rightarrow 0,
$$
holds as weak convergence of currents,
we obtain that
\begin{equation}\label{bdry_nematic_region}
\partial\big{[\big[}\partial[[E^+]]\big\lfloor~\Omega\big{]\big]}=[[\Sigma^0]].
\end{equation}
It follows from \eqref{nematic_region1}, \eqref{bdry_nematic_region}, and the area minimality of $\Gamma$
that
\begin{equation}\label{minimality1}
\mathcal{H}^2\big(\partial^*\big\{x\in\Omega: \ \psi(x)\ge t\big\}\lfloor~\Omega\big)\ge \mathcal{H}^2(\Gamma),
\ \forall t\in [0, \alpha_0).
\end{equation}
Here $\partial^*E$ denotes the reduced boundary of a set $E$ of finite perimeter.
By the co-area formula for BV functions and \eqref{minimality1}, we then have
\begin{eqnarray}\label{LB2}
\Lambda=\lim_{\epsilon_i\to 0} \epsilon\mathcal{E}(s_{\epsilon_i}, n_{\epsilon_i})\ge |D\psi|(\Omega)&=&\int_{\R} \mathcal{H}^2\big(\partial^*\big\{x\in\Omega: \ \psi(x)>t\big\}\lfloor \Omega\big)\,dt\nonumber\\
&\ge&\int_{0}^{\alpha_0} \mathcal{H}^2\big(\partial^*\big\{x\in\Omega: \ \psi(x)>t\big\}\lfloor \Omega\big)\,dt
\nonumber\\
&\ge& \alpha_0 \mathcal{H}^2(\Gamma).
\end{eqnarray}
This proves the part ``$\ge$"  of  \eqref{minimality0} in  Theorem \ref{existence}, when $\Omega$ is a bounded
domain in $\R^3$.

\subsubsection{$\Omega=\R^3$ and $\big|\{x\in\R^3: \ s_\epsilon(x)\ge {s_+}\}\big|=|B_1|$}

First notice that
%the conditions \eqref{W-cond} and \eqref{b_cond1.2} imply that|
%$0\le s_\epsilon<1$ in $\R^3$, and for any $0<\delta<s_+$ there exists $C(\delta)>0$ such that
%$$\Big|\big\{x\in \R^3: \ \delta\le s_\epsilon(x)\le s_+-\delta\big\}\Big|
%\le \frac{1}{C(\delta)}\int_{\R^3} W(s_\epsilon)\,dx\le \frac{\Lambda\epsilon}{C(\delta)}\rightarrow 0,\ {\rm{as}}\ \epsilon\to 0.$$
%This, combined with $\big|\{x\in\R^3: s_\epsilon(x)>\frac{s_+}2\}\big|=|B_1|$, implies that
for any $0<\delta\le s_+$,
\begin{equation}\label{vol_constr1}
\Big|\big\{x\in\R^3: s_\epsilon(x)\ge\delta\big\}\Big|\ge\Big|\big\{x\in\R^3: s_\epsilon(x)\ge s_+\big\}\Big|\ge |B_1|
%+O(\frac{\epsilon}{C(\delta)}).
\end{equation}
As in \eqref{BV_bound1} of the previous subsection, we can obtain that for any $0<\delta<\frac{s_+}2$,
\begin{eqnarray*}
\Lambda+o(1)&\ge& \int_{\R^3} \epsilon \widetilde{\mathcal{W}}(s_\epsilon, n_\epsilon, \nabla s_\epsilon,
\nabla n_\epsilon)\,dx\\
&\ge&{2}\int_{\R^3}\sqrt{\beta W(s_\epsilon)}|\nabla s_\epsilon|\,dx\\
&\ge& {2}\int_{0}^{s_+} \sqrt{\beta W(\tau)}\mathcal{H}^2\big(\partial^*\big\{x\in\R^3: \ s_\epsilon(x)\ge\tau\big\}\big)\,d\tau
\end{eqnarray*}
By the isoperimetric inequality, we have that for any $0<\tau\le s_+$,
\begin{eqnarray*}
\mathcal{H}^2\big(\partial^*\big\{x\in\R^3: \ s_\epsilon(x)\ge\tau\big\}\big)
&\ge& (36\pi)^\frac13 \Big|\big\{x\in\R^3: \ s_\epsilon(x)\ge\tau\big\}\Big|^\frac23\\
&\ge&(36\pi)^\frac13|B_1|^\frac23\\
%&=&(36\pi)^\frac13\big(|B_1|+O(\frac{\epsilon}{C(\delta)})\big)^\frac23.
\end{eqnarray*}
Hence we obtain that
\begin{eqnarray*}
\Lambda+o(1)&\ge& s_+(36\pi)^\frac13|B_1|^\frac23=s_+ \mathcal{H}^2(\partial B_1).
\end{eqnarray*}
%Taking $\epsilon\to 0$, this yields
%$$\Lambda\ge (\phi(s_+-\delta)-\phi(\delta)) (36\pi)^\frac13|B_1|^\frac23=(\phi(s_+-\delta)-\phi(\delta)) \mathcal{H}^2(\partial B_1).$$
%Sending $\delta\to 0$,
This proves ``$\ge$" of Theorem \ref{existence}
for $\Omega=\R^3$. \qed
\subsection{Upper bound estimates}

The upper bound estimates are based on concrete constructions, similar to that by \cite{LPW}.
We will first discuss the construction for a bounded domain $\Omega\subset\R^3$.

\subsubsection{$\Omega\subset\R^3$ is a bounded domain}

We need to introduce some notations.
Fix a large constant $L>0$, whose value will be determined later, we may assume for simplicity that
$Q_{L\epsilon}^\pm\approx\Sigma_{L\epsilon}^\pm\times [0,L\epsilon]$.
We will construct
a function $\hat{s}_\epsilon$ in $Q_{L\epsilon}^-$, that is a linear interpolation of
of $t_\epsilon\big|_{\Sigma_{L\epsilon}^-\times \{0\}}$
and $0\big|_{\Sigma_{L\epsilon}^-\times \{L\epsilon\}}$, i.e.,
\begin{equation}\label{ext1}
\hat{s}_\epsilon(x, t)=\frac{L\epsilon-t}{L\epsilon}t_\epsilon(x), \ (x,t)\in \Sigma_{L\epsilon}^-\times [0,L\epsilon].
\end{equation}
Similarly, a function $\hat{s}_\epsilon$ in $Q_{L\epsilon}^+$ is constructed by a linear interpolation
of $t_\epsilon\big|_{\Sigma_{L\epsilon}^+\times \{0\}}$
and $s_+\big|_{\Sigma_{L\epsilon}^+\times \{L\epsilon\}}$, i.e.,
\begin{equation}\label{ext2}
\hat{s}_\epsilon(x, t)=\frac{L\epsilon-t}{L\epsilon}t_\epsilon(x)+\frac{t}{L\epsilon} s_+,
\ (x,t)\in \Sigma_{L\epsilon}^+\times [0,L\epsilon].
\end{equation}
%Notice that $\tilde{n}_\epsilon$ does not necessarily map $Q_\delta^\pm$ to $\mathbb S^2$.  However, since
%$\mathbb S^2$ is simply connected, we can apply Hardt-Lin's extension Lemma (see \cite{HL1}) as follows.
%$\exists a\in\R^3$, with $|a|\le \frac12$, such that
%\begin{equation}\label{HL_ext}
%\int_{Q_\delta^\pm} \hat{s}_\epsilon^2\big|\nabla \big((\Pi_a\big|_{\mathbb S^2})^{-1} (\Pi_a (\tilde{n}_\epsilon))\big)\big|^2\,dx\lesssim \int_{Q_\delta^\pm} \hat{s}_\epsilon^2\big|\nabla \tilde{n}_\epsilon\big|^2\,dx
%\end{equation}
%Here $\Pi_a(y)=\frac{y-a}{|y-a|}: B_1^3\to\mathbb S^2$ for $y\in B_1^3$. Denote
Let $\widehat{n}_\epsilon\in H^1(\Omega,\mathbb S^2)$ be given by \eqref{b_cond1.20}. Then,
%\begin{equation}\label{ext3}
%\hat{n}_\epsilon(x)=(\Pi_a\big|_{\mathbb S^2})^{-1} (\Pi_a (\tilde{n}_\epsilon))(x), \ x\in Q_\delta^\pm.
%\end{equation}
by direct calculations and applying \eqref{b_cond1.1} and \eqref{W-cond}, we have that
\begin{eqnarray}\label{ext_est1}
&&\int_{Q_{L\epsilon}^-} \big(\beta|\nabla\hat{s}_\epsilon|^2+\frac{1}{\epsilon^2}W(\hat{s}_\epsilon)\big)\,dx
\lesssim\epsilon\int_{\Sigma^-} \big(|\nabla t_\epsilon|^2
+\frac{1}{\epsilon^2}W({t}_\epsilon)\big)\,d\mathcal{H}^2
+ \frac{1}{\epsilon}\int_{\Sigma^-}|t_\epsilon|^2\,d\mathcal{H}^2,
\end{eqnarray}
and
\begin{eqnarray}\label{ext_est2}
&&\int_{Q_{L\epsilon}^+}(\beta|\nabla\hat{s}_\epsilon|^2+\frac{1}{\epsilon^2}W(\hat{s}_\epsilon))
\lesssim \epsilon\int_{\Sigma^+} (|\nabla t_\epsilon|^2+\frac{1}{\epsilon^2}W({t}_\epsilon))\,d\mathcal{H}^2
+ \frac{1}{\epsilon}\int_{\Sigma^+} |t_\epsilon-s_+|^2\,d\mathcal{H}^2.
\end{eqnarray}
By Fubini's theorem, there exists $L_1\in [L, 2L]$ such that
\begin{equation}\label{slice3}
\begin{cases}\displaystyle
L_1\epsilon\int_{\partial \Gamma_{L_1\epsilon}^-\cap\Omega} \big(\beta |\nabla\hat{s}_\epsilon|^2+\frac{1}{\epsilon^2}W(\hat{s}_\epsilon)\big)\,d\mathcal{H}^2
\lesssim \int_{Q_{L\epsilon}^-} \big(\beta |\nabla\hat{s}_\epsilon|^2+\frac{1}{\epsilon^2}W(\hat{s}_\epsilon)\big)\,dx,\\
\displaystyle L_1\epsilon\int_{\partial \Gamma_{L_1\epsilon}^+\cap\Omega} \big(\beta|\nabla\hat{s}_\epsilon|^2+\frac{1}{\epsilon^2}W(\hat{s}_\epsilon)\big)\,d\mathcal{H}^2
\lesssim \int_{Q_{L\epsilon}^+} \big(\beta|\nabla\hat{s}_\epsilon|^2+\frac{1}{\epsilon^2}W(\hat{s}_\epsilon)\big)\,dx.
\end{cases}
\end{equation}
It follows from the regularity theorem of area minimizing surfaces (see \cite{Federer} and \cite{HS})
that $\Gamma\in C^\infty(\overline\Omega)$. Let $\xi\in C^\infty\big([-L_1,
L_1]\big)$ be
an almost minimal $1$-dimensional connecting orbit, i.e., $\xi(-L_1)=0$,
$\xi(L_1)=s_+$, and
\begin{equation}\label{approx_minimal}
\int_{-L_1}^{L_1} \big(\beta\dot\xi^2+W(\xi)\big)\,d\tau= \alpha_0+o_L(1), \ \lim_{L\to\infty}o_L(1)=0.
\end{equation}
Define $\hat{s}_\epsilon: \Omega_{L\epsilon}\cap \Gamma_{L_1\epsilon}\to \R$
by letting
\begin{equation}\label{ext4}
\hat{s}_\epsilon(x)=\xi\big(\frac{d_\Gamma(x)}{\epsilon}\big),
\ \ x\in \Omega_{L\epsilon}\cap \Gamma_{L_1\epsilon}.
\end{equation}
%it is easy to see that
%\begin{equation}\label{ext_est3}
%\int_{\Omega_{\delta}\cap \Gamma_{\delta_1}}\hat{s}_\epsilon^2|\nabla\hat{n}_\epsilon|^2\,dx=0,
%\end{equation}
By the co-area formula and the fact that $|\nabla d_{\Gamma}(x)|=1$ for $x\in \Omega_{L\epsilon}\cap \Gamma_{L_1\epsilon}$, we can estimate
\begin{eqnarray}\label{ext_est4}
&&\int_{\Omega_{L\epsilon}\cap \Gamma_{L_1\epsilon}}
\big(\beta|\nabla\hat{s}_\epsilon|^2+\frac{1}{\epsilon^2}W(\hat{s}_\epsilon)\big)\,dx
=\frac{1}{\epsilon^2}\int_{\Omega_{L\epsilon}\cap \Gamma_{L_1\epsilon}}
\big(\beta\dot{\xi}^2+W(\xi)\big)(\frac{d_\Gamma(x)}{\epsilon})\,dx\nonumber\\
&&=\frac{1}{\epsilon}\int_{-L_1}^{L_1}\big(\beta\dot{\xi}^2+W(\xi)\big)(\tau)\mathcal{H}^2\big(\big\{x\in\Omega_{L\epsilon}: \ d_{\Gamma}(x)=\epsilon\tau\big\}\big)\,d\tau\nonumber\\
&&\le \frac{1}{\epsilon}(\mathcal{H}^2(\Gamma)+o_\epsilon(1))
\int_{-L_1}^{L_1}\big(\beta\dot{\xi}^2+W(\xi)\big)\,d\tau%\nonumber\\
\le \frac{1}{\epsilon}(\mathcal{H}^2(\Gamma)+o_\epsilon(1))(\alpha_0+o_L(1)),
\end{eqnarray}
where we have the fact that the surface $\big\{x\in\Omega: \ d_\Gamma(x)=\tau\big\}$ converges
to $\Gamma$ in $C^2$-norm, as $\tau\rightarrow 0$. Hence it holds
that
$$
\mathcal{H}^2\big(\big\{x\in\Omega_{L\epsilon}: \ d_{\Gamma}(x)=\epsilon\tau\big\}\big)\le
\mathcal{H}^2(\Gamma)+o_\epsilon(1), \ -L_1<\tau< L_1, \ {\rm{here}}\ \lim_{\epsilon\to 0}o_\epsilon(1)=0.
$$
It is not hard to check that
\begin{equation}\label{slice4}
\int_{\partial\Omega_{L\epsilon}\cap\Gamma_{L_1\epsilon}}\big(\beta|\nabla_{\rm{tan}}\hat{s}_\epsilon|^2
+\frac{1}{\epsilon^2}W(\hat{s}_\epsilon)\big)\,d\mathcal{H}^2
\le \frac{C}{\epsilon}(\alpha_0+o_L(1)).
\end{equation}
In the regions $\Omega_{L\epsilon}^\pm\setminus \Gamma_{L_1\epsilon}$, we simply define
\begin{equation}\label{ext5}
\hat{s}_\epsilon= 0  \ {\rm{in}}\ \Omega_{L\epsilon}^-\setminus \Gamma_{L_1\epsilon}; \ \
\hat{s}_\epsilon=s_+  \ {\rm{in}}\ \Omega_{L\epsilon}^+\setminus \Gamma_{L_1\epsilon},
\end{equation}
so that
\begin{equation}\label{ext_est5}
\int_{\Omega_{L\epsilon}^\pm\setminus \Gamma_{L_1\epsilon}}
\big(\beta|\nabla\hat{s}_\epsilon|^2+\frac{1}{\epsilon^2}W(\hat{s}_\epsilon)\big)\,dx=0.
\end{equation}
It remains to construct $\hat{s}_\epsilon$ in the region $(\partial\Omega)_{L\epsilon}\cap \Gamma_{L_1\epsilon}$, which can be roughly viewed as a ball of radius $L_1\epsilon$ centered at $x_*\in\Omega$.
Hence we can do a homogeneous of degree zero extension of
$\hat{s}_\epsilon\ $% from  $\partial\big((\partial\Omega)_\delta\cap \Gamma_{\delta_1}\big)$
with respect to the center $x_*$, i.e.,
\begin{equation}\label{ext6}
\hat{s}_\epsilon(x)=\hat{s}_\epsilon\big(L_1\epsilon\frac{x-x_*}{|x-x_*|}\big),
\ x\in (\partial\Omega)_{L\epsilon}\cap \Gamma_{L_1\epsilon}.
\end{equation}
Hence by \eqref{slice3}, \eqref{slice4}, \eqref{ext_est1}, and \eqref{ext_est2},
we have that
\begin{eqnarray}\label{ext_est6}
&&\int_{(\partial\Omega)_{L\epsilon}\cap \Gamma_{L_1\epsilon}}
\big(\beta|\nabla\hat{s}_\epsilon|^2+\frac{1}{\epsilon^2}W(\hat{s}_\epsilon)\big)\,dx%\nonumber\\
\le L\epsilon \int_{\partial[(\partial\Omega)_{L\epsilon}\cap \Gamma_{L_1\epsilon}]}
\big(\beta|\nabla_{\rm{tan}}\hat{s}_\epsilon|^2
+\frac{1}{\epsilon^2}W(\hat{s}_\epsilon)\big)\,d\mathcal{H}^2\nonumber\\
&&=L\epsilon\Big\{\int_{\partial\Omega_{L\epsilon}\cap\Gamma_{L_1\epsilon}}
+\int_{\partial\Omega\cap\Gamma_{L_1\epsilon}}+\int_{\partial\Gamma_{L_1\epsilon}^-\cap\Omega}
+\int_{\partial\Gamma_{L_1\epsilon}^+\cap\Omega}\Big\}
\big(\beta|\nabla_{\rm{tan}}\hat{s}_\epsilon|^2
+\frac{1}{\epsilon^2}W(\hat{s}_\epsilon)\big)\,d\mathcal{H}^2\nonumber\\
&&\le CL(\alpha_0+o_L(1))+
\Big\{\int_{Q_{L\epsilon}^{-}} +\int_{Q_{L\epsilon}^+}\Big\}\big(\beta|\nabla\hat{s}_\epsilon|^2
+\frac{1}{\epsilon^2}W(\hat{s}_\epsilon)\big)\,dx\nonumber\\
&&\quad+CL\epsilon\int_{\partial\Omega\cap\Gamma_{L_1\epsilon}}  \big(\beta|\nabla_{\rm{tan}}t_\epsilon|^2
+\frac{1}{\epsilon^2}W({t}_\epsilon)\big)\,d\mathcal{H}^2\nonumber\\
&&\le CL+CL\epsilon\int_{\partial \Omega} \big(\beta|\nabla_{\rm{tan}}t_\epsilon|^2
+\frac{1}{\epsilon^2}W({t}_\epsilon)\big)\,d\mathcal{H}^2\nonumber\\
&&\quad+\frac{C}{\epsilon}\Big(\int_{\Sigma^-}|t_\epsilon|^2\,d\mathcal{H}^2
+\int_{\Sigma^+}|t_\epsilon-s_+|^2\,d\mathcal{H}^2\Big)\nonumber\\
%&&\quad+\frac{C}{\delta}\Big(\int_{\Sigma_{\delta}^-}|g_\epsilon-p_0|^2\,d\mathcal{H}^2+\int_{\Sigma_{\delta}^+}|g_\epsilon-p_0|^2\,d\mathcal{H}^2\Big)\nonumber\\
&&\le C+\frac{o_\epsilon(1)}{\epsilon},
\end{eqnarray}
where we have used the conditions \eqref{b_cond1.1} and \eqref{b_cond1.2} in the last step.

Finally, by putting together \eqref{ext1}, \eqref{ext2}, \eqref{ext4},  \eqref{ext5},
and \eqref{ext6}, we find an
extension map $\hat{s}_\epsilon:\Omega\to\R$.
Furthermore, by the estimates \eqref{ext_est1}, \eqref{ext_est2},
\eqref{ext_est4}, \eqref{ext_est5}, \eqref{ext_est6}, we see that $\hat{s}_\epsilon$
satisfies the estimate:
\begin{eqnarray}\label{ext_est11}
\int_{\Omega}
\big(\beta|\nabla\hat{s}_\epsilon|^2+\frac{1}{\epsilon^2}W(\hat{s}_\epsilon)\big)\,dx
\le \frac{1}{\epsilon}(\mathcal{H}^2(\Gamma)+o_\epsilon(1))(\alpha_0+o_L(1)) +C
+\frac{o_\epsilon(1)}{\epsilon}.
\end{eqnarray}
It follows from \eqref{b_cond1.20} that
\begin{eqnarray}\label{ext_est12}
\int_\Omega \hat{s}_\epsilon^2(W_{OF}(\widehat{n}_\epsilon, \nabla\widehat{n}_\epsilon)
+\alpha |\nabla\widehat{n}_\epsilon|^2)\,dx\le Cs_+^2\int_\Omega |\nabla\widehat{n}_\epsilon|^2\,dx
\le \frac{o_\epsilon(1)}{\epsilon}.
\end{eqnarray}
To estimate the contributions from the terms involving the interactive energies between
$\nabla \hat{s}_\epsilon$ and $\widehat{n}_\epsilon$, we proceed as follows:
\begin{itemize}
\item [({\bf A})] $L_1>L_2=L_4=0$. We can estimate
\begin{eqnarray*}
&&\big|\int_\Omega \big(L_1(\nabla\hat{s}_\epsilon\cdot \widehat{n}_\epsilon)^2
+L_3(\nabla\hat{s}_\epsilon\cdot \widehat{n}_\epsilon)(\hat{s}_\epsilon{\rm{div}}\widehat{n}_\epsilon)\big)\,dx\big|\\
&&\le C\big(\int_{\Omega_{L\epsilon}\cap\Gamma_{L_1\epsilon}}|\nabla \hat{s}_\epsilon\cdot\widehat{n}_\epsilon|^2\,dx+\int_{\Omega\setminus\Omega_{L\epsilon}}|\nabla \hat{s}_\epsilon|^2\,dx+\int_\Omega |\nabla\widehat{n}_\epsilon|^2\,dx\big).
\end{eqnarray*}
Notice that \eqref{ext_est1}, \eqref{ext_est2}, and \eqref{ext_est6} imply
\begin{eqnarray*}
\int_{\Omega\setminus\Omega_{L\epsilon}}|\nabla \hat{s}_\epsilon|^2\,dx&\le&
CL\epsilon \int_{\partial\Omega} |\nabla t_\epsilon|^2\,d\mathcal{H}^2+
C\epsilon^{-1}\big(\int_{\Sigma^-}|t_\epsilon|^2\,d\mathcal{H}^2+\int_{\Sigma^+}|t_\epsilon-s_+|^2\,d\mathcal{H}^2\big)\\
&\le& C+\frac{o_\epsilon(1)}{\epsilon}.
\end{eqnarray*}
Since $|\dot{\xi}(t)|\le C_L$ in $t\in [-L, L]$ and $\widehat{n}_\epsilon\cdot\nu_\Gamma=0$
on $\Gamma$, we can estimate
\begin{eqnarray}\label{ext_est19}
&&\int_{\Omega_{L\epsilon}\cap\Gamma_{L_1\epsilon}}|\nabla \hat{s}_\epsilon\cdot\widehat{n}_\epsilon|^2\,dx
=\frac{1}{\epsilon^2}\int_{\Omega_{L\epsilon}\cap\Gamma_{L_1\epsilon}}\dot{\xi}^2\big(\frac{d_\Gamma(x)}{\epsilon}\big)
(\nabla d_\Gamma(x)\cdot\widehat{n}_\epsilon(x))^2\,dx\nonumber\\
&&\le C\epsilon^{-2}\int_{\Omega_{L\epsilon}\cap\Gamma_{L_1\epsilon}}(\nabla d_\Gamma(x)\cdot\widehat{n}_\epsilon(x))^2\,dx\nonumber\\
&&= C\epsilon^{-2}\int_{\Omega_{L\epsilon}\cap\Gamma_{L_1\epsilon}}(\nabla d_\Gamma(x)\cdot\widehat{n}_\epsilon(x)-\nabla d_{\Gamma}(\Pi_\Gamma(x))\cdot\widehat{n}_\epsilon(\Pi_\Gamma(x)))^2\,dx\nonumber\\
&&\le C\epsilon^{-2}\int_{\Omega_{L\epsilon}\cap\Gamma_{L_1\epsilon}}
\big(|\nabla d_\Gamma(x)-\nabla d_{\Gamma}(\Pi_\Gamma(x))|^2
+|\widehat{n}_\epsilon(x)-\widehat{n}_\epsilon(\Pi_\Gamma(x))|^2\big)\,dx\nonumber\\
&&\le C\int_{\Omega_{L\epsilon}\cap\Gamma_{L_1\epsilon}}
\big(|\nabla^2 d_\Gamma(x)|^2+|\nabla\widehat{n}_\epsilon(x)|^2\big)\,dx\nonumber\\
&&\le C\epsilon+\frac{o_\epsilon(1)}{\epsilon},
\end{eqnarray}
where $\Pi_{\Gamma}: \Gamma_{L_1\epsilon}\to\Gamma$ is the (smooth) nearest point projection.
Hence we have that
\begin{eqnarray}\label{ext_est13}
\big|\int_\Omega \big(L_1(\nabla\hat{s}_\epsilon\cdot \widehat{n}_\epsilon)^2
+L_3(\nabla\hat{s}_\epsilon\cdot \widehat{n}_\epsilon)(\hat{s}_\epsilon{\rm{div}}\widehat{n}_\epsilon)\big)\,dx\big|
\le C(1+\epsilon)+\frac{o_\epsilon(1)}{\epsilon}.
\end{eqnarray}
Adding \eqref{ext_est11}, \eqref{ext_est12}, and \eqref{ext_est13} together, we arrive at
\begin{equation}\label{SUB1.0}
\epsilon\int_\Omega \widetilde{\mathcal{W}}(\hat{s}_\epsilon, \widehat{n}_\epsilon,
\nabla\hat{s}_\epsilon, \nabla\widehat{n}_\epsilon)\,dx
\le (\alpha_0+o_L(1))(\mathcal{H}^2(\Gamma)+o_\epsilon(1))+C\epsilon +o_\epsilon(1).
\end{equation}
After first sending $\epsilon\to 0$  and  then $L\to\infty$, and using the minimality of
$(s_\epsilon, n_\epsilon)$ we have  that
\begin{equation}\label{SUB1.1}
\limsup_{\epsilon\to 0}\int_\Omega \epsilon\widetilde{\mathcal{W}}({s}_\epsilon, {n}_\epsilon,
\nabla {s}_\epsilon, \nabla {n}_\epsilon)\,dx\le
\limsup_{\epsilon\to 0}\int_\Omega \epsilon\widetilde{\mathcal{W}}(\hat{s}_\epsilon, \widehat{n}_\epsilon,
\nabla\hat{s}_\epsilon, \nabla\widehat{n}_\epsilon)\,dx
\le \alpha_0\mathcal{H}^2(\Gamma).
\end{equation}

\item [({\bf B})] $L_2>L_1=L_3=0$. We can bound
\begin{eqnarray*}
&&\big|\int_\Omega\big(L_2|\nabla \hat{s}_\epsilon\wedge \widehat{n}_\epsilon|^2
+L_4 \hat{s}_\epsilon\nabla \hat{s}_\epsilon\cdot (\nabla \widehat{n}_\epsilon)\widehat{n}_\epsilon\big)\,dx\big|\\
&&\le C\big(\int_{\Omega_{L\epsilon}\cap\Gamma_{L_1\epsilon}}
|\nabla \hat{s}_\epsilon\wedge \widehat{n}_\epsilon|^2\,dx+\int_{\Omega\setminus\Omega_{L\epsilon}}|\nabla\hat{s}_\epsilon|^2\,dx+\int_\Omega |\nabla \widehat{n}_\epsilon|^2\,dx\big)\\
&&\le C(1+\epsilon)+\frac{o_\epsilon(1)}{\epsilon}+C\int_{\Omega_{L\epsilon}\cap\Gamma_{L_1\epsilon}}
|\nabla \hat{s}_\epsilon\wedge \widehat{n}_\epsilon|^2\,dx.
\end{eqnarray*}
Since $\widehat{n}_\epsilon\wedge\nu_\Gamma=0$ on $\Gamma$,
the last term in the right hand side can be estimated as follows.
\begin{eqnarray*}
&&\int_{\Omega_{L\epsilon}\cap\Gamma_{L_1\epsilon}}|\nabla \hat{s}_\epsilon\wedge\widehat{n}_\epsilon|^2\,dx
=\frac{1}{\epsilon^2}\int_{\Omega_{L\epsilon}\cap\Gamma_{L_1\epsilon}}\dot{\xi}^2\big(\frac{d_\Gamma(x)}{\epsilon}\big)
|\nabla d_\Gamma(x)\wedge\widehat{n}_\epsilon(x)|^2\,dx\\
&&\le C\epsilon^{-2}\int_{\Omega_{L\epsilon}\cap\Gamma_{L_1\epsilon}}|\nabla d_\Gamma(x)\wedge\widehat{n}_\epsilon(x)|^2\,dx\\
&&= C\epsilon^{-2}\int_{\Omega_{L\epsilon}\cap\Gamma_{L_1\epsilon}}|\nabla d_\Gamma(x)\wedge\widehat{n}_\epsilon(x)-\nabla d_{\Gamma}(\Pi_\Gamma(x))\wedge\widehat{n}_\epsilon(\Pi_\Gamma(x))|^2\,dx\\
&&\le C\epsilon^{-2}\int_{\Omega_{L\epsilon}\cap\Gamma_{L_1\epsilon}}
\big(|\nabla d_\Gamma(x)-\nabla d_{\Gamma}(\Pi_\Gamma(x))|^2
+|\widehat{n}_\epsilon(x)-\widehat{n}_\epsilon(\Pi_\Gamma(x))|^2\big)\,dx\\
&&\le C\int_{\Omega_{L\epsilon}\cap\Gamma_{L_1\epsilon}}
\big(|\nabla^2 d_\Gamma(x)|^2+|\nabla\widehat{n}_\epsilon(x)|^2\big)\,dx\\
&&\le C\epsilon+\frac{o_\epsilon(1)}{\epsilon},
\end{eqnarray*}
which yields that
\begin{equation}\label{ext_est14}
\big|\int_\Omega\big(L_2|\nabla \hat{s}_\epsilon\wedge \widehat{n}_\epsilon|^2
+L_4 \hat{s}_\epsilon\nabla \hat{s}_\epsilon\cdot (\nabla \widehat{n}_\epsilon)\widehat{n}_\epsilon\big)\,dx\big|
\le C(1+\epsilon)+\frac{o_\epsilon(1)}{\epsilon}.
\end{equation}
Adding \eqref{ext_est11}, \eqref{ext_est12}, and \eqref{ext_est14} together yields \eqref{SUB1.0}.
This, combined with the minimality of $(s_\epsilon, n_\epsilon)$, implies
\begin{equation}\label{SUB2.0}
\limsup_{\epsilon\to 0}\int_\Omega \epsilon\widetilde{\mathcal{W}}({s}_\epsilon, {n}_\epsilon,
\nabla {s}_\epsilon, \nabla {n}_\epsilon)\,dx
\le \alpha_0\mathcal{H}^2(\Gamma).
\end{equation}
\item [({\bf C})] $L_1=L_2=L_3=L_4=0$. In this case, it is readily seen that
\eqref{SUB2.0} follows directly from \eqref{ext_est11} and \eqref{ext_est12}.
\end{itemize}
Therefore the ``$\le$" part of Theorem \ref{existence} is proven, when $\Omega$ is a bounded domain
in $\R^3$.
\qed
\subsubsection{$\Omega=\R^3$ and $\big|\{x\in\R^3: \ s_\epsilon(x)\ge {s_+}\}\big|=|B_1|$}

The construction for the upper bound estimates for $\Omega=\R^3$ is rather simple. Here we sketch it
as follows.

For a sufficiently large $L>0$, let $\xi_L\in C^\infty\big([{-L}, L]\big)$ be  an almost
$1$-dimensional minimal connecting orbit,
i.e., $\xi_L({-L})=s_+$, $\xi_L({L})=0$, and
\begin{equation}\label{approx_minimal1}
\int_{{-L}}^{{L}} \big(\beta\dot\xi_L^2+W(\xi_L)\big)\,d\tau= \alpha_0+o_L(1).
\end{equation}
Define $\hat{s}_\epsilon:\R^3\to \R_+$ by  letting
\begin{equation*}
\hat{s}_\epsilon(x)=\begin{cases} s_+ & |x|\le 1, \\
\xi_L(\frac{|x|-(1+L\epsilon)}{\epsilon}) & 1\le |x|\le 1+2L\epsilon,\\
0 & |x|\ge 1+2L\epsilon.
\end{cases}
\end{equation*}
Notice that
$$\big\{x\in\R^3:  \hat{s}_\epsilon(x)\ge {s_+}\big\}= B_{1}.$$
Direct calculations imply that
\begin{equation}\label{ext_est15}
\int_{\R^3} \big(|\nabla \hat{s}_\epsilon|^2+\frac{1}{\epsilon^2}W(\hat{s}_\epsilon)\big)\,dx
\le \frac{1}{\epsilon}(\mathcal{H}^2(\Gamma)+C\epsilon)(\alpha_0+o_\epsilon(1)).
\end{equation}

Next, we will construct a  map $\widehat{n}_\epsilon\in H^1(B_1,\mathbb S^2)$ as follows:
\begin{itemize}
\item [({\bf A})] If $L_1>L_2=L_4=0$, then it is well-known that there exists a map $\widehat{n}_\epsilon:
B_{1+2L\epsilon}\to\mathbb S^2$ such that $\widehat{n}_\epsilon(x)\cdot x=0$ on $\partial B_1$, and
$$
\int_{B_{1+2L\epsilon}}|\nabla \widehat{n}_\epsilon|^2\,dx\le C(L).
$$
Hence
\begin{equation}\label{ext_est16}
\int_{\R^3} \hat{s}_\epsilon^2 \big(W_{OF}(\widehat{n}_\epsilon, \nabla \widehat{n}_\epsilon)+\alpha |\nabla\widehat{n}_\epsilon|^2\big)\,dx\le C(L)s_+^2.
\end{equation}
While \begin{eqnarray*}
&&\big|\int_{\R^3} \big(L_1(\nabla\hat{s}_\epsilon\cdot \widehat{n}_\epsilon)^2
+L_3(\nabla\hat{s}_\epsilon\cdot \widehat{n}_\epsilon)(\hat{s}_\epsilon{\rm{div}}\widehat{n}_\epsilon)\big)\,dx\big|\\
&&\le C\big(\int_{B_{1+2L\epsilon}\setminus B_1}|\nabla \hat{s}_\epsilon\cdot\widehat{n}_\epsilon|^2\,dx
+\int_{B_{1+2L\epsilon}} |\nabla\widehat{n}_\epsilon|^2\,dx\big)\\
&&\le C(L)+C\epsilon^{-2}\int_{B_{1+2L\epsilon}\setminus B_1}
\big|\widehat{n}_\epsilon(x)-\widehat{n}_\epsilon(\frac{x}{|x|})\big|^2\,dx\\
&&\le C(L)+C\int_{B_{1+2L\epsilon}\setminus B_1}\big|\nabla\widehat{n}_\epsilon\big|^2\,dx\le C(L).
\end{eqnarray*}
Therefore we arrive at
\begin{eqnarray}\label{SUB2.1}
\limsup_{\epsilon\to 0}\int_{\R^3} \epsilon\widetilde{\mathcal{W}}({s}_\epsilon, {n}_\epsilon,
\nabla {s}_\epsilon, \nabla {n}_\epsilon)\,dx&\le&
\limsup_{\epsilon\to 0}\int_{\R^3} \epsilon\widetilde{\mathcal{W}}(\hat{s}_\epsilon, \widehat{n}_\epsilon,
\nabla \hat{s}_\epsilon, \nabla \widehat{n}_\epsilon)\,dx\nonumber\\
&\le& \alpha_0\mathcal{H}^2(\partial B_1).
\end{eqnarray}
\item [({\bf B})] If $L_2>L_1=L_3=0$, then we simply let
$$\widehat{n}_\epsilon(x)=\frac{x}{|x|}, \ x\in B_{1+2L\epsilon}.$$
Then it is straightforward to check that
\begin{equation}\label{ext_est17}
\int_{\R^3} \hat{s}_\epsilon^2 \big(W_{OF}(\widehat{n}_\epsilon, \nabla \widehat{n}_\epsilon)+\alpha |\nabla\widehat{n}_\epsilon|^2\big)\,dx\le Cs_+^2 L.
\end{equation}
While
\begin{eqnarray}\label{ext_est18}
\big|\int_{\R^3}\big(L_2|\nabla \hat{s}_\epsilon\wedge \widehat{n}_\epsilon|^2
+L_4 \hat{s}_\epsilon\nabla \hat{s}_\epsilon\cdot (\nabla \widehat{n}_\epsilon)\widehat{n}_\epsilon\big)\,dx\big|
\le C\epsilon^{-2}\int_{B_{1+2L\epsilon}\setminus B_1} |\frac{x}{|x|}\wedge \frac{x}{|x|}|^2\,dx=0.
\end{eqnarray}
Hence \eqref{SUB2.1} holds.
\item [({\bf C})] If $L_1=L_2=L_3=L_4=0$, then we simply set $\widehat{n}_\epsilon\equiv (0,0,1)\in \mathbb S^2$.
In this case, it is easy to see that
\begin{eqnarray*}
\int_{\R^3}\epsilon \widetilde{\mathcal{W}}(\hat{s}_\epsilon, \widehat{n}_\epsilon,
\nabla \hat{s}_\epsilon, \nabla \widehat{n}_\epsilon)\,dx
&=&\int_{B_{1+2L\epsilon}\setminus B_1}
\big(\epsilon|\nabla \hat{s}_\epsilon|^2+\frac{1}{\epsilon} W(\hat{s}_\epsilon)\big)\,dx\\
&\le &\alpha_0 (\mathcal{H}^2(\partial B_1)+C\epsilon).
\end{eqnarray*}
Hence  \eqref{SUB2.1} holds.
\end{itemize}
The ``$\le$" part of Theorem \ref{existence} is proven, when $\Omega=\R^3$.
Combining these two subsections, we prove Theorem \ref{existence}.
\qed

\section{Proof of Theorem \ref{sharp}}
\setcounter{equation}{0}
\setcounter{theorem}{1}
This section is devoted to the proof of Theorem \ref{sharp} for the case that $\Omega$ is a bounded
smooth domain in $\R^3$. It involves refined estimates of both upper bounds and lower bounds
of the total energy $\mathcal{E}(s_\epsilon, n_\epsilon)$ for minimizers $(s_\epsilon, n_\epsilon)$,
in which the strict stability of $\Gamma$
plays a crucial role in a perfect matching of the coefficients of leading order term or $O(\frac{1}{\epsilon})$ term in the expansion of $\mathcal{E}(s_\epsilon, n_\epsilon)$.

\subsection{Refined energy upper bounds}
\setcounter{equation}{0}
\setcounter{theorem}{0}

In this subsection, we will  prove an optimal upper bound for the energy
$\mathcal{E}(s_\epsilon, n_\epsilon)$ of minimizers $(s_\epsilon, n_\epsilon)$.
This is done by utilizing the additional assumption on the boundary value
$(t_\epsilon, g_\epsilon)$ to construct  a comparison map such that
$s_\epsilon$ is approximately a minimal connecting orbit in the transition
region of $\Gamma$ of width of $O(\epsilon^\gamma)$, and $n_\epsilon$ is
approximately a minimizing harmonic map in the corresponding configuration
spaces in $\Omega_+$.

We divide the estimates of refined upper bounds for the cases (A), (B), and (C) into three  separate Lemmas.

\begin{lemma}\label{Upper_Est1} Assume $\Gamma$ and the boundary values
$(t_\epsilon, g_\epsilon)$ satisfy the same assumptions as in Theorem \ref{sharp}.
If $L_1>L_2=L_4=0$, then
\begin{eqnarray}\label{A}
&&\inf\Big\{\int_\Omega \mathcal{W}_\epsilon(s,n,\nabla s, \nabla n)\,dx\ \big| \ (s_\epsilon, n_\epsilon)=(t_\epsilon, g_\epsilon) \ {\rm{on}}\ \partial\Omega\Big\}\nonumber\\
&&\le \frac{1}{\epsilon}\int_{0}^{s_+}
\sqrt{2\beta W(\tau)} \mathcal{H}^2(\Gamma(\epsilon\xi_{\epsilon,\gamma}^{-1}(\tau)) \,d\tau+\mathcal{D}_A+o_\epsilon(1),
\end{eqnarray}
where $\mathcal{D}_A$ is given by \eqref{OFM1} and \eqref{bdry1} of Theorem \ref{sharp}.
\end{lemma}
\pf We first construct an extension of $s_\epsilon$ from $\partial\Omega$ to $\Omega$ as follows.
For $\gamma\in (0,1)$, let $\xi_{\epsilon,\gamma}\in C^\infty([-\epsilon^{\gamma-1},\epsilon^{\gamma-1}], \R)$ be given
by \eqref{xi_cond} and $t_\epsilon:\partial\Omega\to\R$ satisfy  \eqref{b_cond1.3}.
Define $s_\epsilon$ in the fast transition region $\Gamma_{\epsilon^\gamma}$ by letting
$$s_\epsilon(x)=\xi_{\epsilon,\gamma}\big(\frac{d_\Gamma(x)}{\epsilon}\big), \ \forall x\in \Gamma_{\epsilon^\gamma}.$$
In the off-transition region $\Omega^+\setminus \Gamma_{\epsilon^\gamma}$, we perform a linear interpolation between
$s_\epsilon$ and $s_+$ in a $\epsilon^\gamma$-neighborhood of $\Sigma^+$. More precisely,
decompose
\begin{eqnarray*}\Omega^\pm\setminus \Gamma_{\epsilon^\gamma}&=&
\Big((\Omega^\pm\setminus\Gamma_{\epsilon^\gamma})\cap \{x\in\Omega: d(x,\Sigma^\pm)<\epsilon^\gamma\}\Big)
\cup \Big((\Omega^\pm\setminus\Gamma_{\epsilon^\gamma})\setminus \{x\in\Omega: d(x,\Sigma^\pm)<\epsilon^\gamma\}\Big)\\
&=&E_{\epsilon^\gamma}^\pm\cup F_{\epsilon^\gamma}^\pm.
\end{eqnarray*}
Define
$$s_\epsilon(x)=\begin{cases} s_+, & x\in F_{\epsilon^\gamma}^+,\\
\mbox{linear interpolation of } \ s_\epsilon\big|_{\Sigma^+\setminus\Gamma_{\epsilon^\gamma}} \ \mbox{and}\  s_+\big|
_{\{x\in\Omega^+\setminus\Gamma_{\epsilon^\gamma}: d(x,\Sigma^+)=\epsilon^\gamma\}},  & x\in E_{\epsilon^\gamma}^+.
\end{cases}
$$
Similarly, define
$$s_\epsilon(x)=\begin{cases} 0, & x\in F_{\epsilon^\gamma}^-,\\
\mbox{linear interpolation of } \ s_\epsilon\big|_{\Sigma^-\setminus\Gamma_{\epsilon^\gamma}} \ \mbox{and}\  0\big|
_{\{x\in\Omega^-\setminus\Gamma_{\epsilon^\gamma}: d(x,\Sigma^-)=\epsilon^\gamma\}},  & x\in E_{\epsilon^\gamma}^-.
\end{cases}
$$
Then by the co-area formula and direct calculations (see Maggi \cite{Maggi2}) we can estimate
\begin{eqnarray}\label{s-UB1}
&&\int_\Omega \big(\beta |\nabla s_\epsilon|^2+\frac{1}{\epsilon^2}W(s_\epsilon)\big)\,dx\nonumber\\
&&=\Big\{\int_{\Gamma_{\epsilon^\gamma}}+\int_{E_{\epsilon^\gamma}^-}+\int_{E_{\epsilon^\gamma}^+}\Big\}
 \big(\beta |\nabla s_\epsilon|^2+\frac{1}{\epsilon^2}W(s_\epsilon)\big)\,dx\nonumber\\
 &&\le\frac{1}{\epsilon}\int_{-\epsilon^{\gamma-1}}^{\epsilon^{\gamma-1}}\big(\beta |\xi'_{\epsilon,\gamma}(t)|^2+W(\xi_{\epsilon,\gamma}(t))\big) \mathcal{H}^2(\Gamma(\epsilon t)) \,dt
 +C\epsilon^\gamma \int_{\partial\Omega\setminus\Gamma_{\epsilon^\gamma}} |\nabla_{\rm{tan}}t_\epsilon|^2\,d\mathcal{H}^2\nonumber\\
 &&\quad+C\epsilon^{-\gamma}\big(\int_{\Sigma^+\setminus\Gamma_{\epsilon^\gamma}}
 |t_\epsilon-s_+|^2\,d\mathcal{H}^2
 +\int_{\Sigma^-\setminus\Gamma_{\epsilon^\gamma}}|t_\epsilon|^2\,d\mathcal{H}^2\big)\nonumber\\
 &&\quad+C\int_{E_{\epsilon^\gamma}^+\cup E_{\epsilon^\gamma}^-} \frac{1}{\epsilon^2} W(s_\epsilon)\,dx\nonumber\\
 &&\le \frac{1}{\epsilon}\int_{-\epsilon^{\gamma-1}}^{\epsilon^{\gamma-1}}\big(\beta |\xi'_{\epsilon,\gamma}(t)|^2+W(\xi_{\epsilon,\gamma}(t))\big) \mathcal{H}^2(\Gamma(\epsilon t)) \,dt
 +C\epsilon^\gamma+o_\epsilon(1).
\end{eqnarray}
Here we have applied \eqref{W-cond} and \eqref{b_cond1.3} in the last step, which ensures
\begin{eqnarray*}
\int_{E_{\epsilon^\gamma}^+\cup E_{\epsilon^\gamma}^-} \frac{1}{\epsilon^2} W(s_\epsilon)\,dx
&\le& C\epsilon^\gamma\Big(\int_{\Sigma^+\setminus\Gamma_{\epsilon^\gamma}}\frac{|t_\epsilon-s_+|^2}{\epsilon^2}\,d\mathcal{H}^2+\int_{\Sigma^-\setminus\Gamma_{\epsilon^\gamma}}\frac{|t_\epsilon|^2}{\epsilon^2}\,d\mathcal{H}^2\Big)\\
&\le& C\epsilon^\gamma.
\end{eqnarray*}
Notice that the condition \eqref{xi_cond} implies that
\begin{eqnarray}\label{sharp_path}
&&\frac{1}{\epsilon}\int_{-\epsilon^{\gamma-1}}^{\epsilon^{\gamma-1}}\big(\beta |\xi'_{\epsilon,\gamma}(t)|^2+W(\xi_{\epsilon,\gamma}(t))\big) \mathcal{H}^2(\Gamma(\epsilon t)) \,dt\nonumber\\
&&\le \frac{1}{\epsilon}\int_{-\epsilon^{\gamma-1}}^{\epsilon^{\gamma-1}}
\sqrt{2\beta W(\xi_{\epsilon,\gamma}(t))}|\xi'_{\epsilon,\gamma}(t)| \mathcal{H}^2(\Gamma(\epsilon t)) \,dt
+C_2\epsilon^{\gamma-2}e^{-C_1\epsilon^{\gamma-1}}\nonumber\\
&&\le \frac{1}{\epsilon}\int_{0}^{s_+}
\sqrt{2\beta W(\tau)} \mathcal{H}^2(\Gamma(\epsilon\xi_{\epsilon,\gamma}^{-1}(\tau)) \,d\tau
+o_{\epsilon}(1).
\end{eqnarray}

Next we want to construct an extension map $n_\epsilon:\Omega\to\mathbb S^2$ from $g_\epsilon:\partial\Omega\to\mathbb S^2$.  To do it, let $n\in H^1(\Omega^+, \mathbb S^2)$ achieve
$\mathcal{D}_A$, i.e.,
$n=g$ on $\Sigma^+$ and $n\cdot\nu_\Gamma=0$ on $\Gamma$, and
$$
E(n; \Omega^+)=\mathcal{D}_A.
$$
Recall that there exists a map $\widehat{n}_\epsilon$
in  the region $\Omega^+\setminus\big(\Gamma_{2\epsilon^\gamma}\cup \Omega_{\epsilon^\gamma}^+\big)$
that is a linear interpolation between $g_\epsilon\big|_{\Sigma^+\setminus\Gamma_{2\epsilon^\gamma}}$
and $n\big|_{\partial\Omega_{\epsilon^\gamma}^+\setminus\Gamma_{2\epsilon^\gamma}}$, i.e.,
$$
\widehat{n}_\epsilon(r,\theta)=\frac{r}{\epsilon^\gamma}g_\epsilon(\theta)+\frac{\epsilon^\gamma-r}{\epsilon^\gamma} n(\epsilon^\gamma, \theta),
\ \ x=(r,\theta)\in (\Sigma^+\setminus\Gamma_{2\epsilon^\gamma})\times [0,\epsilon^\gamma]\approx
\Omega^+\setminus\big(\Gamma_{2\epsilon^\gamma}\cup \Omega_{\epsilon^\gamma}^+\big).
$$
Since $\widehat{n}_\epsilon$ may not map into $\mathbb S^2$, we need to apply Hardt-Lin's extension Lemma
to find a point $a\in \R^3$, with $|a|\le \frac12$, such that the map
$\Psi_a=\big(\Pi_a\big|_{\mathbb S^2}\big)^{-1}\circ\Pi_a$, with $\Pi_a(y)=\frac{y-a}{|y-a|}:\R^3\to\mathbb S^2$, satisfies
$$
\int_{\Omega^+\setminus(\Gamma_{2\epsilon^\gamma}\cup\Omega^+_{\epsilon^\gamma})} s_\epsilon^2\big|\nabla (\Psi_a(\widehat{n}_\epsilon))\big|^2\,dx
\le C\int_{\Omega^+\setminus(\Gamma_{2\epsilon^\gamma}\cup\Omega^+_{\epsilon^\gamma})} s_\epsilon^2\big|\nabla \widehat{n}_\epsilon\big|^2\,dx.
$$
Now we define ${n}_\epsilon:\Omega^+\to\mathbb{S}^2$ as follows. First, we define
\begin{equation*}
{n}_{\epsilon}(x)=\begin{cases} n(x) & x\in (\Omega^+\cap\Gamma_{\epsilon^\gamma})\cup\Omega^+_{\epsilon^\gamma},\\
\Psi_a\big(\widehat{n}_\epsilon(x)\big) & x\in \Omega^+\setminus \big(\Gamma_{2\epsilon^\gamma}
\cup\Omega^+_{\epsilon^\gamma}\big).
\end{cases}
\end{equation*}

Since $(\Omega^+\setminus \Omega^+_{\epsilon^\gamma})\cap (\Gamma_{2\epsilon^\gamma}\setminus \Gamma_{\epsilon^\gamma})$
can be identified as a ball of radius $\epsilon^\gamma$, centered at a point $x_*$, we can define
$n_\epsilon: (\Omega^+\setminus \Omega^+_{\epsilon^\gamma})\cap (\Gamma_{2\epsilon^\gamma}\setminus \Gamma_{\epsilon^\gamma})\to\mathbb S^2$ as the homogeneous degree zero extension, with respect to $x_*$,
of the value of $n_\epsilon$ on $\partial \big((\Omega^+\setminus \Omega^+_{\epsilon^\gamma})\cap (\Gamma_{2\epsilon^\gamma}\setminus \Gamma_{\epsilon^\gamma})\big)$.

Then we can calculate
\begin{eqnarray}\label{n-UB1.1}
&&\int_{(\Omega^+\cap\Gamma_{\epsilon^\gamma})\cup\Omega^+_{\epsilon^\gamma}}\big(s_\epsilon^2W_{OF}(n_\epsilon, \nabla n_\epsilon)
+\alpha s_\epsilon^2|\nabla{n}_\epsilon|^2\big)\,dx\nonumber\\
&&\le (1+o_\epsilon(1)) s_+^2\int_{\Omega^+}(W_{OF}(n,\nabla n)+\alpha|\nabla n|^2)\,dx
\le (1+o_\epsilon(1))\mathcal{D}_A,
\end{eqnarray}
and
\begin{eqnarray}\label{n-UB1.2}
&&\int_{\Omega^+\setminus(\Gamma_{2\epsilon^\gamma}\cup\Omega^+_{\epsilon^\gamma})}\big(s_\epsilon^2W_{OF}(n_\epsilon, \nabla n_\epsilon)
+\alpha s_\epsilon^2|\nabla{n}_\epsilon|^2\big)\,dx
\le C\int_{\Omega^+\setminus(\Gamma_{2\epsilon^\gamma}\cup\Omega^+_{\epsilon^\gamma})}
s_\epsilon^2|\nabla n_\epsilon|^2\,dx\nonumber\\
&&\le Cs_+^2\epsilon^\gamma\big(\int_{\Sigma^+}|\nabla_{\rm{tan}}g_\epsilon|^2\,d\mathcal{H}^2
+\int_{\partial\Omega^+_{\epsilon^\gamma}\cap\Omega^+}|\nabla_{\rm{tan}}n|^2\,d\mathcal{H}^2\big)
+\frac{Cs_+^2}{\epsilon^\gamma}\int_{\Sigma^+\setminus\Gamma_{\epsilon^\gamma}} |g_\epsilon(\theta)-n(\epsilon^\gamma, \theta)|^2\,d\mathcal{H}^2\nonumber\\
&&\le C\epsilon^\gamma+C\int_{\Omega^+\setminus \Omega^+_{2\epsilon^\gamma}}|\nabla n|^2
+\frac{C}{\epsilon^\gamma}\int_{\Sigma^+\setminus\Gamma_{\epsilon^\gamma}}|g_\epsilon-g|^2\,d\mathcal{H}^2
+C\int_{\Sigma^+}|g(\theta)-n(\epsilon^\gamma, \theta)|^2\,d\mathcal{H}^2\nonumber\\
&&\le C\epsilon^\gamma+o_\epsilon(1)+C\int_{\Omega^+\setminus \Omega^+_{2\epsilon^\gamma}}|\nabla n|^2\le C\epsilon^\gamma+o_\epsilon(1).
\end{eqnarray}
where we have used \eqref{b_cond1.4} and the absolute continuity of $\displaystyle\int |\nabla n|^2\,dx$, and
the inequality
$$
\int_{\Sigma^+}|g(\theta)-n(\epsilon^\gamma, \theta)|^2\,d\mathcal{H}^2
\le\epsilon^\gamma\int_{\Omega^+\setminus\Omega^+_{\epsilon^\gamma}}|\nabla n|^2\,dx\le C\epsilon^\gamma.
$$
.

While
\begin{eqnarray}\label{n-UB1.20}
&&\int_{(\Omega^+\setminus \Omega^+_{\epsilon^\gamma})\cap (\Gamma_{2\epsilon^\gamma}\setminus \Gamma_{\epsilon^\gamma})}
\big(s_\epsilon^2W_{OF}(n_\epsilon, \nabla n_\epsilon)
+\alpha s_\epsilon^2|\nabla{n}_\epsilon|^2\big)\,dx\nonumber\\
&&\le Cs_+^2\int_{(\Omega^+\setminus \Omega^+_{\epsilon^\gamma})\cap (\Gamma_{2\epsilon^\gamma}\setminus \Gamma_{\epsilon^\gamma})}|\nabla{n}_\epsilon|^2\,dx\nonumber\\
&&\le C\epsilon^\gamma\int_{\partial\big((\Omega^+\setminus \Omega^+_{\epsilon^\gamma})\cap (\Gamma_{2\epsilon^\gamma}\setminus \Gamma_{\epsilon^\gamma})\big)} |\nabla_{\rm{tan}} n_\epsilon|^2\,d\mathcal{H}^2\nonumber\\
&&\le C\epsilon^\gamma\Big\{\int_{\Sigma^+\cap\Gamma_{2\epsilon^\gamma}}
+\int_{\partial\Gamma_{\epsilon^\gamma}\cap(\Omega^+\setminus\Omega^+_{\epsilon^\gamma})}
+\int_{\partial\Omega^+_{\epsilon^\gamma}\cap(\Gamma_{2\epsilon^\gamma}\setminus\Gamma_{\epsilon^\gamma})}
+\int_{\partial\Gamma_{2\epsilon^\gamma}\cap(\Omega^+\setminus\Omega^+_{\epsilon^\gamma})}\Big\}
|\nabla_{\rm{tan}} n_\epsilon|^2\,d\mathcal{H}^2\nonumber\\
&&\le C\Big(\epsilon^\gamma\int_{\Sigma^+}|\nabla g_\epsilon|^2\,d\mathcal{H}^2+\int_{\Omega^+\setminus\Omega^+_{\epsilon^\gamma}}|\nabla n|^2\,dx
+\int_{\Omega^+\setminus(\Gamma_{2\epsilon^\gamma}\cup\Omega^+_{\epsilon^\gamma})}|\nabla n_\epsilon|^2\,dx
+\int_{\Omega^+_{\epsilon^\gamma}\setminus\Omega^+_{2\epsilon^\gamma}}|\nabla n|^2\,dx\Big)\nonumber\\
&&\le C(\epsilon^\gamma+o_\epsilon(1)+\int_{\Omega^+\setminus\Omega^+_{2\epsilon^\gamma}}|\nabla n|^2\,dx)
\le C(\epsilon^\gamma+o_\epsilon(1)),
\end{eqnarray}
where we have used \eqref{n-UB1.20} in the last step, and we have applied Fubini's theorem
which guarantees the following inequalities:
$$
\epsilon^\gamma\int_{\partial\Gamma_{\epsilon^\gamma}\cap(\Omega^+\setminus\Omega^+_{\epsilon^\gamma})}
|\nabla_{\rm{tan}}n_\epsilon|^2\,d\mathcal{H}^2
\le C\int_{\Omega^+\setminus\Omega_{\epsilon^\gamma}^+}|\nabla n|^2\,dx=o_\epsilon(1),
$$
$$
\epsilon^\gamma\int_{\partial\Gamma_{2\epsilon^\gamma}\cap(\Omega^+\setminus\Omega^+_{\epsilon^\gamma})}
|\nabla_{\rm{tan}}n_\epsilon|^2\,d\mathcal{H}^2\le C
\int_{\Omega^+\setminus(\Gamma_{2\epsilon^\gamma}\cup\Omega^+_{\epsilon^\gamma})}|\nabla n_\epsilon|^2\,dx
\le C(\epsilon^\gamma+o_\epsilon(1)),
$$
and
$$
\epsilon^\gamma\int_{\partial\Omega^+_{\epsilon^\gamma}\cap(\Gamma_{2\epsilon^\gamma}\setminus\Gamma_{\epsilon^\gamma})}|\nabla_{\rm{tan}}n_\epsilon|^2\,d\mathcal{H}^2
\le C\int_{\Omega^+\setminus\Omega^+_{2\epsilon^\gamma}}|\nabla n|^2\,dx=o_\epsilon(1).
$$

Hence
\begin{eqnarray}\label{n-UB1.3}
\int_{\Omega^+}\big(s_\epsilon^2W_{OF}(n_\epsilon, \nabla n_\epsilon)
+\alpha s_\epsilon^2|\nabla{n}_\epsilon|^2\big)\,dx\le (1+o_\epsilon(1))\mathcal{D}_A+C\big(\epsilon^\gamma+o_\epsilon(1)\big).
\end{eqnarray}
The most difficult term to estimate is the interactive energy between $\nabla s_\epsilon$ and $n_\epsilon$. To do it,
we proceed as follows.
\begin{eqnarray*}\label{n-UB1.4}
&&\int_{\Omega^+}\big(L_1|\nabla s_\epsilon\cdot n_\epsilon|^2
+L_3(\nabla s_\epsilon\cdot n_\epsilon)(s_\epsilon{\rm{div}}n_\epsilon)\big)\,dx\nonumber\\
&&\le
\Big\{\int_{\Omega^+\cap\Gamma_{\epsilon^\gamma}}+\int_{\Omega^+\setminus(\Omega_{\epsilon^\gamma}^+\cup\Gamma_{\epsilon^\gamma})}\Big\} \big(L_1|\nabla s_\epsilon\cdot n_\epsilon|^2
+L_3(\nabla s_\epsilon\cdot n_\epsilon)(s_\epsilon{\rm{div}}n_\epsilon\big)\,dx\nonumber\\
&&\le C\Big(\int_{\Omega^+\cap\Gamma_{\epsilon^\gamma}}|\nabla s_\epsilon\cdot n_\epsilon|^2
+\big(\int_{\Omega^+\cap\Gamma_{\epsilon^\gamma}}|\nabla s_\epsilon\cdot n_\epsilon|^2\big)^\frac12\Big)\nonumber\\
&&\quad+C\Big(\int_{E_{\epsilon^\gamma}^+}|\nabla s_\epsilon|^2\,dx+\big(\int_{E_{\epsilon^\gamma}^+}|\nabla s_\epsilon|^2\,dx\big)^\frac12\Big)=I+II.
\end{eqnarray*}
From the estimate \eqref{s-UB1}, we can see that
$$\int_{E_{\epsilon^\gamma}^+}|\nabla s_\epsilon|^2\,dx\le C(\epsilon^\gamma+o_\epsilon(1))$$
so that
$$II\le C(\epsilon^{\frac{\gamma}2}+o_\epsilon(1)).$$
To estimate $I$, let $\Pi_\Gamma: \Gamma_{\epsilon^\gamma}\to\Gamma$ be the smooth nearest point
projection map. Since $\nabla d_{\Gamma}(x)=\nu_{\Gamma}(x)$ for $x\in\Gamma$, we have that
$\nabla d_{\Gamma}(x)\cdot n(x)=0$ for $x\in\Gamma$. Hence
\begin{eqnarray*}
&&\int_{\Omega^+\cap\Gamma_{\epsilon^\gamma}}|\nabla s_\epsilon\cdot n_\epsilon|^2\,dx
\le \frac{1}{\epsilon^2}\int_{\Omega^+\cap\Gamma_{\epsilon^\gamma}}
(\xi'_{\epsilon,\gamma})^2(\frac{d_\Gamma(x)}{\epsilon})|\nabla d_{\Gamma}(x)\cdot n(x)|^2\,dx\\
&&\le \frac{1}{\epsilon^2}\int_{\Omega^+\cap\Gamma_{\epsilon^\gamma}}(\xi'_{\epsilon, \gamma})^2(\frac{d_\Gamma(x)}{\epsilon}) |\nabla d_{\Gamma}(x)\cdot n(x)
-\nabla d_{\Gamma}(\Pi_{\Gamma}(x))\cdot n(\Pi_{\Gamma}(x))|^2\,dx\\
&&= \frac{1}{\epsilon^2}\big\{\int_{\Omega^+\cap\Gamma_{L\epsilon}}+\int_{\Omega^+\cap(\Gamma_{\epsilon^\gamma}\setminus\Gamma_{L\epsilon})}\big\}
(\xi'_{\epsilon,\gamma})^2(\frac{d_\Gamma(x)}{\epsilon})\\
&&\qquad\qquad\qquad\qquad\qquad\cdot\big(|\nabla d_{\Gamma}(x)-\nabla d_{\Gamma}(\Pi_\Gamma(x))|^2
+|n(x)-n(\Pi_\Gamma(x))|^2\big)\,dx\\
&&=III+IV.
\end{eqnarray*}
$III$ can be estimated similarly to \eqref{ext_est19} so that
$$III\le C\epsilon+o_\epsilon(1).$$
While we can utilize the decay property of $|\xi'|$ to estimate $IV$ as follows.
\begin{eqnarray*}
IV&\le&
\frac{C}{\epsilon^2}\int_{\Omega^+\cap(\Gamma_{\epsilon^\gamma}\setminus\Gamma_{L\epsilon})}
e^{-\frac{C}{\epsilon}d_\Gamma(x)}\big(|\nabla d_{\Gamma}(x)-\nabla d_{\Gamma}(\Pi_\Gamma(x))|^2
+|n(x)-n(\Pi_\Gamma(x))|^2\big)\,dx\\
&\le& C\big\|\nabla^2d_{\Gamma}\big\|_{\Omega^+\cap\Gamma_{\epsilon^\gamma}} \epsilon^\gamma
+C\epsilon^{-2} \int_{\Omega^+\cap\Gamma_{\epsilon^\gamma}}
e^{-\frac{C}{\epsilon}d_\Gamma(x)}|n(x)-n(\Pi_\Gamma(x))|^2\,dx.
\end{eqnarray*}
Notice that by identifying $\Omega^+\cap\Gamma_{\epsilon^\gamma}$ with $\Gamma\times [0,\epsilon^\gamma]$,
we can bound
\begin{eqnarray*}
&&\int_{\Omega^+\cap\Gamma_{\epsilon^\gamma}}
e^{-\frac{C}{\epsilon}d_\Gamma(x)}|n(x)-n(\Pi_\Gamma(x))|^2\,dx\\
&&\le C\int_{\Gamma} \int_0^{\epsilon^\gamma} e^{-\frac{Ct}{\epsilon}} |n(t,\theta)-n(0,\theta)|^2\,dtd\mathcal{H}^2\\
&&\le C\int_{\Gamma} \int_0^{\epsilon^\gamma} te^{-\frac{Ct}{\epsilon}} (\int_0^{\epsilon^\gamma} |n_t|^2(\tau,\theta)\,d\tau)\,dtd\mathcal{H}^2\\
&&\le C\big(\int_0^{\epsilon^\gamma} t e^{-\frac{Ct}{\epsilon}}\,dt\big)\int_{\Omega^+\cap\Gamma_{\epsilon^\gamma}} |\nabla n|^2\,dx\\
&&\le C\epsilon^2\big(\int_0^{\infty} t e^{-{Ct}}\,dt\big)\int_{\Omega^+\cap\Gamma_{\epsilon^\gamma}} |\nabla n|^2\,dx\\
&&\le C\epsilon^2 \int_{\Omega^+\cap\Gamma_{\epsilon^\gamma}} |\nabla n|^2\,dx.
\end{eqnarray*}
Therefore we obtain that
$$
IV\le C\epsilon^\gamma+C \int_{\Omega^+\cap\Gamma_{\epsilon^\gamma}} |\nabla n|^2\,dx
\le C\epsilon^\gamma+o_\epsilon(1).
$$
From the estimates of $III$ and $IV$, we obtain that
\begin{equation}\label{n-UB1.5}
\int_{\Omega^+}\big(L_1|\nabla s_\epsilon\cdot n_\epsilon|^2
+L_3(\nabla s_\epsilon\cdot n_\epsilon)(s_\epsilon{\rm{div}}n_\epsilon)\big)\,dx\le
C\big(\epsilon^{\frac{\gamma}{2}}+o_\epsilon(1)\big).
\end{equation}

Next we want to construct a map $n_\epsilon:\Omega^{-}\to\mathbb S^2$ such that
it enjoys a upper bound estimate similar to that in $\Omega^+$.
First let $\widehat{n}_\epsilon:\Omega^-\to\R^3$ be such that
$$
\begin{cases}
\Delta \widehat{n}_\epsilon=0 &\ \ {\rm{in}}\ \ \Omega^-, \\
\widehat{n}_\epsilon=g_\epsilon &\ \ {\rm{on}}\ \Sigma^-,\\
\widehat{n}_\epsilon=n & \ \ {\rm{on}}\ \ \Gamma.
\end{cases}
$$
Then it is well-known that
\begin{eqnarray*}
\int_{\Omega^-}|\nabla \widehat{n}_\epsilon|^2\,dx\le C\big\|\widehat{n}_\epsilon\big\|_{H^\frac12(\partial\Omega^-)}^2
&\le& C\big(\int_{\Sigma^-} |\nabla g_\epsilon|^2\,d\mathcal{H}^2+\int_{\Omega^+} |\nabla n|^2\,dx\big)\\
&\le& C\big(1+\int_{\Sigma^-} |\nabla g|^2\,d\mathcal{H}^2+\int_{\Omega^+} |\nabla n|^2\,dx\big).
\end{eqnarray*}
Moreover, since $g_\epsilon\rightarrow g$ in $H^1(\Sigma^-)$,  there exists
$\widetilde{n}\in H^1(\Omega^-)$ such that
\begin{equation}\label{strong-h1}
\widehat{n}_\epsilon\rightarrow \widetilde{n} \ \ {\rm{in}}\ \ H^1(\Omega^-).
\end{equation}
Applying Hardt-Lin's extension Lemma again, there exists $a\in \R^3$ with $|a|\le \frac12$ such that
for $n_\epsilon=\Psi_a(\widehat{n}_\epsilon)$, where $\Psi_a=\big(\Pi_a\big|_{\mathbb S^2}\big)^{-1}\circ\Pi_a$, with $\Pi_a(y)=\frac{y-a}{|y-a|}:\R^3\to\mathbb S^2$, satisfies
\begin{equation*}
\int_{\Omega^-}|\nabla {n}_\epsilon|^2\,dx
\le C\int_{\Omega^-}|\nabla \hat{n}_\epsilon|^2\,dx\le C\big(1+\int_{\Sigma^-} |\nabla g|^2\,d\mathcal{H}^2+\int_{\Omega^+} |\nabla n|^2\,dx\big),
\end{equation*}
and
\begin{equation}\label{strong-h2}
{n}_\epsilon\rightarrow \Psi_a(\widetilde{n}) \ \ {\rm{in}}\ \ H^1(\Omega^-).
\end{equation}
With the help of \eqref{strong-h2}, we can estimate
\begin{eqnarray}\label{n-UB1.6}
&&\int_{\Omega^-} \big(s_\epsilon^2W_{OF}({n}_\epsilon, \nabla {n}_\epsilon)+\alpha s_\epsilon^2|\nabla {n}_\epsilon|^2\big)\,dx\le C\int_{\Omega^-}s_\epsilon^2|\nabla {n}_\epsilon|^2\,dx\nonumber\\
&&\le C\big(\int_{\Omega^-\cap\Gamma_{\epsilon^\epsilon}}s_\epsilon^2|\nabla n_{\epsilon}|^2\,dx
+\int_{E_{\epsilon^\gamma}^-}s_\epsilon^2|\nabla n_{\epsilon}|^2\,dx\big)\nonumber\\
&&\le  C\big(\int_{\Omega^-\cap\Gamma_{\epsilon^\gamma}}|\nabla n_{\epsilon}|^2\,dx
+\int_{E_{\epsilon^\gamma}^-}|\nabla n_{\epsilon}|^2\,dx\big)\nonumber\\
&&\le C\int_{\Omega^-}|\nabla (n_{\epsilon}-\Phi_a(\widetilde{n}))|^2\,dx
+C\big(\int_{\Omega^-\cap\Gamma_{\epsilon^\gamma}}|\nabla \widetilde{n}|^2\,dx
+\int_{E_{\epsilon^\gamma}^-}|\nabla \widetilde{n}|^2\,dx\big)\nonumber\\
&&=o_\epsilon(1).
\end{eqnarray}
While
\begin{eqnarray}\label{n-UB1.7}
&&\int_{\Omega^-}\big(L_1|\nabla s_\epsilon\cdot n_\epsilon|^2
+L_3(\nabla s_\epsilon\cdot n_\epsilon)(s_\epsilon{\rm{div}}n_\epsilon)\big)\,dx\nonumber\\
&&\le
\Big\{\int_{\Omega^-\cap\Gamma_{\epsilon^\gamma}}+\int_{\Omega^-\setminus(\Omega_{\epsilon^\gamma}^-\cup\Gamma_{\epsilon^\gamma})}\Big\} \big(L_1|\nabla s_\epsilon\cdot n_\epsilon|^2
+L_3(\nabla s_\epsilon\cdot n_\epsilon)(s_\epsilon{\rm{div}}n_\epsilon)\big)\,dx\nonumber\\
&&\le C\Big(\int_{\Omega^-\cap\Gamma_{\epsilon^\gamma}}|\nabla s_\epsilon\cdot n_\epsilon|^2
+\big(\int_{\Omega^-\cap\Gamma_{\epsilon^\gamma}}|\nabla s_\epsilon\cdot n_\epsilon|^2\big)^\frac12\Big)\nonumber\\
&&\quad+C\Big(\int_{E_{\epsilon^\gamma}^-}|\nabla s_\epsilon|^2\,dx+\big(\int_{E_{\epsilon^\gamma}^-}|\nabla s_\epsilon|^2\,dx\big)^\frac12\Big)=V+VI.
\end{eqnarray}
Again from the estimate \eqref{s-UB1}, we see that
$$
VI\le C\big(\epsilon^{\frac{\gamma}2}+o(1)\big).
$$
$VI$ can be estimated similarly to that of $I$. In fact,
\begin{eqnarray*}
&&\int_{\Omega^-\cap\Gamma_{\epsilon^\gamma}}|\nabla s_\epsilon\cdot n_\epsilon|^2\\
&&\le C\epsilon^{-2}\Big(\epsilon^{2+\gamma}\big\|\nabla^2d_\Gamma\big\|_{L^\infty(\Gamma_{\epsilon^\gamma})}
+\int_{\Omega^-\cap\Gamma_{\epsilon^\gamma}} |\xi'|^2(\frac{d_\Gamma(x)}{\epsilon})
|n(x)-n(\Pi_\Gamma(x))|^2\,dx\Big)\\
&&\le C\epsilon^\gamma+o_\epsilon(1)+C\epsilon^{-2}\int_{\Omega^-\cap(\Gamma_{\epsilon^\gamma}\setminus\Gamma_{L\epsilon}}|n(x)-n(\Pi_\Gamma(x))|^2\,dx\\
&&\le C\epsilon^\gamma+o_\epsilon(1).
\end{eqnarray*}
Hence
$$VI\le C\epsilon^{\frac{\gamma}2}+o_\epsilon(1).$$
Substituting the estimates of $III$ and $IV$ into \eqref{n-UB1.7}, we obtain
\begin{equation}\label{n-UB1.8}
\int_{\Omega^-}\big(L_1|\nabla s_\epsilon\cdot n_\epsilon|^2
+L_3(\nabla s_\epsilon\cdot n_\epsilon)(s_\epsilon{\rm{div}}n_\epsilon)\big)\,dx
\le C\big(\epsilon^{\frac{\gamma}2}+o_\epsilon(1)\big).
\end{equation}
Combining \eqref{s-UB1}, \eqref{n-UB1.3}, \eqref{n-UB1.5}, \eqref{n-UB1.6},
with \eqref{n-UB1.8} and \eqref{sharp_path}, we obtain the upper bound \eqref{A}. \qed

\medskip
\begin{lemma}\label{Upper_Est2} Assume $\Gamma$ and the boundary value
$(t_\epsilon, g_\epsilon)$ satisfies the same assumptions as in Theorem \ref{sharp}.
If $L_2>L_1=L_3=0$, then
\begin{eqnarray}\label{B}
&&\inf\Big\{\int_\Omega \mathcal{W}_\epsilon(s,n,\nabla s, \nabla n)\,dx: \ (s_\epsilon, n_\epsilon)=(t_\epsilon, g_\epsilon) \ {\rm{on}}\ \partial\Omega\Big\}\nonumber\\
&&\le \frac{1}{\epsilon}\int_{0}^{s_+}
\sqrt{2\beta W(\tau)} \mathcal{H}^2(\Gamma(\epsilon\xi_{\epsilon,\gamma}^{-1}(\tau)) \,d\tau+\mathcal{D}_B+o_\epsilon(1),
\end{eqnarray}
where $\mathcal{D}_B$ is given by \eqref{OFM2} and \eqref{bdry2} of Theorem \ref{sharp}.
\end{lemma}
\pf The proof of \eqref{B} can be done almost exactly as in Lemma \ref{Upper_Est1}.
In fact, the construction of $s_\epsilon$ is exactly same as in Lemma \ref{Upper_Est1}.
While the construction $n_\epsilon$ also follows the same procedure,
except that we replace the map $n$, that is a minimizer of $\mathcal{D}_A$ in Lemma \ref{Upper_Est1},
by a map $n$ that minimizes $\mathcal{D}_B$. Namely,
$n\in H^1(\Omega^+,\mathbb S^2)$ satisfies $n=g$ on $\Sigma^-$, $n\wedge \nu_\Gamma=0$ on $\Gamma$,
and
$$
E(n; \Omega^+)=\mathcal{D}_B.
$$

Since every other term in the integral
$\displaystyle\int_\Omega \epsilon \widetilde{\mathcal{W}}_\epsilon(s_\epsilon, n_\epsilon, \nabla s_\epsilon,\nabla n_\epsilon)\,dx$
can be estimated in the same way as in Lemma \ref{Upper_Est1}, it suffices to sketch the estimate of the term
$$\int_{\Omega^+\cap\Gamma_{\epsilon^\gamma}} (L_2|\nabla s_\epsilon\wedge n_\epsilon|^2
+L_4s_\epsilon\nabla s_\epsilon\cdot (\nabla n_\epsilon) n_\epsilon\big)\,dx.
$$
Recall from the condition $n\wedge \nu_\Gamma=0$ on $\Gamma$ that
$$\nabla d_{\Gamma}(\Pi_\Gamma(x))\wedge n(\Pi_\Gamma(x))=0,
\ \forall x\in\Omega^+\cap\Gamma_{\epsilon^\gamma}.$$
From the construction of $(s_\epsilon, n_\epsilon)$, we know that
\begin{eqnarray*}
&&\int_{\Omega^+\cap\Gamma_{\epsilon^\gamma}} L_2|\nabla s_\epsilon\wedge n_\epsilon|^2\,dx
\le C\epsilon^{-2}\int_{\Omega^+\cap\Gamma_{\epsilon^\gamma}}
|\xi_{\epsilon,\gamma}'|^2(\frac{d_\Gamma(x)}{\epsilon}) |\nabla d_{\Gamma}(x)\wedge n(x)|^2\,dx\\
&&\le C\epsilon^{-2}\int_{\Omega^+\cap\Gamma_{\epsilon^\gamma}}
|\xi_{\epsilon,\gamma}'|^2(\frac{d_\Gamma(x)}{\epsilon})
|\nabla d_{\Gamma}(x)\wedge n(x)-\nabla d_{\Gamma}(\Pi_\Gamma(x))\wedge n(\Pi_\Gamma(x))|^2\,dx\\
&&
\le C\epsilon^{-2} \int_{\Omega^+\cap\Gamma_{\epsilon^\gamma}}
|\xi_{\epsilon,\gamma}'|^2(\frac{d_\Gamma(x)}{\epsilon})
(|\nabla d_{\Gamma}(x)-\nabla d_{\Gamma}(\Pi_\Gamma(x))|^2+|n(x)-n(\Pi_\Gamma(x))|^2)\,dx\\
&&\le C\epsilon^\gamma+o_\epsilon(1).
\end{eqnarray*}
This implies that
\begin{equation*}
|\int_{\Omega^+\cap\Gamma_{\epsilon^\gamma}}
L_4s_\epsilon\nabla s_\epsilon\cdot (\nabla n_\epsilon) n_\epsilon\,dx|
\le C\big(\int_{\Omega^+\cap\Gamma_{\epsilon^\gamma}} |\nabla s_\epsilon\wedge n_\epsilon|^2\,dx\big)^\frac12
\big(\int_{\Omega^+\cap\Gamma_{\epsilon^\gamma}} |\nabla n|^2\big)^\frac12=o_\epsilon(1).
\end{equation*}
Thus the estimate \eqref{B} holds.\qed

\medskip
\begin{lemma}\label{Upper_Est3} Assume $\Gamma$ and the boundary values
$(s_\epsilon, n_\epsilon)$ satisfy the same assumptions as in Theorem \ref{sharp}.
If $L_1=L_2=L_3=L_4=0$, then
\begin{eqnarray}\label{C}
&&\inf\Big\{\int_\Omega \mathcal{W}_\epsilon(s,n,\nabla s, \nabla n)\,dx:
\ (s_\epsilon, n_\epsilon)=(t_\epsilon, g_\epsilon) \ {\rm{on}}\ \partial\Omega\Big\}\nonumber\\
&&\le \frac{1}{\epsilon}\int_{0}^{s_+}
\sqrt{2\beta W(\tau)} \mathcal{H}^2(\Gamma(\epsilon\xi_{\epsilon,\gamma}^{-1}(\tau)) \,d\tau+\mathcal{D}_C+o_\epsilon(1),
\end{eqnarray}
where $\mathcal{D}_C$ is given by \eqref{OFM3} and \eqref{bdry3} of Theorem \ref{sharp}.
\end{lemma}
\pf The proof of \eqref{C} can be done almost exactly as in Lemma \ref{Upper_Est1}. In fact,
the construction and estimate of $s_\epsilon$ is exactly same as in that in Lemma \ref{Upper_Est1}.
While in the construction of $n_\epsilon$, we simply replace the minimizer $n$ of
$\mathcal{D}_A$ in Lemma \ref{Upper_Est1}
by a map $n$ that minimizes $\mathcal{D}_C$. Namely, $n\in H^1(\Omega^+, \mathbb S^2)$
satisfies $n=g$ on $\Sigma^+$, $\frac{\partial n}{\partial \nu_\Gamma}=0$ on $\Gamma$, and
$$
E(n;\Omega^+)=\mathcal{D}_C.
$$
Since $L_1=L_2=L_3=L_4=0$,
$$\int_\Omega \widetilde{\mathcal{W}}_\epsilon(s_\epsilon, n_\epsilon, \nabla s_\epsilon, \nabla n_\epsilon)\,dx
=\int_\Omega \big(s_\epsilon^2 W_{OF}(n_\epsilon, \nabla n_\epsilon)+\alpha s_\epsilon^2|\nabla n_\epsilon|^2
+\beta|\nabla s_\epsilon|^2+ \frac{1}{\epsilon^2} W(s_\epsilon)\big)\,dx$$
%\end{document}
can be estimated as in \eqref{s-UB1}, \eqref{n-UB1.3}, and \eqref{n-UB1.6} of Lemma \ref{Upper_Est1}.
\qed

\subsection{Refined energy lower bounds for the case (A) \fbox{$L_1>0=L_2=L_4$,\ $3L_3^2<4L_1\alpha$}}

In this subsection, we will establish an improved lower bound of energy that matches the refined
upper bound of energy, which ensures the planar anchoring condition for the limiting director field
on the sharp interface $\Gamma$.

First, it follows from $3L_3^2<4L_1\alpha$ that there exists a positive number $\mu>0$ such that
$3L_3^2\le 4(L_1-\mu)(\alpha-\mu)$ so that by Cauchy-Schwarz inequality we have
\begin{equation}\label{CS3}
|L_3(\nabla s_\epsilon\cdot n_\epsilon)s_\epsilon {\rm{div}} n_\epsilon|
\le \sqrt{3}|L_3| |\nabla s_\epsilon\cdot n_\epsilon| |s_\epsilon| |\nabla n_\epsilon|
\le (L_1 -\mu) (\nabla s_\epsilon\cdot n_\epsilon)^2+(\alpha-\mu) s_\epsilon^2 |\nabla n_\epsilon|^2.
\end{equation}
This implies that
\begin{eqnarray*}
\widetilde{\mathcal{W}}_\epsilon(s_\epsilon, n_\epsilon,\nabla s_\epsilon, \nabla n_\epsilon)
&=&( \beta |\nabla s_\epsilon|^2+\frac{1}{\epsilon^2}W(s_\epsilon))+s_\epsilon^2 W_{OF}(n_\epsilon, \nabla n_\epsilon)\\
&&+\big(\alpha s_\epsilon^2|\nabla n_\epsilon|^2
+L_1(\nabla s_\epsilon\cdot n_\epsilon)^2
+L_3(\nabla s_\epsilon\cdot n_\epsilon)s_\epsilon {\rm{div}} n_\epsilon\big)\\
&\ge& \beta |\nabla s_\epsilon|^2+\frac{1}{\epsilon^2}{W(s_\epsilon)}
+\mu s_\epsilon^2 |\nabla n_\epsilon|^2+
\mu|\nabla s_\epsilon\cdot n_\epsilon|^2\\
&=&\big[(\beta+\gamma \cos^2\theta_\epsilon) |\nabla s_\epsilon|^2+\frac{1}{\epsilon^2}W(s_\epsilon)\big]
+\mu s_\epsilon^2 |\nabla n_\epsilon|^2\\
&\ge& \frac{2}{\epsilon}|\nabla s_\epsilon|\sqrt{W(s_\epsilon)}\sqrt{\beta+\mu\cos^2\theta_\epsilon}
+\mu s_\epsilon^2 |\nabla n_\epsilon|^2,
\end{eqnarray*}
where $\cos\theta_\epsilon=\frac{\nabla s_\epsilon}{|\nabla s_\epsilon|}\cdot n_\epsilon$.

Thus by the co-area formula we obtain
\begin{eqnarray}\label{LB3}
&&\int_\Omega \widetilde{\mathcal{W}}_\epsilon(s_\epsilon, n_\epsilon,\nabla s_\epsilon, \nabla n_\epsilon)\,dx\nonumber\\
&&\ge \int_\Omega \big(\frac{2}{\epsilon}|\nabla s_\epsilon|\sqrt{W(s_\epsilon)}\sqrt{\beta+\mu\cos^2\theta_\epsilon}+\mu s_\epsilon^2 |\nabla n_\epsilon|^2\big)\,dx\nonumber\\
&& \ge \mu  \int_\Omega s_\epsilon^2 |\nabla n_\epsilon|^2\,dx
+\frac{2}{\epsilon}\int_0^{s_+}\sqrt{W(\tau)} \int_{\partial^*\mathcal{S}_\epsilon(\tau)} \sqrt{\beta+\mu\cos^2\theta_\epsilon}\,d\mathcal{H}^2 \,d\tau\nonumber\\
&&=\mu  \int_\Omega s_\epsilon^2 |\nabla n_\epsilon|^2\,dx
+\frac{2}{\epsilon}\int_0^{s_+}\sqrt{W(\tau)} \int_{\partial^*\mathcal{S}_\epsilon(\tau)} \big(\sqrt{\beta+\mu\cos^2\theta_\epsilon}-\sqrt{\beta}\big)\,d\mathcal{H}^2 \,d\tau\nonumber\\
&&\ \ \ +\frac{2}{\epsilon}\int_0^{s_+}\sqrt{\beta W(\tau)} \mathcal{H}^2(\partial^*\mathcal{S}_\epsilon(\tau))\,d\tau,
\end{eqnarray}
where $\mathcal{S}_\epsilon(\tau)=\big\{x\in\Omega: \ s_\epsilon(x)\ge\tau\big\}$.

It follows from the assumption of $s_\epsilon$ on $\partial\Omega$ that  for any $0<\tau<s_+$,
the enclosed surface $\mathcal{T}_\epsilon(\tau)\subset\partial\Omega$
between $\partial^*\mathcal{S}_\epsilon(\tau)$ and $\Sigma$ is
a strip with width at most $C\epsilon$. Hence by the area minimality of $\Gamma$,
we have
\begin{eqnarray}\label{area_comp1}
\mathcal{H}^2(\Gamma)&\le& \mathcal{H}^2\big(\partial^*\mathcal{S}_\epsilon(\tau)\cup \mathcal{T}_\epsilon(\tau)\big)
= \mathcal{H}^2\big(\partial^*\mathcal{S}_\epsilon(\tau)\big)+\mathcal{H}^2\big(\mathcal{T}_\epsilon(\tau)\big)\nonumber\\
&\le& \mathcal{H}^2(\partial^*\mathcal{S}_\epsilon(\tau))+C\epsilon, \ 0<\tau<s_+.
\end{eqnarray}
This implies that
\begin{equation}\label{LB4}
\frac{2}{\epsilon}\int_0^{s_+}\sqrt{\beta W(\tau)} \mathcal{H}^2(\partial^*\mathcal{S}_\epsilon(\tau))\,d\tau
\ge\frac{1}{\epsilon} \alpha_0\mathcal{H}^2(\Gamma)-C.
\end{equation}
Notice that
\begin{eqnarray*}
\sqrt{\beta+\mu\cos^2\theta_\epsilon}-\sqrt{\beta}&=&\frac{\mu\cos^2\theta_\epsilon}{\sqrt{\beta+\mu\cos^2\theta_\epsilon}+\sqrt{\beta}}\\
&\ge&\mu_*\cos^2\theta_\epsilon,
\end{eqnarray*}
where $\displaystyle\mu_*=\frac{\mu}{\sqrt{\beta+\mu}+\sqrt{\beta}}>0$.

Hence, by matching the refined upper bound \eqref{A} with \eqref{LB3}, we conclude that
\begin{equation}\label{UB5}
\mu \int_\Omega (s_\epsilon^2 |\nabla n_\epsilon|^2+|\nabla s_\epsilon\cdot n_\epsilon|^2)\,dx \le \mathcal{D}_A+C+o_\epsilon(1),
\end{equation}
\begin{equation}\label{UB6}
\frac{\mu_*}{\epsilon}\int_0^{s_+}\sqrt{W(\tau)}
 \int_{\partial^*\mathcal{S}_\epsilon(\tau)} \cos^2\theta_\epsilon\,d\mathcal{H}^2 \,d\tau
 \le \mathcal{D}_A+C+o_\epsilon(1),
\end{equation}
and
\begin{equation}\label{UB7}
\frac{2}{\epsilon}\int_0^{s_+}\sqrt{\beta W(\tau)} \mathcal{H}^2(\partial^*\mathcal{S}_\epsilon(\tau))\,d\tau
\le \frac{1}{\epsilon}\alpha_0\mathcal{H}^2(\Gamma)+\mathcal{D}_A+o_\epsilon(1).
\end{equation}

%Hence by using the refined upper bound \eqref{A} and the co-area formula, we obtain
%\begin{eqnarray}\label{matchA}
%&&\frac{\alpha}2 \int_\Omega s_\epsilon^2 |\nabla n_\epsilon|^2\,dx
%+\frac{\gamma}{\epsilon(\sqrt{\beta+\gamma}+\sqrt{\beta})}\int_0^{s_+}
%\sqrt{W_0(\delta)}\int_{\Gamma_\delta^\epsilon}\sqrt{\beta+\frac{L_1}{2}\cos^2\theta_\epsilon}\,dH^2d\delta\nonumber\\
%&&\le\frac{\alpha}2 \int_\Omega s_\epsilon^2 |\nabla n_\epsilon|^2\,dx
%+\frac{2}{\epsilon}\int_\Omega |\nabla s_\epsilon|\sqrt{W_0(s_\epsilon)}
%\sqrt{\beta+\frac{L_1}{2}\cos^2\theta_\epsilon}\,dx\nonumber\\
%&&\le \int_\Omega \mathcal{W}^\epsilon(s_\epsilon, n_\epsilon,\nabla s_\epsilon, \nabla n_\epsilon)\,dx\nonumber\\
%&&\le \frac{\alpha_0}{\epsilon} H^2(\Gamma)+D_A+o_\epsilon(1).
%\end{eqnarray}
%Here $\Gamma_\delta^\epsilon=\big\{x\in\Omega: \ s_\epsilon(x)=\delta\big\}$.
%With some additional estimates, we can have
%\begin{equation}\label{matchA1}
%\int_0^{s_+}2\sqrt{\beta W_0(\delta)}H^2(\Gamma_\delta^\epsilon)\,d\delta=\alpha_0 H^2(\Gamma)+\epsilon o_\epsilon(1).
%\end{equation}
%It follows from \eqref{matchA} and \eqref{matchA1} that
%\begin{eqnarray}\label{matchA2}
%\frac{1}{\epsilon}\int_0^{s_+}
%2\sqrt{\beta W_0(\delta)}\int_{\Gamma_\delta^\epsilon}\big[\sqrt{1+\frac{L_1}{2\beta}\cos^2\theta_\epsilon}-1\big]\,dH^2d\delta
%+\frac{\alpha}2 \int_\Omega s_\epsilon^2 |\nabla n_\epsilon|^2\,dx\le D_A+o_\epsilon(1).
%\end{eqnarray}
For any fixed $\delta>0$ and for any $\epsilon\in (0,1)$, applying Fubini's theorem to \eqref{UB6} and \eqref{UB7}
we obtain that there exist $C>0$, that is independent of $\delta$, and $\delta_\epsilon\in (\delta, 2\delta)$ such that
 \begin{equation}\label{slice6}
\displaystyle \int_{\partial^*\mathcal{S}_\epsilon(s_+-{\delta_\epsilon})} \cos^2 \theta_\epsilon \,d\mathcal{H}^2
\le C(\beta, \mathcal{D}_A)\frac{\epsilon}{\delta},
\end{equation}
\begin{equation}\label{slice6.0}
\displaystyle \mathcal{H}^2(\partial^*\mathcal{S}_\epsilon(s_+ -{\delta_\epsilon}))\le
\mathcal{H}^2(\Gamma)+C\frac{\epsilon}{\delta}.
 \end{equation}
Since it is  straightforward to get \eqref{slice6}, we only sketch how to obtain \eqref{slice6.0}.
From \eqref{area_comp1}, we have
\begin{eqnarray*}
&&\frac{2}{\epsilon} \Big\{\int_{0}^{s_+-2\delta}+\int_{s_+ -\delta}^{s_+} \Big\}
\sqrt{\beta W(\tau)}\mathcal{H}^2\big(\partial^*{\mathcal{S}_\epsilon(\tau)}\big)\,d\tau\\
&&\ge \frac{2}{\epsilon} \Big(\int_{[0, {s_+-2\delta}]\cup [s_- -\delta, s_+]} \sqrt{\beta W(\tau)}\,d\tau\Big)
\mathcal{H}^2(\Gamma) -C,
\end{eqnarray*}
so that
\begin{eqnarray*}
&&\Big(\inf_{s_+-2\delta\le\tau\le s_+-\tau} \mathcal{H}^2(\partial^*\mathcal{S}_\epsilon(\tau))\Big)
 \frac{2}{\epsilon} \int_{s_+-2\delta}^{s_+ -\delta} \sqrt{\beta W(\tau)}\,d\tau\\
&&\le \frac{2}{\epsilon} \int_{s_+-2\delta}^{s_+ -\delta}
\sqrt{\beta W(\tau)}\mathcal{H}^2\big(\partial^*{\mathcal{S}_\epsilon(\tau)}\big)\,d\tau\\
&&\le \Big(\frac{2}{\epsilon} \int_{s_+-2\delta}^{s_+ -\delta} \sqrt{\beta W(\tau)}\,d\tau\Big)
\mathcal{H}^2(\Gamma) +\mathcal{D}_A+C+o_\epsilon(1).
\end{eqnarray*}
This, combined with the estimate
$\displaystyle\int_{s_+-2\delta}^{s_+ -\delta} \sqrt{\beta W(\tau)}\,d\tau\approx \delta$, yields
\eqref{slice6.0}.

If we set
$$\Omega_{\epsilon,\delta}^-=\Big\{x\in\Omega: \ s_\epsilon(x)<\delta\Big\},
\ \ \Omega_{\epsilon,\delta}^+=\Big\{x\in\Omega: \ s_\epsilon(x)> s_+-\delta\Big\},$$
and
$$E_{\epsilon, \delta}=\Big\{x\in\Omega: \ \delta\le s_\epsilon(x)\le s_+-\delta\Big\},$$
then it follows from \eqref{W-cond} and \eqref{UB5} that
\begin{equation}\label{slice7}
 \begin{cases} \displaystyle |E_{\epsilon,{\delta_\epsilon}}|\le \frac{C\epsilon}{\delta^2}, \\
\displaystyle \int_{\Omega_{\epsilon,\delta_\epsilon}^+} |\nabla n_\epsilon|^2
\le \frac{C(\alpha, \mathcal{D}_A)}{(s_+-\delta)^2}.
 \end{cases}
 \end{equation}
From \eqref{W-cond} and
$$\int_{\Omega_{\epsilon, \delta_\epsilon}^+} \frac{1}{\epsilon^2}W(s_\epsilon)\,dx\leq \frac{C}{\epsilon},$$
we see that for a.e. $x\in \Omega_{\epsilon, \delta_\epsilon}^+$, $s_+-\delta_\epsilon\le s_\epsilon(x)\le 1$ and
$s_\epsilon(x)\rightarrow s_+$. Hence for any $1<p<\infty$, it holds that
\begin{equation}\label{UB8}
\int_{\Omega_{\epsilon,\delta_\epsilon}^+}|s_\epsilon(x)-s_+|^p\,dx\rightarrow 0.
\end{equation}

Furthermore, from the Cauchy-Schwarz inequality and \eqref{LB4}, we also have
\begin{eqnarray*}
\int_{\Omega}( \beta |\nabla s_\epsilon|^2+\frac{1}{\epsilon^2}W(s_\epsilon))\,dx
&\ge& \frac{2}{\epsilon}\int_0^{s_+} \sqrt{\beta W(\tau)} \mathcal{H}^2(\partial^*\mathcal{S}_\epsilon(\tau))\,d\tau\\
&\ge&\frac{1}{\epsilon} \alpha_0\mathcal{H}^2(\Gamma)-C.
\end{eqnarray*}
Matching with the refined upper bound \eqref{A},  this also implies that
\begin{eqnarray}\label{UB9}
&&\int_{\Omega}\big(s_\epsilon^2 W_{OF}(n_\epsilon, \nabla n_\epsilon)+\alpha s_\epsilon^2|\nabla n_\epsilon|^2
+L_1(\nabla s_\epsilon\cdot n_\epsilon)^2
+L_3(\nabla s_\epsilon\cdot n_\epsilon)s_\epsilon {\rm{div}} n_\epsilon\big)\,dx\nonumber\\
&&\le \mathcal{D}_A+C+o_\epsilon(1).
\end{eqnarray}

From \eqref{slice6} and \eqref{slice7}, there exist $\epsilon_i\to 0$, $\delta_i=\delta_{\epsilon_i}\in (\epsilon_i^\frac14, 2\epsilon_i^\frac14)$,
a set $\Omega_*\subset\Omega$ of finite perimeter
such that for $\Omega_i=\Omega_{\epsilon_i,\delta_i}^+$, it holds
that
\begin{itemize}
\item [(a)] $\chi_{\Omega_i}\rightharpoonup \chi_{\Omega_*}$ in $BV(\Omega)$,
and $\chi_{\Omega_i}\rightarrow\chi_{\Omega_*}$ in $L^1(\R^3)$.
\item [(b)] We also view $\Gamma_i=\partial\Omega_i\lfloor\Omega$ and $\Gamma_*=\partial^*\Omega_*\lfloor\Omega$
as oriented boundaries and integral rectifiable 2-currents, and use the same notations,
i.e., $\Gamma_i=[[\partial\Omega_i \lfloor \Omega]]$
and $\Gamma_*=[[\partial\Omega_*\lfloor\Omega]]$. Then
$$\Gamma_i\rightharpoonup \Gamma_*$$
weakly converges as oriented boundaries and integral rectifiable 2-currents. By the lower
semicontinuity, we have that
\begin{equation}\label{LSC1}
\mathcal{H}^2(\Gamma)\le \mathcal{H}^2(\Gamma_*)\le \lim_i \mathcal{H}^2(\Gamma_i)=\mathcal{H}^2(\Gamma),
\end{equation}
where the first inequality follows from the area minimality of $\Gamma$, since
$\partial \Gamma=\partial\Gamma_*=[[\Sigma^0]]$.
\end{itemize}

Since $\Gamma$ is assumed to be a unique area minimizing surface spanned by $\Sigma^0$, we have that
$\Gamma_*=\Gamma$.  Also, since
$$\partial \Omega_i\lfloor (\R^3\setminus\Omega)\rightarrow \Sigma^+$$
as convergence of currents, we conclude that
$$\partial \Omega_*\lfloor (\R^3\setminus\Omega)=\Sigma^+.$$
Therefore $\Omega_*=\Omega^+$.  Next we need to show

\medskip
\noindent{\it Claim 1}. {\it For any $\eta>0$, it holds
\begin{equation}\label{no_hole}
\lim_{\epsilon_i\rightarrow 0} \mathcal{H}^2(\partial\Omega_i\cap\big\{x\in\Omega: d_\Gamma(x)\ge \eta\big\}=0.
\end{equation}}
It follows from \eqref{LSC1} that
$$\mathcal{H}^2\lfloor (\partial\Omega_i\cap \Omega)\rightharpoonup
\mathcal{H}^2\lfloor \Gamma$$
as weak convergence of Radon measures. Hence by the lower semicontinuity,
\begin{eqnarray*}
\mathcal{H}^2(\Gamma)&=&\mathcal{H}^2\big(\Gamma\cap\big\{x\in\Omega: d_\Gamma(x)< \eta\big\}\big)\\
&\le &\liminf_{\epsilon_i\to 0} \mathcal{H}^2\big(\partial\Omega_i\cap \big\{x\in\Omega: d_\Gamma(x)< \eta\big\}\big)\\
&\le& \liminf_{\epsilon_i\to 0} \mathcal{H}^2(\partial\Omega_i\cap\Omega)=\mathcal{H}^2(\Gamma).
\end{eqnarray*}
This clearly implies \eqref{no_hole}.

\medskip
\noindent{\it Claim 2}. {\it  There exists a map $n\in {\rm{SBV}}(\Omega^+,\mathbb S^2)$ \footnote{Here SBV$(\Omega)$ denotes the space of all BV (or bounded variations) functions such that the Cantor part of the distributional derivatives is zero. See Ambrosio \cite{Ambrosio} for more discussions.}
such that after passing to a subsequence,
$$n_{\epsilon_i}\chi_{\Omega_i} \rightharpoonup n\chi_{\Omega^+} \ \ {\rm{in}}\ \ {\rm{BV}}(\Omega),
\ \ {\rm{and}}\ \  s_{\epsilon_i}\chi_{\Omega_i} \rightarrow s_+\chi_{\Omega^+} \ \ {\rm{in}}\ \ L^2(\Omega).
$$
Furthermore, $n\in H^1(\Omega^+, \mathbb S^2)$.}

To show Claim 2, we first observe that the absolutely continuous part of the distributional derivative
of $v_i=n_{\epsilon_i}\chi_{\Omega_i}$ is $\nabla v_i=\nabla n_{\epsilon_i}\chi_{\Omega_i}$, which is
uniformly bounded in $L^2$, i.e.
$$\int_\Omega |\nabla v_i|^2\,dx=\int_{\Omega_i} |\nabla n_i|^2\,dx\le C.$$
The jump part $J_{v_i}$ of $v_i$ satisfies
$$J_{v_i}\subset\partial\Omega_i\cap\Omega=\Gamma_i,$$
so that
$$
\mathcal{H}^2(J_{v_i})\le \mathcal{H}^2(\Gamma_i)\le 2\mathcal{H}^2(\Gamma).
$$
Moreover, we have that
$$\|v_i\|_{L^\infty(\Omega)}\le 1.$$
Thus it follows from \cite{Ambrosio}
that $\{v_i\}\subset {\rm{SBV}}(\Omega)$ is a weakly compact sequence in ${\rm{SBV}}(\Omega)$.
There exists a $n\in {\rm{SBV}}(\Omega)$ such that $v_i\rightharpoonup n$ in ${\rm{BV}}(\Omega)$
and strongly in $L^1(\Omega)$. Since $|v_i|=1$ in $\Omega_i$
and $|v_i|=0$ in $\Omega\setminus\Omega_i$, it follows that $|n|=1$ in $\Omega^+$ and $|n|=0$
in $\Omega\setminus\Omega^+$ so that $n\in {\rm{SBV}}(\Omega^+,\mathbb S^2)$. From the lower
semicontinuity, we have that
\begin{equation}\label{LSC5}
\int_{\Omega^+}|\nabla n|^2\,dx\le\liminf_{\epsilon_i\to 0} \int_{\Omega_i} |\nabla n_{\epsilon_i}|^2\,dx
\end{equation}
Now we want to show its jump set has $\mathcal{H}^2$-measure zero. This follows
from \eqref{no_hole} and the lower semicontinuity:
\begin{eqnarray*}
\mathcal{H}^2\big(J_n\cap\big\{x\in\Omega^+: \ d_\Gamma(x)>\eta\big\}\big)
&\le&\liminf_{\epsilon_i\to 0} \mathcal{H}^2\big(J_{n_{\epsilon_i}}\cap\big\{x\in\Omega^+: \ d_\Gamma(x)>\eta\big\}\big)\\
&\le&\liminf_{\epsilon_i\to 0} \mathcal{H}^2\big(\partial\Omega_i\cap\big\{x\in\Omega^+: \ d_\Gamma(x)>\eta\big\}\big)=0.
\end{eqnarray*}
This, after sending $\eta\to 0$, yields $\mathcal{H}^2\big(J_n\cap\Omega^+)=0$. Hence $n\in H^1(\Omega^+)$.

\bigskip
It follows from (\ref{LSC5}) that $\nabla n_{\epsilon_i}\chi_{\Omega_i}\rightharpoonup \nabla n\chi_{\Omega^+}$ in
$L^2(\Omega)$, and hence
$$s_{\epsilon_i} \nabla n_{\epsilon_i}\chi_{\Omega_i}\rightharpoonup s_+ \nabla n\chi_{\Omega^+}
\ \ {\rm{in}}\ \ L^1(\Omega).$$
This, combined with the uniform $H^1$ bound \eqref{slice7}, further implies
\begin{equation}\label{alpha-conv}
s_{\epsilon_i} \nabla n_{\epsilon_i}\chi_{\Omega_i}\rightharpoonup s_+ \nabla n\chi_{\Omega^+}
\ \ {\rm{in}}\ \ L^2(\Omega).
\end{equation}
We claim that
\begin{equation}\label{distri_limit}
\nabla s_{\epsilon_i}\cdot n_{\epsilon_i}\rightarrow 0 \ \ {\rm{in}}\ \ \mathcal{D}'(\Omega^+).
\end{equation}
To see this, let $\phi\in C^\infty_0(\Omega^+)$. Then by integration by parts we have
\begin{eqnarray*}
\int_{\Omega^+} \nabla s_{\epsilon_i}\cdot n_{\epsilon_i} \phi\,dx
&=&-\int_{\Omega^+} s_{\epsilon_i} \big({\rm{div}}n_{\epsilon_i} \phi+n_{\epsilon_i}\nabla\phi\big)\,dx\\
&\rightarrow& -\int_{\Omega^+} s_{+} \big({\rm{div}}n \phi+n\nabla\phi\big)\,dx\\
&=& -\int_{\Omega^+} s_{+} {\rm{div}}(n \phi)\,dx=0, \ {\rm{as}}\ \epsilon_i\rightarrow 0.
\end{eqnarray*}
Since $\displaystyle\int_{\Omega_i} |\nabla s_{\epsilon_i}\cdot n_{\epsilon_i}|^2\,dx$ is uniformly bounded,
it follows from \eqref{distri_limit} that
\begin{equation}\label{L_1-conv}
\nabla s_{\epsilon_i}\cdot n_{\epsilon_i}\chi_{\Omega_i}\rightharpoonup 0 \ \ {\rm{in}}\ \ L^2(\Omega).
\end{equation}
It follows from \eqref{alpha-conv} and the lower semicontinuity of $\displaystyle\int_{\Omega^+}
s_{\epsilon_i}^2W_{OF}(n_{\epsilon_i}, \nabla n_{\epsilon_i})\,dx$ that
\begin{eqnarray}\label{lsc3}
\int_{\Omega^+} s_{+}^2W_{OF}(n, \nabla n)\,dx&=&\int_{\Omega} s_{+}^2W_{OF}(n, \nabla n)\chi_{\Omega^+}\,dx\nonumber\\
&\le&\liminf_{\epsilon_i\rightarrow 0}
\int_{\Omega} s_{\epsilon_i}^2W_{OF}(n_{\epsilon_i}, \nabla n_{\epsilon_i})\chi_{\Omega_i}\,dx\nonumber\\
&=&\liminf_{\epsilon_i\rightarrow 0}
\int_{\Omega_i} s_{\epsilon_i}^2W_{OF}(n_{\epsilon_i}, \nabla n_{\epsilon_i})\,dx.
\end{eqnarray}
Next we claim that
\begin{equation}\label{lsc4}
\int_{\Omega^+} \alpha s_+^2|\nabla n|^2\,dx
\le\liminf_{\epsilon_i\to 0}
\int_{\Omega_i}\big(\alpha s_{\epsilon_i}^2|\nabla n_{\epsilon_i}|^2+L_1(\nabla s_{\epsilon_i}\cdot n_{\epsilon_i})^2
+L_3(\nabla s_{\epsilon_i}\cdot n_{\epsilon_i})s_{\epsilon_i} {\rm{div}} n_{\epsilon_i}\big)\,dx.
\end{equation}
For $\eta>0$, define $\Omega^+_\eta=\big\{x\in\Omega^+: \ d(x,\partial\Omega^+)>\eta\big\}$.
Since $\Omega_i\rightarrow \Omega^+$ in Hausdorff distance, we may assume that
for $i$ sufficiently large, $\Omega^+_{\eta}\subset\Omega_i$ and hence
\begin{equation}\label{alpha-L1}
s_{\epsilon_i}\nabla n_{\epsilon_i}\rightharpoonup s_+\nabla n\ {\rm{in}}\ L^2(\Omega^+_\eta),
\ \ \nabla s_{\epsilon_i}\cdot n_{\epsilon_i}\rightharpoonup 0 \ {\rm{in}}\ L^2(\Omega^+_\eta).
\end{equation}
This and \eqref{CS3} imply that
\begin{eqnarray*}
&&\mathcal{D}\equiv\liminf_{\epsilon_i\to 0}
\int_{\Omega_i}\big(\alpha s_{\epsilon_i}^2|\nabla n_{\epsilon_i}|^2+L_1(\nabla s_{\epsilon_i}\cdot n_{\epsilon_i})^2
+L_3(\nabla s_{\epsilon_i}\cdot n_{\epsilon_i})s_{\epsilon_i} {\rm{div}} n_{\epsilon_i}\big)\,dx\\
&&\ge\liminf_{\epsilon_i\to 0}
\int_{\Omega^+_\eta}\big(\alpha s_{\epsilon_i}^2|\nabla n_{\epsilon_i}|^2+L_1(\nabla s_{\epsilon_i}\cdot n_{\epsilon_i})^2
+L_3(\nabla s_{\epsilon_i}\cdot n_{\epsilon_i})s_{\epsilon_i} {\rm{div}} n_{\epsilon_i}\big)\,dx\\
&&
=\liminf_{\epsilon_i\to 0}
\int_{\Omega^+_\eta}\Big[\alpha |s_{\epsilon_i}(\nabla n_{\epsilon_i}-\nabla n)+s_{\epsilon_i}\nabla n|^2
+L_1(\nabla s_{\epsilon_i}\cdot n_{\epsilon_i})^2\\
&&\qquad\qquad\qquad+L_3(\nabla s_{\epsilon_i}\cdot n_{\epsilon_i})(s_{\epsilon_i} ({\rm{div}} n_{\epsilon_i}-{\rm{div}}n)+s_{\epsilon_i}{\rm{div}}n)\Big]\,dx\\
&&=\liminf_{\epsilon_i\to 0}
\Big\{\int_{\Omega^+_\eta}\alpha s_{\epsilon_i}^2|\nabla n|^2\,dx\\
&&\qquad\qquad+\int_{\Omega^+_\eta}\Big[\alpha s_{\epsilon_i}^2|\nabla (n_{\epsilon_i}-n)|^2
+L_1(\nabla s_{\epsilon_i}\cdot n_{\epsilon_i})^2
+L_3(\nabla s_{\epsilon_i}\cdot n_{\epsilon_i})
s_{\epsilon_i} {\rm{div}} (n_{\epsilon_i}-n)\Big]\,dx\\
&&\qquad\qquad+\int_{\Omega^+_\eta}\Big(2\alpha s_{\epsilon_i} \nabla (n_{\epsilon_i}-n) (s_{\epsilon_i}\nabla n)
+L_3(\nabla s_{\epsilon_i}\cdot n_{\epsilon_i}) (s_{\epsilon_i}{\rm{div}} n)\Big)\,dx\Big\}\\
&&=\liminf_{\epsilon_i\to 0} (A_i+B_i+C_i).
\end{eqnarray*}
Applying \eqref{CS3} with $n_\epsilon$ replaced by $n_{\epsilon_i}-n$, we have that $B_i\ge 0$.
From \eqref{UB8} and \eqref{alpha-L1}, we see that
$C_i\rightarrow 0$. On the other hand, since $s_{\epsilon_i}\nabla n\rightharpoonup s_+\nabla n$
in $L^2(\Omega)$, we have
$$\liminf_{\epsilon_i\rightarrow 0} A_i\ge \int_{\Omega^+_\eta}\alpha s_+^2|\nabla n|^2\,dx.$$
Therefore we obtain that
$$\mathcal{D}\ge \int_{\Omega^+_\eta}\alpha s_+^2|\nabla n|^2\,dx.$$
Sending $\eta$ to zero, this implies that
\begin{equation}\label{lsc5}
\mathcal{D}\ge \int_{\Omega^+}\alpha s_+^2|\nabla n|^2\,dx.
\end{equation}
Combining \eqref{lsc5}
\begin{eqnarray}\label{lsc6}
&&E(n,\Omega^+)=\int_{\Omega^+} s_+^2\big(W_{OF}(n,\nabla n)+\alpha|\nabla n|^2\big)\,dx\\
&&\le\liminf_{\epsilon_i\rightarrow 0}\int_{\Omega}\big(s_{\epsilon_i}^2 W_{OF}(n_{\epsilon_i}, \nabla n_{\epsilon_i})+\alpha s_{\epsilon_i}^2|\nabla n_{\epsilon_i}|^2
+L_1(\nabla s_{\epsilon_i}\cdot n_{\epsilon_i})^2
+L_3(\nabla s_{\epsilon_i}\cdot n_{\epsilon_i})s_{\epsilon_i} {\rm{div}} n_{\epsilon_i}\big)\,dx.\nonumber
\end{eqnarray}

From the assumption on $g_{\epsilon_i}$ on $\Sigma^+$, we see that
$n=g$ on $\Sigma^+$. Next we want to show the trace of $n$ on $\Gamma$ satisfies the planar anchoring
condition:
\begin{equation}\label{trace}
n\cdot \nu_{\Gamma} =0\ \ {\rm{on}}\ \ \Gamma.
\end{equation}

%We will leave the proof of the boundary condition of $n$ on $\Gamma$ in the subsection below.
%\subsection{Planar anchoring condition on $\Gamma_*$}

\noindent\underline{\it Sketch of proof of \eqref{trace}}. We will  show the planar anchoring condition of
$n$ on $\Gamma$ as follows. For simplicity, write $n_i=n_{\epsilon_i}$.
First it is not hard to show that as $i\to\infty$, $\Omega_i\to \Omega^+\ \mbox{in Hausdorff distance},$
$$d\mathcal{H}^{2}\lfloor \Gamma_i \rightharpoonup d\mathcal{H}^{2}\lfloor \Gamma$$
as convergence of Radon measures,  and
$$n_i\rightharpoonup n\ \ {\rm{in}}\ \ H^1(\Omega^+,\mathbb S^2).$$
This implies that
$${\rm{div}}(n_i)\rightharpoonup {\rm{div}}(n) \ {\rm{in}}\ L^2(\Omega),\ \
\chi_{\Omega_i}\rightarrow \chi_{\Omega^+} \ {\rm{in}}\ L^2(\mathbb R^3).$$
Therefore
\begin{eqnarray}\label{bdry_conv}
\int_{\Gamma_i}n_i\cdot \nu_{\Gamma_i}\,d\mathcal{H}^2
&=&\int_{\partial\Omega_i} n_i\cdot\nu_{\partial\Omega_i}\,d\mathcal{H}^2
-\int_{\partial\Omega_i\cap\partial\Omega} g_i\cdot\nu_{\partial\Omega}\,d\mathcal{H}^2\nonumber\\
&=&\int_{\Omega_i} {\rm{div}}(n_i)\,dx
-\int_{\partial\Omega_i\cap\partial\Omega} g_i\cdot\nu_{\partial\Omega}\,d\mathcal{H}^2\nonumber\\
&=&\int_{\mathbb R^3} {\rm{div}}(n_i)\chi_{\Omega_i}\,dx
-\int_{\partial\Omega_i\cap\partial\Omega} g_i\cdot\nu_{\partial\Omega}\,d\mathcal{H}^2\nonumber\\
&\to &\int_{\mathbb R^3} {\rm{div}}(n)\chi_{\Omega^+}\,dx-\int_{\partial\Omega^+\cap\partial\Omega} g\cdot\nu_{\partial\Omega}\,d\mathcal{H}^2
\nonumber\\
&=& \int_{\partial\Omega^+} n\cdot\nu_{\partial\Omega^+}\,d\mathcal{H}^2-\int_{\partial\Omega^+\cap\partial\Omega} g\cdot\nu_{\partial\Omega}\,d\mathcal{H}^2
\nonumber\\
&=& \int_{\Gamma} n\cdot\nu_{\Gamma}\,d\mathcal{H}^2.
\end{eqnarray}

It is readily seen that the planar anchoring condition of $n$ on $\Gamma$
follows from \eqref{bdry_conv}
and the following lower semicontinuity property: for any nonnegative convex function $f:\mathbb R\to\mathbb R_+$, it holds that
\begin{equation}\label{lsc9}
\int_{\Gamma} f\big(n\cdot\nu_{\Gamma}\big)\,d\mathcal{H}^2 \le\liminf_{i\to\infty} \int_{\Gamma_i}f\big(n_i\cdot \nu_{\Gamma_i}\big)\,d\mathcal{H}^2.
\end{equation}
Indeed, if we choose $f(\theta)=\theta^2$, then \eqref{lsc9} and \eqref{slice6} imply that
\begin{eqnarray*}
\int_{\Gamma} \big(n\cdot\nu_{\Gamma}\big)^2\,d\mathcal{H}^2
\le\liminf_{i\to\infty} \int_{\Gamma_i}\big(n_i\cdot \nu_{\Gamma_i}\big)^2\,d\mathcal{H}^2
= \liminf_{i\to\infty}\int_{\Gamma_{i}} \cos^2 \theta_\epsilon \,d\mathcal{H}^2=0.
\end{eqnarray*}
This implies that $n\cdot\nu_{\Gamma}=0$ $\mathcal{H}^2$ a.e. on $\Gamma$.

Now we want to show \eqref{lsc9} as follows.  Define a family of Radon measures
$${\Theta}_i(A)=\mathcal{H}^2(\Gamma_i\cap A) \ {\rm{for}}\ i\ge 1; \ {\Theta}(A)=\mathcal{H}^2(\Gamma\cap A),$$
and
$$\mu_i(A)=\int_A f\big(n_i\cdot \nu_{\Gamma_i}\big)\,d\Theta_i$$
for any measurable set $A\subset\mathbb R^3$.

It is readily seen that there exists a nonnegative Radon measure $\mu$ such that, after passing to a subsequence,
$$\Theta_i\rightharpoonup \Theta\ \ {\rm{and}}\ \ \mu_i\rightharpoonup \mu,$$
as convergence of Radon measures in $\mathbb R^3$. By the Radon-Nikodym theorem, we can decompose
$$\mu=(D_\Theta\mu) \Theta +\mu^s, \ \ {\rm{with}}\ \ \mu^s\perp\Theta.$$
Then we have
$$\int_{A} D_\Theta\mu \,d\Theta\le \mu(A)\le \liminf_{i\to\infty}\mu_i(A), $$
for any open set $A\subset\mathbb R^3$. Hence \eqref{lsc9} follows, if we can show
\begin{equation}\label{lsc1}
f(\nu\cdot\nu_{\Gamma})(x)\le (D_\Theta\mu)(x),  \ \Theta-{\rm{a.e.}}\ x\in {\rm{supp}}(\Theta)=\Gamma.
\end{equation}

From the convexity of $f$, there exist $a_k, b_k\in\mathbb R$ such that
\begin{equation}\label{convexity}
f(\theta)=\sup_{k} (a_k \theta+b_k).
\end{equation}
For $x\in\Gamma$,  we can find $r_j\to 0$ such that for each $j$, it holds that
\begin{equation}\label{fubini}
\lim_{i\rightarrow\infty}\int_{\partial B_{r_j}(x)\cap\Omega_i} n_i\cdot \frac{y-x}{|y-x|}\,d\mathcal{H}^2=
\lim_{i\rightarrow\infty}\int_{\partial B_{r_j}(x)\cap\Omega^+} n\cdot \frac{y-x}{|y-x|}\,d\mathcal{H}^2,
\ \ \ \mu\big(\partial B_{r_j}(x)\big)=0.
\end{equation}
Therefore we have that
\begin{equation}\label{fubini1}
\lim_{i\rightarrow\infty} \int_{B_{r_j}(x)\cap\Gamma_i} f(n_i\cdot\nu_{\Gamma_i})\,d\mathcal{H}^2=
\lim_{i\to\infty}\mu_i(B_{r_j}(x))=\mu(B_{r_j}(x)).
\end{equation}
Similar to \eqref{bdry_conv}, we have that for each $j$,
\begin{eqnarray*}
\lim_{i\to\infty}\int_{\partial (\Omega_i\cap B_{r_j}(x))} n_i\cdot \nu_{\partial(\Omega_i\cap B_{r_j}(x))}\,d\mathcal{H}^2
&=&\lim_{i\to\infty}\int_{\Omega_i\cap B_{r_j}(x)} {\rm{div}} n_i\,dx\\
&=& \int_{\Omega_i\cap B_{r_j}(x)} {\rm{div}} n\,dx\\
&=&\int_{\partial (\Omega^+\cap B_{r_j}(x))} n\cdot \nu_{\partial(\Omega^+\cap B_{r_j}(x))}\,d\mathcal{H}^2.
\end{eqnarray*}
This, combined with \eqref{fubini}, implies that for each $j$ it holds that
$$\lim_{i\to\infty}\int_{\Gamma_i\cap B_{r_j}(x)} n_i\cdot \nu_{\Gamma_i}\,d\mathcal{H}^2
=\int_{\Gamma\cap B_{r_j}(x)} n\cdot \nu_{\Gamma}\,d\mathcal{H}^2.
$$
Recall that  for $\Theta$ a.e. $x\in\Gamma$, it holds that
$$D_\Theta \mu(x)=\lim_{j\rightarrow\infty}\frac{\mu(B_{r_j}(x))}{\Theta(B_{r_j}(x))},$$
and
$$\lim_{j\to\infty} \frac{\int_{\Gamma\cap B_{r_j}(x)} n\cdot\nu_{\Gamma}\,d\mathcal{H}^2}{\Theta(B_{r_j}(x))}=(n\cdot\nu_{\Gamma})(x).$$
Applying \eqref{fubini1}, we obtain that for any fixed $k$,
\begin{eqnarray*}
D_\Theta\mu(x)&=&\displaystyle\lim_{j\rightarrow\infty}\lim_{i\to\infty}\frac{\int_{\Gamma_i\cap B_{r_j}(x)}f(n_i\cdot\nu_{\Gamma_i})\,d\mathcal{H}^2}
{\Theta(B_{r_j}(x))}\\
&\ge& \displaystyle\lim_{j\rightarrow\infty}\lim_{i\to\infty}\frac{\int_{\Gamma_i\cap B_{r_j}(x)} (a_kn_i\cdot\nu_{\Gamma_i}+b_k)\,d\mathcal{H}^2}
{\Theta(B_{r_j}(x))}\\
&\ge& \lim_{j\rightarrow\infty}\lim_{i\to\infty}\Big[a_k \frac{\int_{\Gamma_i\cap B_{r_j}(x)} n_i\cdot\nu_{\Gamma_i}\,d\mathcal{H}^2}
{\Theta(B_{r_j}(x))}+  b_k \frac{\Theta_i(B_{r_j}(x))}{\Theta(B_{r_j}(x))}\Big]\\
&\ge& \lim_{j\rightarrow\infty}\Big[a_k \frac{\int_{\Gamma\cap B_{r_j}(x)} n\cdot\nu_{\Gamma}\,d\mathcal{H}^2}
{\Theta(B_{r_j}(x))}+  b_k \Big]\\
&=& a_k (n\cdot\nu_{\Gamma})(x)+b_k.
\end{eqnarray*}
Taking supremum over $k\ge 1$, we conclude that for $\Theta$ a.e. $x\in\Gamma$,
$$D_\Theta\mu(x)\ge \sup_{k} \big(a_k (n\cdot\nu_{\Gamma})(x)+b_k\big)=f\big(n\cdot\nu_{\Gamma}\big)(x).$$
This yields \eqref{lsc9}.

For $0<\tau<s_+$, let $\mathcal{T}(\tau)\subset\Omega$ be an area minimizing surface spanned
by $\Sigma(\tau)=\big\{x\in\partial\Omega: \ d_{\Gamma}(x)=t\big\}=\partial\Gamma(t)$:
$$
\mathcal{H}^2(\mathcal{T}(\tau))=\min\Big\{\mathcal{H}^2(S): \ S\ \mbox{is an integral}\ 2\mbox{-current in\ } \Omega,
\ \partial S=\Sigma(t)\Big\}.
$$
Then by putting all the above estimates together, we obtain the following lower bound:
\begin{eqnarray}\label{LB7}
\liminf_{\epsilon\to 0}\int_\Omega \widetilde{\mathcal{W}}_\epsilon(s_\epsilon, n_\epsilon, \nabla s_\epsilon, \nabla n_\epsilon)\,dx\ge \frac{2}{\epsilon}\int_0^{s_+} \sqrt{\beta W(\tau)} \mathcal{H}^2(\mathcal{T}(\epsilon\xi_{\epsilon,\gamma}^{-1}(\tau))\,d\tau
+{E}(n, \Omega^+),
\end{eqnarray}
where $E(n, \Omega^+)$ is the Oseen-Frank energy given by
$$
{E}(n, \Omega^+)=s_+^2\int_{\Omega^+}\big(W_{OF}(n, \nabla n)+\alpha|\nabla n|^2\big)\,dx.
$$

Finally, we want to show that $n\in H^1(\Omega^+, \mathbb S^2)$ is a minimizer of the Oseen-Frank
energy $E(\cdot, \Omega^+)$, subject to the boundary condition: $n=g$ on $\Sigma^+$ and
$n\cdot\nu_\Gamma=0$ on $\Gamma$, i.e.,
\begin{equation}\label{energy_min1}
\displaystyle {E}(n, \Omega^+)=\mathcal{D}_A.
\end{equation}

In order to prove \eqref{energy_min1}, we need to show that the leading order term in the lower bound estimate
(\ref{LB7}) exactly matches that in the refined upper bound estimate \eqref{A}, i.e.,
 \begin{equation}\label{LBUB}
 \int_0^{s_+} \sqrt{\beta W(\tau)} \mathcal{H}^2(\Gamma(\epsilon\xi_{\epsilon,\gamma}^{-1}(\tau))\,d\tau
\le \int_0^{s_+} \sqrt{\beta W(\tau)} \mathcal{H}^2(\mathcal{T}(\epsilon\xi_{\epsilon,\gamma}^{-1}(\tau))\,d\tau
+C\epsilon^2.
 \end{equation}

The validity of \eqref{LBUB} is a consequence of the strict stability of $\Gamma$, which ensures

\medskip
\noindent{\it Claim 3}.  {\it  Under the condition that $\Gamma$ is a strictly stable, area minimizing surface, there exist $\eta_0>0$ and $C_0>0$,
depending only on $\Gamma$ and $\Omega$, such that}
\begin{equation}\label{square_comp}
\mathcal{H}^2(\Gamma(\lambda))\le \mathcal{H}^2(\mathcal{T}(\lambda))+C_0\lambda^2, \ \forall \lambda\in [-\eta_0, \eta_0].
\end{equation}
The proof of \eqref{square_comp} is based on the second variation of surface areas and the strict stability of $\Gamma$, we refer the reader to \cite{LPW} page 45-47.

It follows from \eqref{square_comp} that
\begin{eqnarray*}
&&\int_0^{s_+} \sqrt{\beta W(\tau)} \mathcal{H}^2(\Gamma(\epsilon\xi_{\epsilon,\gamma}^{-1}(\tau))\,d\tau\\
&&\le\int_0^{s_+} \sqrt{\beta W(\tau)}\big( \mathcal{H}^2(\mathcal{T}(\epsilon\xi_{\epsilon,\gamma}^{-1}(\tau))
+C(\epsilon\xi_{\epsilon,\gamma}^{-1}(\tau))^2\big)\,d\tau\\
&&\le \int_0^{s_+} \sqrt{\beta W(\tau)}\mathcal{H}^2(\mathcal{T}(\epsilon\xi_{\epsilon,\gamma}^{-1}(\tau))\,d\tau
+C\epsilon^2\int_0^{s_+}\sqrt{\beta W(\tau)} \big(\xi_{\epsilon,\gamma}^{-1}(\tau)\big)^2\,d\tau\\
&&\le \int_0^{s_+} \sqrt{\beta W(\tau)}\mathcal{H}^2(\mathcal{T}(\epsilon\xi_{\epsilon,\gamma}^{-1}(\tau))\,d\tau
+C\epsilon^2,
\end{eqnarray*}
since there exists $C>0$, that is independent of $\epsilon$, such that
$$
\int_0^{s_+}\sqrt{\beta W(\tau)} \big(\xi_{\epsilon,\gamma}^{-1}(\tau)\big)^2\,d\tau\le C\alpha_0.
$$
Hence \eqref{LBUB} holds.
Now it is readily seen that the refined upper bound \eqref{A}, the refined lower bound \eqref{LB7}, and \eqref{LBUB}
imply $E(n,\Omega^+)\le \mathcal{D}_A.$
On the other hand, since $n=g$ on $\Sigma^+$ and $n\cdot\nu_\Gamma=0$ on $\Gamma$, we automatically have
$E(n,\Omega_+)\ge\mathcal{D}_A$. Thus \eqref{energy_min1} holds, and
\begin{equation}\label{perfect_match1}
\int_\Omega \widetilde{\mathcal{W}}_\epsilon(s_\epsilon, n_\epsilon, \nabla s_\epsilon, \nabla n_\epsilon)\,dx=\frac{2}{\epsilon}\int_0^{s_+} \sqrt{\beta W(\tau)} \mathcal{H}^2(\Gamma(\epsilon\xi^{-1}(\tau))\,d\tau
+\mathcal{D}_A+o_\epsilon(1).
\end{equation}
Finally, we claim that (see also \cite{LPW} page 46)
\begin{equation}\label{alternate_expansion}
\frac{2}{\epsilon}\int_0^{s_+} \sqrt{\beta W(\tau)} \mathcal{H}^2(\Gamma(\epsilon\xi_{\epsilon,\gamma}^{-1}(\tau))\,d\tau
=\frac{\alpha_0}{\epsilon}\mathcal{H}^2(\Gamma)+o_\epsilon(1).
\end{equation}
In fact, it is not hard to see that
$$\mathcal{H}^2(\Gamma(\lambda))
= \mathcal{H}^2(\Gamma)+a\lambda+ O(\lambda^2),\ \lambda\in (-\epsilon^\gamma, \epsilon^\gamma),$$
where $a=\frac{d}{d\lambda}\big|_{\lambda=0} \mathcal{H}^2(\Gamma(\lambda))$. Thus we have that
\begin{eqnarray*}
&&\frac{2}{\epsilon}\int_0^{s_+} \sqrt{\beta W(\tau)} \mathcal{H}^2(\Gamma(\epsilon\xi_{\epsilon,\gamma}^{-1}(\tau))\,d\tau=\frac{2}{\epsilon}\int_{-\epsilon^{\gamma-1}}^{\epsilon^{\gamma-1}} \sqrt{\beta W(\xi_{\epsilon,\gamma}(t))} \mathcal{H}^2(\Gamma(\epsilon t))|\xi'_{\epsilon,\gamma}(t)|\,dt\\
&&=\frac{2}{\epsilon}\int_{-\epsilon^{\gamma-1}}^{\epsilon^{\gamma-1}} \sqrt{\beta W(\xi_{\epsilon,\gamma}(t))} \big(\mathcal{H}^2(\Gamma)+a\epsilon t+O(\epsilon^2 t^2)\big)|\xi'_{\epsilon,\gamma}(t)|\,dt\\
&&=\frac{\alpha_0}{\epsilon} \mathcal{H}^2(\Gamma)+2a \int_{-\epsilon^{\gamma-1}}^{\epsilon^{\gamma-1}} \sqrt{\beta W(\xi_{\epsilon,\gamma}(t))} |\xi'_{\epsilon,\gamma}(t)| t\,dt+O(\epsilon)\\
&&=\frac{\alpha_0}{\epsilon} \mathcal{H}^2(\Gamma)+O(\epsilon),
\end{eqnarray*}
where we have the fact that $\sqrt{\beta W(\xi_{\epsilon,\gamma}(t))} |\xi'_{\epsilon,\gamma}(t)|$ is an even function
so that
$$\int_{-\epsilon^{\gamma-1}}^{\epsilon^{\gamma-1}} \sqrt{\beta W(\xi_{\epsilon,\gamma}(t))} |\xi'_{\epsilon,\gamma}(t)| t\,dt=0.$$
Thus we arrive at
\begin{equation}\label{perfect_match2}
\int_\Omega \widetilde{\mathcal{W}}_\epsilon(s_\epsilon, n_\epsilon, \nabla s_\epsilon, \nabla n_\epsilon)\,dx
=\frac{\alpha_0}{\epsilon}\mathcal{H}^2(\Gamma)+\mathcal{D}_A+o_\epsilon(1).
\end{equation}
This proves part ({\bf A}) of Theorem \ref{sharp}. \qed

\subsection{Refined lower bound for the case (B) \fbox{$L_2>0=L_1=L_3, \ L_4^2<4L_2\alpha$}}

The refined lower bound in this case can be done similarly to that of the case (A).
The major difference arises in showing the homeotropic boundary condition of
the limiting map $n$ on the sharp interface $\Gamma$, on which we will focus.

First, from $L_4^2<4L_2\alpha$ we can find a positive number $\mu>0$ such that
$L_4^2\le 4(L_2-\mu)(\alpha-\mu)$. Hence by Cauchy-Schwarz inequality we have
\begin{eqnarray}\label{CS4}
&&|L_4s_\epsilon\nabla s_\epsilon \cdot(\nabla n_\epsilon)n_\epsilon|\nonumber\\
&&=|L_4 (\nabla s_\epsilon-(\nabla s_\epsilon\cdot n_\epsilon)n_\epsilon) \cdot (s_\epsilon\nabla n_\epsilon)|
\nonumber\\
&&\le (L_2 -\mu) |\nabla s_\epsilon-(\nabla s_\epsilon\cdot n_\epsilon)n_\epsilon|^2
+(\alpha-\mu) s_\epsilon^2 |\nabla n_\epsilon|^2\nonumber\\
&&= (L_2 -\mu) |\nabla s_\epsilon\wedge n_\epsilon|^2
+(\alpha-\mu) s_\epsilon^2 |\nabla n_\epsilon|^2.
\end{eqnarray}
This implies that
\begin{eqnarray*}
\widetilde{\mathcal{W}}_\epsilon(s_\epsilon, n_\epsilon,\nabla s_\epsilon, \nabla n_\epsilon)
&=&( \beta |\nabla s_\epsilon|^2+\frac{1}{\epsilon^2}W(s_\epsilon))+s_\epsilon^2 W_{OF}(n_\epsilon, \nabla n_\epsilon)\\
&+&\big(\alpha s_\epsilon^2|\nabla n_\epsilon|^2
+L_2|\nabla s_\epsilon\wedge n_\epsilon|^2
+L_4(s_\epsilon\nabla s_\epsilon)\cdot (\nabla n_\epsilon) n_\epsilon\big)\\
&\ge& \beta |\nabla s_\epsilon|^2+\frac{1}{\epsilon^2}{W(s_\epsilon)}
+\mu s_\epsilon^2 |\nabla n_\epsilon|^2+
\mu|\nabla s_\epsilon\wedge n_\epsilon|^2\\
&=&\big[(\beta+\mu \sin^2\theta_\epsilon) |\nabla s_\epsilon|^2+\frac{1}{\epsilon^2}W(s_\epsilon)\big]
+\mu s_\epsilon^2 |\nabla n_\epsilon|^2\\
&\ge& \frac{2}{\epsilon}|\nabla s_\epsilon|\sqrt{W(s_\epsilon)}\sqrt{\beta+\mu\sin^2\theta_\epsilon}
+\mu s_\epsilon^2 |\nabla n_\epsilon|^2,
\end{eqnarray*}
where $\sin\theta_\epsilon=\frac{\nabla s_\epsilon}{|\nabla s_\epsilon|}\wedge n_\epsilon$.

Note that we also have
\begin{eqnarray*}
\widetilde{\mathcal{W}}_\epsilon(s_\epsilon, n_\epsilon,\nabla s_\epsilon, \nabla n_\epsilon)
&=&( \beta |\nabla s_\epsilon|^2+\frac{1}{\epsilon^2}W(s_\epsilon))+s_\epsilon^2 W_{OF}(n_\epsilon, \nabla n_\epsilon)\\
&&+\big(\alpha s_\epsilon^2|\nabla n_\epsilon|^2
+L_2|\nabla s_\epsilon\wedge n_\epsilon|^2
+L_4(s_\epsilon\nabla s_\epsilon)\cdot (\nabla n_\epsilon) n_\epsilon\big)\\
&\ge& \beta |\nabla s_\epsilon|^2+\frac{1}{\epsilon^2}{W(s_\epsilon)}
+\mu s_\epsilon^2 |\nabla n_\epsilon|^2+
\mu|\nabla s_\epsilon\wedge n_\epsilon|^2\\
&\ge& \frac{2}{\epsilon}|\nabla s_\epsilon|\sqrt{\beta W(s_\epsilon)}
+\mu s_\epsilon^2 |\nabla n_\epsilon|^2+\mu|\nabla s_\epsilon\wedge n_\epsilon|^2.
\end{eqnarray*}

Hence,  by the co-area formula, this implies that
\begin{eqnarray}\label{LB3.1}
&&\int_\Omega \widetilde{\mathcal{W}}_\epsilon(s_\epsilon, n_\epsilon,\nabla s_\epsilon, \nabla n_\epsilon)\,dx\nonumber\\
&&\ge \int_\Omega \big(\frac{2}{\epsilon}|\nabla s_\epsilon|\sqrt{W(s_\epsilon)}\sqrt{\beta+\mu\sin^2\theta_\epsilon}+\mu s_\epsilon^2 |\nabla n_\epsilon|^2\big)\,dx\nonumber\\
&&\ge\mu \int_\Omega s_\epsilon^2 |\nabla n_\epsilon|^2\,dx
+\frac{2\mu_*}{\epsilon}\int_0^{s_+}\sqrt{W(\tau)} \int_{\partial^*\mathcal{S}_\epsilon(\tau)} \sin^2\theta_\epsilon\,d\mathcal{H}^2 \,d\tau\nonumber\\
&&\ \ \ +\frac{2}{\epsilon}\int_0^{s_+}\sqrt{\beta W(\tau)} \mathcal{H}^2(\partial^*\mathcal{S}_\epsilon(\tau))\,d\tau,
\end{eqnarray}
where $\displaystyle\mu_*=\frac{\mu}{\sqrt{\mu+\beta}+\sqrt{\beta}}$,
and
\begin{eqnarray}\label{LB3.2}
&&\int_\Omega \widetilde{\mathcal{W}}_\epsilon(s_\epsilon, n_\epsilon,\nabla s_\epsilon, \nabla n_\epsilon)\,dx\nonumber\\
&&\ge\mu \int_\Omega \big(s_\epsilon^2 |\nabla n_\epsilon|^2+|\nabla s_\epsilon\wedge n_\epsilon|^2\big)\,dx
+\frac{2}{\epsilon}\int_0^{s_+}\sqrt{\beta W(\tau)} \mathcal{H}^2(\partial^*\mathcal{S}_\epsilon(\tau))\,d\tau.
\end{eqnarray}

Matching \eqref{LB3.1} and \eqref{LB3.2} with the upper bound \eqref{B}, we conclude that
\begin{equation}\label{UB5.1}
\mu \int_\Omega (s_\epsilon^2 |\nabla n_\epsilon|^2+|\nabla s_\epsilon\wedge n_\epsilon|^2)\,dx
\le \mathcal{D}_B+C+o_\epsilon(1),
\end{equation}
\begin{equation}\label{UB6.1}
\frac{\mu_*}{\epsilon}\int_0^{s_+}\sqrt{W(\tau)}
 \int_{\partial^*\mathcal{S}_\epsilon(\tau)} \sin^2\theta_\epsilon\,d\mathcal{H}^2 \,d\tau
 \le \mathcal{D}_B+C+o_\epsilon(1),
\end{equation}
and
\begin{equation}\label{UB7.1}
\frac{2}{\epsilon}\int_0^{s_+}\sqrt{\beta W(\tau)} \mathcal{H}^2(\partial^*\mathcal{S}_\epsilon(\tau))\,d\tau
\le \frac{1}{\epsilon}\alpha_0\mathcal{H}^2(\Gamma)+\mathcal{D}_B+o_\epsilon(1).
\end{equation}
As in the previous section, for any small $\delta>0$ there exist $C>0$, independent of
$\epsilon$, and $\delta_\epsilon\in (\delta, 2\delta)$ such that
 \begin{equation}\label{slice10}
\displaystyle \int_{\partial^*\mathcal{S}_\epsilon(s_+-{\delta_\epsilon})} \sin^2 \theta_\epsilon \,d\mathcal{H}^2
\le C\frac{\epsilon}{\delta},
\end{equation}
\begin{equation}\label{slice11}
\displaystyle \mathcal{H}^2(\partial^*\mathcal{S}_\epsilon(s_+ -{\delta_\epsilon}))\le
\mathcal{H}^2(\Gamma)+C\frac{\epsilon}{\delta}.
 \end{equation}
Moreover, from \eqref{UB6.1} we have
\begin{equation}\label{n-UB10}
\int_{\Omega_{\epsilon,\delta_\epsilon}^+} \big(|\nabla n_\epsilon|^2+|\nabla s_\epsilon\wedge n_\epsilon|^2\big)\,dx
\le C.
\end{equation}

As in the previous section, there exists $\epsilon_i\to 0$ and $\delta_i\in (\epsilon_i^\frac14, 2\epsilon_i^\frac14)$
such that $\Omega_i=\Omega^+_{\epsilon_i,\delta_i}$ converges to $\Omega^+$ weakly in $BV(\R^3)$,
$$\Gamma_i=\partial\Omega_i\lfloor \Omega\rightharpoonup \Gamma=\partial\Omega^+\lfloor\Omega,$$
as convergence of measures and integral $2$-currents.
Furthermore, there exists a map $n\in H^1(\Omega^+,\mathbb S^2)$ such that
$$n_{\epsilon_i}\chi_{\Omega_i} \rightharpoonup n\chi_{\Omega^+} \ \ {\rm{in}}\ \ {\rm{BV}}(\Omega),
\ \ {\rm{and}}\ \  s_{\epsilon_i}\chi_{\Omega_i} \rightarrow s_+\chi_{\Omega^+} \ \ {\rm{in}}\ \ L^2(\Omega),
$$
and hence
$$
\nabla n_{\epsilon_i}\chi_{\Omega_i}\rightharpoonup \nabla n\chi_{\Omega^+},
\ {\rm{and}}\ \ s_{\epsilon_i}\nabla n_{\epsilon_i}\chi_{\Omega_i}\rightharpoonup s_+\nabla n\chi_{\Omega^+}
\ \ {\rm{in}}\ \ L^2(\Omega).
$$
As a consequence of these weak convergences and \eqref{n-UB10}, we can deduce
$$\nabla s_{\epsilon_i}\wedge n_{\epsilon_i}\rightharpoonup 0 \ \ {\rm{in}}\ \ L^2(\Omega).$$

Similar to the proof of \eqref{lsc3} and \eqref{lsc4}, we can obtain
\begin{eqnarray}\label{lsc11}
\int_{\Omega^+} s_{+}^2W_{OF}(n, \nabla n)\,dx
\le\liminf_{\epsilon_i\rightarrow 0}
\int_{\Omega_i} s_{\epsilon_i}^2W_{OF}(n_{\epsilon_i}, \nabla n_{\epsilon_i})\,dx,
\end{eqnarray}
and
\begin{equation}\label{lsc12}
\int_{\Omega^+} \alpha s_+^2|\nabla n|^2\,dx
\le\liminf_{\epsilon_i\to 0}
\int_{\Omega_i}\big(\alpha s_{\epsilon_i}^2|\nabla n_{\epsilon_i}|^2+L_2|\nabla s_{\epsilon_i}\wedge n_{\epsilon_i}|^2
+L_4(s_{\epsilon_i}\nabla s_{\epsilon_i}\cdot (\nabla n_\epsilon)n_\epsilon\big)\,dx.
\end{equation}
Adding \eqref{lsc11} with \eqref{lsc12}, we arrive at
\begin{eqnarray}\label{lsc13}
&&E(n,\Omega^+)=\int_{\Omega^+} s_+^2\big(W_{OF}(n,\nabla n)+\alpha|\nabla n|^2\big)\,dx\\
&&\le\liminf_{\epsilon_i\rightarrow 0}\int_{\Omega}\big(s_{\epsilon_i}^2 W_{OF}(n_{\epsilon_i}, \nabla n_{\epsilon_i})+\alpha s_{\epsilon_i}^2|\nabla n_{\epsilon_i}|^2
+L_2|\nabla s_{\epsilon_i}\wedge n_{\epsilon_i}|^2
+L_4 s_{\epsilon_i} \nabla s_{\epsilon_i}\cdot (\nabla n_{\epsilon_i}) n_{\epsilon_i}\big)\,dx.\nonumber
\end{eqnarray}

Now we want to show the homeotropic condition of $n$ on $\Gamma$, i.e.,
\begin{equation}\label{hemotropic}
n\wedge\nu_\Gamma=0 \ {\rm{on}}\ \Gamma.
\end{equation}
In order to show \eqref{hemotropic}, we first want to prove
\begin{equation}\label{wedge}
\int_{\Gamma} n\wedge\nu_{\Gamma}\,d\mathcal{H}^2
=\lim_{i\to\infty} \int_{\Gamma_i} n_i\wedge\nu_{\Gamma_i}\,d\mathcal{H}^2.
\end{equation}
In fact, since
$$\nabla\times n_i\rightharpoonup \nabla\times n \ {\rm{in}}\ L^2(\Omega^+),
\ {\rm{and}}\ \chi_{\Omega_i}\rightarrow \chi_{\Omega^+} \ {\rm{in}}\ L^2(\mathbb R^3),$$
we must have
$$\lim_{i\rightarrow\infty}\int_{\Omega_i} \nabla\times n_i\,dx=\int_{\Omega^+} \nabla\times n\,dx.$$
This, combined with the divergence theorem, implies that
\begin{eqnarray*}
&&\int_{\partial\Omega_i\cap\partial\Omega}
n_i\wedge \nu_{\partial\Omega_i}\,d\mathcal{H}^2+\int_{\Gamma_i}
n_i\wedge \nu_{\Gamma_i}\,d\mathcal{H}^2\\
&&=\int_{(\partial\Omega_i\cap\partial\Omega)\cup\Gamma_i}
n_i\wedge \nu_{\partial\Omega_i}\,d\mathcal{H}^2=\int_{\Omega_i} \nabla\times n_i\,dx\\
&&\rightarrow\int_{\Omega^+} \nabla\times n\,dx
=\int_{(\partial\Omega^+\cap\partial\Omega)\cup\Gamma} n\wedge\nu_{\partial\Omega^+}\,d\mathcal{H}^2\\
&&=\int_{\partial\Omega^+\cap\partial\Omega} n\wedge\nu_{\partial\Omega}\,d\mathcal{H}^2+
\int_{\Gamma} n\wedge\nu_{\Gamma}\,d\mathcal{H}^2.
\end{eqnarray*}
On other hand, it follows from the assumption on $g_{\epsilon_i}$ on $\partial\Omega$ that
\begin{eqnarray*}
&&\int_{\partial\Omega_i\cap\partial\Omega}
n_i\wedge \nu_{\partial\Omega_i}\,d\mathcal{H}^2=\int_{\partial\Omega_i\cap\partial\Omega}
g_{\epsilon_i}\wedge \nu_{\partial\Omega_i}\,d\mathcal{H}^2\\
&&\rightarrow\int_{\partial\Omega^+\cap\partial\Omega}
g\wedge \nu_{\partial\Omega}\,d\mathcal{H}^2=\int_{\partial\Omega^+\cap\partial\Omega}
n\wedge \nu_{\partial\Omega}\,d\mathcal{H}^2.
\end{eqnarray*}
Thus \eqref{wedge} holds.

From \eqref{wedge}, we can use a similar argument to show that for any nonnegative convex function
$f:\mathbb R\to\mathbb R_+$, it holds that
\begin{equation}\label{hemotropic1}
\int_{\Gamma} f(n\wedge \nu_{\Gamma})\,d\mathcal{H}^2
\le\liminf_{i\rightarrow\infty} \int_{\Gamma_i} f(n_i\wedge \nu_{\Gamma_i})\,d\mathcal{H}^2.
\end{equation}
In particular, this implies that
\begin{equation}\label{hemotropic2}
\int_{\Gamma} |n\wedge \nu_{\Gamma}|^2\,d\mathcal{H}^2
\le\liminf_{i\rightarrow\infty} \int_{\Gamma_i} |n_i\wedge \nu_{\Gamma_i}|^2\,d\mathcal{H}^2
=\liminf_{i\rightarrow\infty} \int_{\Gamma_i} \sin^2\theta_{\epsilon_i}\,d\mathcal{H}^2=0.
\end{equation}
This yields \eqref{hemotropic}.

For the convenience of readers, we will sketch the proof of \eqref{hemotropic1} as follows.
Define a family of Radon measures
$$\mathcal{M}_i(A)=\int_{A} f(n_i\wedge\nu_{\Gamma_i})\,d\Theta_i,
$$
for any measurable set $A\subset\mathbb R^3$.

Without loss of generality, we can assume that
there exists a Radon measure $\mathcal{M}$ in $\mathbb R^3$ such that
$$\mathcal{M}_i\rightharpoonup \mathcal{M}$$
as convergence of Radon measures on $\mathbb R^3$. Again by Radon-Nikodym theorem, we can decompose
$$\mathcal{M}=(D_\Theta \mathcal{M}) \Theta +\mathcal{M}^s, \ \mathcal{M}^s\perp \mathcal{M}.$$
Hence for any open set $O\subset\mathbb R^3$, it holds that
$$\int_{O} D_\Theta \mathcal{M} \,d\Theta\le \mathcal{M}(O)\le\liminf_{i\to\infty} \mathcal{M}_i(O).$$
Now we want to show
\begin{equation}\label{point_bound}
D_\Theta \mathcal{M}(x)\ge f(n\wedge\nu_{\Gamma})(x), \ \Theta\ {\rm{a.e.}}\ x\in\Gamma.
\end{equation}

First, it is not hard to see that for any $x\in\Gamma$, it holds that for $L^1$ a.e. $r>0$,
$$\lim_{i\to\infty} \int_{\partial B_r(x)\cap\Omega_i} n_i \wedge \frac{y-x}{|y-x|}\,d\mathcal{H}^2
=\int_{\partial B_r(x)\cap\Omega^+} n \wedge \frac{y-x}{|y-x|}\,d\mathcal{H}^2.$$
This, combined with
\begin{eqnarray*}
&&\int_{\partial(B_r(x)\cap\Omega_i)} n_i\wedge\nu_{\partial (B_r(x)\cap\Omega_i)}=\int_{B_r(x)\cap\Omega_i} \nabla\times n_i\\
&&\to \int_{B_r(x)\cap\Omega^+} \nabla\times n=\int_{\partial(B_r(x)\cap\Omega^+)} n\wedge\nu_{\partial (B_r(x)\cap\Omega^+)},
\end{eqnarray*}
yields that
$$
\lim_{i\to\infty} \int_{B_r(x)\cap\Gamma_i} n_i \wedge \nu_{\Gamma_i}\,d\mathcal{H}^2
=\int_{B_r(x)\cap\Gamma} n \wedge \nu_{\Gamma}\, d\mathcal{H}^2.
$$
It is readily seen that
$$\mathcal{M}(\partial B_r(x))=0$$
holds  for $L^1$ a.e. $r>0$.

Therefore, for any given $x\in \Gamma$,
\begin{eqnarray}\label{lsc10}
&&\frac{\mathcal{M}(B_r(x))}{\Theta(B_r(x))}=\lim_{i\to\infty}\frac{\mathcal{M}_i(B_r(x))}{\Theta(B_r(x))}\nonumber\\
&&\ge \lim_{i\to\infty}\frac{\int_{B_r(x)\cap \Gamma_i} (a_k n_i\wedge \nu_{\Gamma_i} +b_k)\,d\mathcal{H}^2}{\Theta(B_r(x))}\nonumber\\
&&= a_k \lim_{i\rightarrow \infty}\displaystyle\frac{\int_{B_r(x)\cap\Gamma_i} n_i\wedge\nu_{\Gamma_i}\,d\mathcal{H}^2}{\Theta(B_r(x))} +b_k\nonumber\\
&&= a_k \displaystyle\frac{\int_{B_r(x)\cap\Gamma} n\wedge\nu_{\Gamma}\,d\mathcal{H}^2}{\Theta(B_r(x))} +b_k
\end{eqnarray}
holds for any $k\ge 1$ and $L^1$ a.e. $r>0$.

Since
$$D_\Theta \mathcal{M}(x)=\lim_{r\to 0}\frac{\mathcal{M}(B_r(x))}{\Theta(B_r(x))},$$
and
$$(n\wedge\nu_{\Gamma})(x)=\lim_{r\to 0} \frac{\int_{B_r(x)\cap\Gamma} n\wedge\nu_{\Gamma}\,d\mathcal{H}^2}{\Theta(B_r(x))} $$
hold for $\Theta$ a.e. $x\in\Gamma$, after passing to the limit in
\eqref{lsc10} we obtain that for $\Theta$ a.e. $x\in\Gamma$,
$$D_\Theta \mathcal{M}(x)\ge a_k (n\wedge\nu_{\Gamma})(x)+b_k, \ \forall k\ge 1.
$$
Taking supremum over $k\ge 1$ and using \eqref{convexity}, this yields \eqref{point_bound}.

From \eqref{hemotropic}, we conclude that by
$$E(n; \ \Omega^+)\ge\mathcal{D}_B,$$
which, combined with \eqref{LB3.2} and \eqref{lsc13}, implies that
\begin{eqnarray}\label{OFLB}
\int_\Omega \widetilde{\mathcal{W}}_\epsilon(s_\epsilon, n_\epsilon,\nabla s_\epsilon, \nabla n_\epsilon)\,dx
\ge \frac{2}{\epsilon}\int_0^{s_+}\sqrt{\beta W(\tau)} \mathcal{H}^2(\mathcal{S}_\epsilon(\tau))\,d\tau+
\mathcal{D}_B+o_\epsilon(1).
\end{eqnarray}
This, combined with the inequality \eqref{LBUB} and the upper bound \eqref{B}, further implies
that
$$E(n;\ \Omega^+)=\mathcal{D}_B,$$
\begin{equation}\label{LBUB1}
\int_\Omega \widetilde{\mathcal{W}}_\epsilon(s_\epsilon, n_\epsilon,\nabla s_\epsilon, \nabla n_\epsilon)\,dx
=\frac{\alpha_0}{\epsilon}\mathcal{H}^2(\Gamma)+\mathcal{D}_B+o_\epsilon(1).
\end{equation}
The part ({\bf B}) of Theorem \ref{sharp} is proven. \qed

\subsection{Refined lower bound for the case (C) \fbox{$L_1=L_2=L_3=L_4=0$}}

This case is the easiest among the three cases we discuss in this paper. In fact, by
\eqref{C} and direct calculations, we obtain  that
\begin{eqnarray*}
&&\frac{\alpha_0}{\epsilon}\mathcal{H}^2(\Gamma)+\mathcal{D}_C+o_\epsilon(1)\\
 &&\ge \int_{\Omega}\widetilde{\mathcal{W}}_\epsilon(s_\epsilon, n_\epsilon,\nabla s_\epsilon, \nabla n_\epsilon)\,dx\nonumber\\
&&=\int_{\Omega} \big[(\beta |\nabla s_\epsilon|^2+\frac{1}{\epsilon^2}W(s_\epsilon))+(s_\epsilon^2 W_{OF}(n_\epsilon, \nabla n_\epsilon)+\alpha s_\epsilon^2|\nabla n_\epsilon|^2)\big]\,dx\\
&&\ge \frac{2}{\epsilon}\int_0^{s_+} \sqrt{\beta {W(\tau)}}\mathcal{H}^2(\mathcal{S}(\epsilon\xi_{\epsilon,\gamma}^{-1}(\tau))\,d\tau
+(s_+-\delta)^2 \int_{\Omega^+_{\epsilon,\delta}} (W_{OF}(n_\epsilon, \nabla n_\epsilon)+\alpha |\nabla n_\epsilon|^2)\,dx\\
&&\ge  \frac{\alpha_0}{\epsilon}\mathcal{H}^2(\Gamma)-C
+(s_+-\delta)^2 \int_{\Omega^+_{\epsilon,\delta}} (W_{OF}(n_\epsilon, \nabla n_\epsilon)+\alpha |\nabla n_\epsilon|^2)\,dx.
\end{eqnarray*}
This implies that for any $\delta>0$, there exists $\delta_\epsilon\in (\delta, 2\delta)$
such that
\begin{equation}\label{slice13}
\int_{\Omega^+_{\epsilon,\delta_\epsilon}} (W_{OF}(n_\epsilon, \nabla n_\epsilon)+\alpha |\nabla n_\epsilon|^2)\,dx
\le \frac{1}{(s_+-\delta_\epsilon)^2}\big(\mathcal{D}_C+C+o_\epsilon(1)\big),
\end{equation}
and
\begin{equation}\label{slice14}
\displaystyle \mathcal{H}^2(\mathcal{S}_\epsilon(s_+ -{\delta_\epsilon}))\le
\mathcal{H}^2(\Gamma)+C\frac{\epsilon}{\delta}.
 \end{equation}
As in the previous two cases, we can argue as follows.
For $\epsilon_i\rightarrow 0$, there exists $\delta_{\epsilon_i}\in (\epsilon_i^\frac12, 2\epsilon_i^\frac12)$
such that $\Omega_i=\Omega_{\epsilon_i, \delta_{\epsilon_i}}^+\rightarrow \Omega^+$ weakly
 in ${\rm{BV}}(\R^3)$,
 $$\mathcal{S}_\epsilon(s_+ -{\delta_{\epsilon_i}})=\partial\Omega_i\lfloor \Omega
 \rightharpoonup \Gamma=\partial\Omega^+\lfloor\Omega,$$
weakly converges as measures and integral currents. Furthermore, there exists $n\in H^1(\Omega^+)$
with $n=g$ on $\Sigma^+$ such that
$$
s_{\epsilon_i} \nabla n_{\epsilon_i}\chi_{\Omega_i}\rightharpoonup s_+\nabla n \chi_{\Omega^+}
\ \ {\rm{in}}\ \ L^2(\Omega).
$$
By the lower semicontinuity, we then have
\begin{equation}\label{lsc20}
\mathcal{D}_C\le\int_{\Omega^+} s_+^2(W_{OF}(n,\nabla n)+\alpha |\nabla n|^2)\,dx
\le\liminf_{\epsilon_i\to 0} \int_\Omega \big(s_{\epsilon_i}^2 W_{OF}(n_{\epsilon_i}, \nabla n_{\epsilon_i})+\alpha s_{\epsilon_i}^2|\nabla n_{\epsilon_i}|^2\big)\,dx.
\end{equation}
On the other hand, it follows from the inequality \eqref{LBUB} that
$$ \frac{2}{\epsilon}\int_0^{s_+} \sqrt{\beta {W(\tau)}}\mathcal{H}^2(\mathcal{S}(\epsilon\xi_{\epsilon,\gamma}^{-1}(\tau))\,d\tau
\ge \frac{\alpha_0}{\epsilon} \mathcal{H}^2(\Gamma)-C\epsilon.$$
Thus we can conclude that
$$E(n;\ \Omega^+)=\mathcal{D}_C,$$
and
$$
\int_{\Omega}\widetilde{\mathcal{W}}_\epsilon(s_\epsilon, n_\epsilon,\nabla s_\epsilon, \nabla n_\epsilon)\,dx
=\frac{\alpha_0}{\epsilon} \mathcal{H}^2(\Gamma)+\mathcal{D}_C+o_\epsilon(1).
$$
This establishes the conclusion of part (C) in Theorem \ref{sharp}. \qed

\section{Proof of Theorem \ref{sharp1}}

In this section, we will consider the asymptotic expansion of
$\mathcal{E}(s_\epsilon, n_\epsilon)$ for $\Omega=\R^3$ and prove Theorem \ref{sharp1}.
We will first provide a sharp upper bound estimate for the cases (A), (B), and (C).

\subsection{Refined upper bounds}
\setcounter{equation}{0}
\setcounter{theorem}{3}

The constructions are similar to those in the section 2.2.2, except that we need to use an almost
minimal $1$-dimensional connecting orbit in the fast transition region of width $O(\epsilon^\gamma)$
around $\partial B_1$.

Let $\xi_{\epsilon,\gamma}\in C^\infty([-\epsilon^\gamma, \epsilon^\gamma])$ be given by \eqref{xi_cond}.
Define $\hat{s}_\epsilon:\R^3\to \R_+$ by  letting
\begin{equation*}
\hat{s}_\epsilon(x)=\begin{cases} s_+ & |x|\le 1, \\
\xi_{\epsilon,\gamma}(\frac{|x|-(1+\epsilon^\gamma)}{\epsilon}) & 1\le |x|\le 1+2\epsilon^\gamma,\\
0 & |x|\ge 1+2\epsilon^\gamma.
\end{cases}
\end{equation*}
Then we have
$$\big\{x\in\R^3:  \hat{s}_\epsilon(x)\ge {s_+}\big\}= B_{1},$$
and
\begin{eqnarray}\label{ext_est21}
&&\int_{\R^3} \big(|\nabla \hat{s}_\epsilon|^2+\frac{1}{\epsilon^2}W(\hat{s}_\epsilon)\big)\,dx
=\frac{1}{\epsilon}\int_{-\epsilon^{\gamma-1}}^{\epsilon^{\gamma-1}}
\big(\beta|\xi_{\epsilon,\gamma}'(t)|^2+W(\xi_{\epsilon,\gamma}(t))\big)\mathcal{H}^2(\partial B_{1+\epsilon^\gamma+\epsilon t})\,dt\nonumber\\
&&=\frac{4\pi}{\epsilon}\int_{-\epsilon^{\gamma-1}}^{\epsilon^{\gamma-1}}
\big(\beta|\xi_{\epsilon,\gamma}'(t)|^2+W(\xi_{\epsilon,\gamma}(t))\big)[(1+\epsilon^\gamma)^2+2\epsilon(1+\epsilon^\gamma)t+\epsilon^2 t^2]\,dt\nonumber\\
&&=\frac{\alpha_0}{\epsilon}\mathcal{H}^2(\partial B_1)+C\epsilon^{2\gamma-1}+o_\epsilon(1)
=\frac{\alpha_0}{\epsilon}\mathcal{H}^2(\partial B_1)+o_\epsilon(1),
\end{eqnarray}
since $\gamma>\frac12$.

\medskip
Now we divide the discussion into three different cases:
\begin{itemize}
\item [Case (A)] $L_1>L_2=L_4=0$.
Let $n\in H^1(B_1, \mathbb S^2)$ be such that $n\cdot\nu_{\partial B_1}=0$ on $\partial B_1$, and
\begin{equation}
E(n; B_1)=\mathcal{D}_A,
\end{equation}
where $\mathcal{D}_A$ is given by \eqref{OFM10} and \eqref{bdry10}.
We may assume that $\bar{n}\in H^1(B_{1+2\epsilon^\gamma},\mathbb S^2)$ be such that
$\bar{n}=n$ in $B_1$, and
$$\int_{B_{1+2\epsilon^\gamma}\setminus B_1}|\nabla\bar{n}|^2\,dx=o_\epsilon(1).$$
Then we have
\begin{eqnarray}\label{ext_est22}
&&\int_{\R^3} \hat{s}_\epsilon^2 \big(W_{OF}(\bar{n}, \nabla \bar{n})+\alpha |\nabla\bar{n}|^2\big)\,dx\nonumber\\
&&\le \int_{B_1} {s}_+^2 \big(W_{OF}({n}, \nabla {n})+\alpha |\nabla{n}|^2\big)\,dx
+Cs_+^2\int_{B_{1+2\epsilon^\gamma}\setminus B_1}|\nabla \bar{n}|^2\,dx
=\mathcal{D}_A+o_\epsilon(1).
\end{eqnarray}
While \begin{eqnarray}\label{ext_est23}
&&\big|\int_{\R^3} \big(L_1(\nabla\hat{s}_\epsilon\cdot \bar{n})^2
+L_3(\nabla\hat{s}_\epsilon\cdot \bar{n})(\hat{s}_\epsilon{\rm{div}}\bar{n})\big)\,dx\big|\nonumber\\
&&=\big|\int_{B_{1+2\epsilon^\gamma}\setminus B_1} \big(L_1(\nabla\hat{s}_\epsilon\cdot \bar{n})^2
+L_3(\nabla\hat{s}_\epsilon\cdot \bar{n})(\hat{s}_\epsilon{\rm{div}}\bar{n})\big)\,dx\big|\nonumber\\
&&\le C\big(\int_{B_{1+2\epsilon^\gamma}\setminus B_1}|\nabla \hat{s}_\epsilon\cdot\bar{n}|^2\,dx
+\int_{B_{1+2\epsilon^\gamma}\setminus B_1} |\nabla\bar{n}|^2\,dx\big)\nonumber\\
&&\le o_\epsilon(1)+C\epsilon^{-2}\int_{B_{1+2\epsilon^\gamma}\setminus B_1}|\xi'_{\epsilon,\gamma}|^2\big(\frac{|x|-(1+\epsilon^\gamma)}{\epsilon}\big)
\big|\bar{n}(x)-\bar{n}(\frac{x}{|x|})\big|^2\,dx\nonumber\\
&&\le o_\epsilon(1)+C\int_{B_{1+2\epsilon^\epsilon}\setminus B_1}\big|\nabla\bar{n}\big|^2\,dx=o_\epsilon(1).
\end{eqnarray}
Therefore by putting together \eqref{ext_est21}, \eqref{ext_est22}, and \eqref{ext_est23}, we arrive at
\begin{eqnarray}\label{SUB3.1}
\limsup_{\epsilon\to 0}\int_{\R^3} \widetilde{\mathcal{W}}({s}_\epsilon, {n}_\epsilon,
\nabla {s}_\epsilon, \nabla {n}_\epsilon)\,dx&\le&
\limsup_{\epsilon\to 0}\int_{\R^3} \widetilde{\mathcal{W}}(\hat{s}_\epsilon, \bar{n},
\nabla \hat{s}_\epsilon, \nabla \bar{n})\,dx\nonumber\\
&\le& \frac{\alpha_0}{\epsilon}\mathcal{H}^2(\partial B_1)+\mathcal{D}_A.
\end{eqnarray}

\item [Case (B)] $L_2>L_1=L_3=0$. Let $n\in H^1(B_1, \mathbb S^2)$ be such that $n\wedge\nu_{\partial B_1}=0$ on $\partial B_1$, and
\begin{equation}
E(n; B_1)=\mathcal{D}_B,
\end{equation}
where $\mathcal{D}_B$ is given by \eqref{OFM20} and \eqref{bdry20}.
We may assume that $\bar{n}\in H^1(B_{1+2\epsilon^\gamma},\mathbb S^2)$ be such that
$\bar{n}=n$ in $B_1$, and
$$\int_{B_{1+2\epsilon^\gamma}\setminus B_1}|\nabla\bar{n}|^2\,dx=o_\epsilon(1).$$
Then we have
\begin{eqnarray}\label{ext_est24}
&&\int_{\R^3} \hat{s}_\epsilon^2 \big(W_{OF}(\bar{n}, \nabla \bar{n})+\alpha |\nabla\bar{n}|^2\big)\,dx\nonumber\\
&&\le \int_{B_1} {s}_+^2 \big(W_{OF}({n}, \nabla {n})+\alpha |\nabla{n}|^2\big)\,dx
+Cs_+^2\int_{B_{1+2\epsilon^\gamma}\setminus B_1}|\nabla \bar{n}|^2\,dx
=\mathcal{D}_B+o_\epsilon(1).
\end{eqnarray}
While \begin{eqnarray}\label{ext_est25}
&&\big|\int_{\R^3} \big(L_2|\nabla\hat{s}_\epsilon\wedge \bar{n}|^2
+L_4(\hat{s}_\epsilon\nabla\hat{s}_\epsilon)(\nabla\bar{n})\bar{n}\big)\,dx\big|\nonumber\\
&&=\big|\int_{B_{1+2\epsilon^\gamma}\setminus B_1} \big(L_2|\nabla\hat{s}_\epsilon\wedge \bar{n}|^2
+L_4(\hat{s}_\epsilon\nabla\hat{s}_\epsilon)(\nabla\bar{n})\bar{n}\big)\,dx\big|\nonumber\\
&&\le C\big(\int_{B_{1+2\epsilon^\gamma}\setminus B_1}|\nabla \hat{s}_\epsilon\wedge\bar{n}|^2\,dx
+\int_{B_{1+2\epsilon^\gamma}\setminus B_1} |\nabla\bar{n}|^2\,dx\big)\nonumber\\
&&\le o_\epsilon(1)+C\epsilon^{-2}\int_{B_{1+2\epsilon^\gamma}\setminus B_1}|\xi'_{\epsilon,\gamma}|^2\big(\frac{|x|-(1+\epsilon^\gamma)}{\epsilon}\big)
\big|\bar{n}(x)-\bar{n}(\frac{x}{|x|})\big|^2\,dx\nonumber\\
&&\le o_\epsilon(1)+C\int_{B_{1+2\epsilon^\epsilon}\setminus B_1}\big|\nabla\bar{n}\big|^2\,dx=o_\epsilon(1).
\end{eqnarray}
Therefore by putting together \eqref{ext_est21}, \eqref{ext_est24}, and \eqref{ext_est25}, we arrive at
\begin{eqnarray}\label{SUB4.1}
\limsup_{\epsilon\to 0}\int_{\R^3} \widetilde{\mathcal{W}}({s}_\epsilon, {n}_\epsilon,
\nabla {s}_\epsilon, \nabla {n}_\epsilon)\,dx&\le&
\limsup_{\epsilon\to 0}\int_{\R^3} \widetilde{\mathcal{W}}(\hat{s}_\epsilon, \bar{n},
\nabla \hat{s}_\epsilon, \nabla \bar{n})\,dx\nonumber\\
&\le& \frac{\alpha_0}{\epsilon}\mathcal{H}^2(\partial B_1)+\mathcal{D}_B.
\end{eqnarray}

\item [Case (C)] $L_1=L_2=L_3=L_4=0$. Let $n\equiv (0,0,1)\in\mathbb S^2$. Then by \eqref{ext_est21}
we see that
\begin{eqnarray}\label{SUB5.1}
\int_{\R^3} \widetilde{\mathcal{W}}({s}_\epsilon, {n}_\epsilon,
\nabla {s}_\epsilon, \nabla {n}_\epsilon)\,dx&\le&
\int_{\R^3} \widetilde{\mathcal{W}}(\hat{s}_\epsilon, \bar{n},
\nabla \hat{s}_\epsilon, \nabla \bar{n})\,dx\nonumber\\
&=&\int_{\R^3} \big(|\nabla \hat{s}_\epsilon|^2+\frac{1}{\epsilon^2}W(\hat{s}_\epsilon)\big)\,dx\nonumber\\
&\le&\frac{\alpha_0}{\epsilon}\mathcal{H}^2(\partial B_1)+o_\epsilon(1).
\end{eqnarray}
\end{itemize}

\subsection{Refined lower bounds}
In this subsection, we will sketch the proof of a sharp a lower bound estimate for the cases (A), (B), and (C).
The ideas are similar to those presented in the section 3 for bounded domain cases, except that we will work on
the entire space $\R^3$ where we only have the weak compactness property of ${\rm{BV}}(\R^3)$ locally
in $\R^3$. We will focus on
the case (A), and only sketch the cases (B) and (C).
\begin{itemize}
\item [Case (A)] First, as in the discussion of case (A) in the section 3.2, there exists $\mu>0$ such that
\begin{eqnarray}\label{LB4.0}
&&\int_{\R^3} \widetilde{\mathcal{W}}_\epsilon(s_\epsilon, n_\epsilon,\nabla s_\epsilon, \nabla n_\epsilon)\,dx
\nonumber\\
&&\ge \frac{2}{\epsilon}\int_0^{s_+}\sqrt{\beta W(\tau)} \mathcal{H}^2(\partial^*\mathcal{S}_\epsilon(\tau))\,d\tau
+\mu\int_{\R^3}\big(s_\epsilon^2|\nabla n_\epsilon|^2+(\nabla s_\epsilon\cdot n_\epsilon)^2\big) \,dx,
\end{eqnarray}
and
 \begin{eqnarray}\label{LB4.1}
&&\int_{\R^3} \widetilde{\mathcal{W}}_\epsilon(s_\epsilon, n_\epsilon,\nabla s_\epsilon, \nabla n_\epsilon)\,dx\nonumber\\
&&\ge \mu  \int_{\R^3} s_\epsilon^2 |\nabla n_\epsilon|^2\,dx
+\frac{2\mu_*}{\epsilon}\int_0^{s_+}\sqrt{\beta W(\tau)} \int_{\partial^*\mathcal{S}_\epsilon(\tau)} \cos^2\theta_\epsilon\,d\mathcal{H}^2 \,d\tau\nonumber\\
&&\ \ +\frac{2}{\epsilon}\int_0^{s_+}\sqrt{\beta W(\tau)} \mathcal{H}^2(\partial^*\mathcal{S}_\epsilon(\tau))\,d\tau,
\end{eqnarray}
where $\mathcal{S}_\epsilon(\tau)=\big\{x\in\R^3: \ s_\epsilon(x)\ge\tau\big\}$,
$\cos\theta_\epsilon=\frac{\nabla s_\epsilon}{|\nabla s_\epsilon|}\cdot n_\epsilon$,
and $\displaystyle\mu_*=\frac{\mu}{\sqrt{\beta+\mu}+\sqrt{\beta}}>0$.

Notice that by the isoperimetric inequality (see, e.g., Case (C) below) we have
\begin{equation}\label{LB4.2}
\frac{2}{\epsilon}\int_0^{s_+}\sqrt{\beta W(\tau)} \mathcal{H}^2(\partial^*\mathcal{S}_\epsilon(\tau))\,d\tau
\ge \frac{\alpha_0}{\epsilon}\mathcal{H}^2(\partial B_1).
\end{equation}
It follows from \eqref{LB4.0}, \eqref{LB4.1}, \eqref{LB4.2} and \eqref{SUB3.1} that
there exists $\tau_{\epsilon}\in (0, s_+)$ such that
\begin{equation}\label{UB5.2}
 \int_{\mathcal{S}_\epsilon(\tau_{\epsilon})} (s_\epsilon^2 |\nabla n_\epsilon|^2+|\nabla s_\epsilon\cdot n_\epsilon|^2)\,dx \le C(\mu,\ \mathcal{D}_A),
\end{equation}
\begin{equation}\label{UB6.2}
 \int_{\mathcal{S}_\epsilon(\tau_{\epsilon})} \cos^2\theta_\epsilon\,d\mathcal{H}^2
 \le C\epsilon,
\end{equation}
and
\begin{equation}\label{UB7.2}
 \mathcal{H}^2(\partial^*\mathcal{S}_\epsilon(\tau_\epsilon))
\le \mathcal{H}^2(\partial B_1)+o_\epsilon(1).
\end{equation}
By the isoperimetric inequality and the volume constraint condition, we have that
\begin{equation}\label{volume1}
|B_1|\le  \big|\mathcal{S}_\epsilon(\tau_\epsilon)\big|\le  |B_1|
\big(\frac{\big|\mathcal{H}^2(\partial^*\mathcal{S}_\epsilon(\tau_\epsilon))\big|}{\mathcal{H}^2(\partial B_1)}\big)^\frac32\le |B_1|(1+o_\epsilon(1)).
\end{equation}

To simplify the presentation, we denote by $E_\epsilon=\mathcal{S}_\epsilon(\tau_\epsilon)$. Although
$E_\epsilon$ may not converge in $L^1(\R^3)$ due to the non-compactness of $\R^3$, we can apply
the quantitative stability theorem by Fusco-Maggi-Pratelli \cite{FMP} (see also \cite{Maggi1})
to show that $E_\epsilon$ does converge
in $L^1(\R^3)$ after suitable translations. In fact, if we set the isoperimetric deficit and Fraenkel asymmetry by
\begin{equation}\label{deficit}
A(E_\epsilon)=\frac{\mathcal{H}^2(E_\epsilon)}{3|B_1|^\frac13|E_\epsilon|^\frac23}-1,
\ {\rm{and}}\ \ \delta(E_\epsilon)=\inf\Big\{ \frac{|E_\epsilon\Delta B_r(x)|}{|E_\epsilon|}: \ |B_r(x)|={|E_\epsilon|}\Big\},
\end{equation}
then it follows from \cite{FMP} that
\begin{equation}\label{qs}
A(E_\epsilon)\le C\delta(E_\epsilon)^\frac12\le o_\epsilon(1),
\end{equation}
where we have used \eqref{UB7.2} and \eqref{volume1} in the second inequality of \eqref{qs}.
From \eqref{volume1}, there exist $x_\epsilon\in \R^3$ and $r_\epsilon=1+o_\epsilon(1)$ such that
after passing to a subsequence,
$$|E_\epsilon\Delta B_{r_\epsilon}(x_\epsilon)|\rightarrow 0, \ \ {\rm{as}}\ \ \epsilon\to 0,$$
or equivalently,
\begin{equation}\label{domain_conv}
\widehat{E}_\epsilon\equiv E_\epsilon\setminus\{x_\epsilon\}\rightarrow B_1 \ \ {\rm{in}}\ \ L^1(\R^3).
\end{equation}
Since the problem is invariant under translations, for simplicity we may assume that $x_\epsilon=0$
so that $E_\epsilon=\widehat{E}_\epsilon$. The rest of argument can be done almost identically to
the case (A) of Theorem \ref{sharp} presented in section 3.2. For instance, we can show that there
exists $n\in H^1(B_1,\mathbb S^2)$ such that
$$s_\epsilon\rightarrow s_+ \ \ {\rm{in}}\ \  L^2(B_1), \ n_\epsilon\rightharpoonup n\ \ {\rm{in}} \ \ H^1(B_1),$$
and
$$n(x)\cdot x=0 \ \  {\rm{on}}\ \ \partial B_1.$$
Moreover, by the lower semicontinuity we have that
\begin{eqnarray*}
&&\mathcal{D}_A\le E(n; B_1)=s_+^2\int_{B_1} (W_{OF}(n,\nabla n)+\alpha|\nabla n|^2)\,dx\\
&&\le\liminf_{\epsilon\to 0} \int_{E_\epsilon}
\big(s_\epsilon^2W_{OF}(n_\epsilon, \nabla n_\epsilon)+\alpha |\nabla n_\epsilon|^2
+L_1(\nabla s_{\epsilon}\cdot n_{\epsilon})^2
+L_3(\nabla s_{\epsilon}\cdot n_{\epsilon})s_{\epsilon} {\rm{div}} n_{\epsilon}\big)\,dx\\
&&=\liminf \Big(\int_{\R^3} \widetilde{\mathcal{W}}(s_\epsilon, n_\epsilon, \nabla s_\epsilon,\nabla n_\epsilon)\,dx
-\int_{\R^3} \big(|\nabla {s}_\epsilon|^2+\frac{1}{\epsilon^2}W({s}_\epsilon)\big)\,dx\Big)\\
&&\le \mathcal{D}_A.
\end{eqnarray*}
Hence
$$E(n; B_1)=\mathcal{D}_A,  \ \ \int_{\R^3} \widetilde{\mathcal{W}}(s_\epsilon, n_\epsilon, \nabla s_\epsilon, \nabla n_\epsilon)\,dx\ge\frac{\alpha_0}{\epsilon}\mathcal{H}^2(\partial B_1)+\mathcal{D}_A+o_\epsilon(1).$$

\item [Case (B)] First, as in the discussion of case (B) in the section 3.2, there exists $\mu>0$ such that
\begin{eqnarray}\label{LB5.1}
&&\int_{\R^3} \widetilde{\mathcal{W}}_\epsilon(s_\epsilon, n_\epsilon,\nabla s_\epsilon, \nabla n_\epsilon)\,dx\nonumber\\
&&\ge
\frac{2\mu_*}{\epsilon}\int_0^{s_+}\sqrt{W(\tau)} \int_{\partial^*\mathcal{S}_\epsilon(\tau)} \sin^2\theta_\epsilon\,d\mathcal{H}^2 +\frac{2}{\epsilon}\int_0^{s_+}\sqrt{\beta W(\tau)} \mathcal{H}^2(\partial^*\mathcal{S}_\epsilon(\tau))\,d\tau,
\end{eqnarray}
where $\displaystyle\mu_*=\frac{\mu}{\sqrt{\mu+\beta}+\sqrt{\beta}}$ and $\sin^2\theta_\epsilon
=|\frac{\nabla s_\epsilon}{|\nabla s_\epsilon|}\wedge n_\epsilon|^2$,
and
\begin{eqnarray}\label{LB5.2}
&&\int_{\R^3} \widetilde{\mathcal{W}}_\epsilon(s_\epsilon, n_\epsilon,\nabla s_\epsilon, \nabla n_\epsilon)\,dx\nonumber\\
&&\ge\mu \int_{\R^3} \big(s_\epsilon^2 |\nabla n_\epsilon|^2+|\nabla s_\epsilon\wedge n_\epsilon|^2\big)\,dx
+\frac{2}{\epsilon}\int_0^{s_+}\sqrt{\beta W(\tau)} \mathcal{H}^2(\partial^*\mathcal{S}_\epsilon(\tau))\,d\tau.
\end{eqnarray}

As in the case (A), we can find $\tau_\epsilon\in (0, s_+)$ such that
\begin{equation}\label{UB5.3}
 \int_{\mathcal{S}_\epsilon(\tau_{\epsilon})} (s_\epsilon^2 |\nabla n_\epsilon|^2+|\nabla s_\epsilon\wedge n_\epsilon|^2)\,dx \le C(\mu,\ \mathcal{D}_B),
\end{equation}
\begin{equation}\label{UB6.3}
 \int_{\mathcal{S}_\epsilon(\tau_{\epsilon})} \sin^2\theta_\epsilon\,d\mathcal{H}^2
 \le C\epsilon,
\end{equation}
and
\begin{equation}\label{UB7.3}
 \mathcal{H}^2(\partial^*\mathcal{S}_\epsilon(\tau_\epsilon))
\le \mathcal{H}^2(\partial B_1)+o_\epsilon(1).
\end{equation}

As in the discussion of Case (A) above, we can apply the quantitative stability theorem of \cite{FMP}
to conclude that   after passing to a subsequence,
$$S_\epsilon(\tau_\epsilon)\rightarrow B_1 \ \ {\rm{in}}\ \ L^1(\R^3).$$
Furthermore, by an argument similar to the case (B) of Theorem \ref{sharp} presented in the section 3.3,
there exists a $n\in H^1(B_1, \mathbb S^2)$ such that
$$s_\epsilon\to s_+ \ \ {\rm{in}}\ \ L^2(B_1), \ \ n_\epsilon\rightharpoonup n\ \ {\rm{in}}\ \ H^1(B_1),$$
$$n(x)\wedge x=0 \ \ {\rm{on}}\ \ \partial B_1,$$
and by the lower semicontinuity,
\begin{eqnarray*}
&&\mathcal{D}_B\le E(n; B_1)=s_+^2\int_{B_1} (W_{OF}(n,\nabla n)+\alpha|\nabla n|^2)\,dx\\
&&\le\liminf_{\epsilon\to 0} \int_{E_\epsilon}
\big(s_\epsilon^2W_{OF}(n_\epsilon, \nabla n_\epsilon)+\alpha |\nabla n_\epsilon|^2
+L_2|\nabla s_{\epsilon}\wedge n_{\epsilon}|^2
+L_4 s_{\epsilon} \nabla s_{\epsilon}\cdot (\nabla n_{\epsilon}) n_{\epsilon}\big)\,dx\\
&&=\liminf_{\epsilon\to 0} \Big(\int_{\R^3} \widetilde{\mathcal{W}}(s_\epsilon, n_\epsilon, \nabla s_\epsilon,\nabla n_\epsilon)\,dx
-\int_{\R^3} \big(|\nabla {s}_\epsilon|^2+\frac{1}{\epsilon^2}W({s}_\epsilon)\big)\,dx\Big)\\
&&\le \lim_{\epsilon\to 0} \Big(\frac{\alpha_0}{\epsilon}\mathcal{H}^2(\partial B_1)+ \mathcal{D}_B+o_\epsilon(1)
-\frac{\alpha_0}{\epsilon}\mathcal{H}^2(\partial B_1)\Big)=\mathcal{D}_B.
\end{eqnarray*}
Hence
$$E(n; B_1)=\mathcal{D}_B,  \ \ \int_{\R^3} \widetilde{\mathcal{W}}(s_\epsilon, n_\epsilon, \nabla s_\epsilon, \nabla n_\epsilon)\,dx\ge\frac{\alpha_0}{\epsilon}\mathcal{H}^2(\partial B_1)+\mathcal{D}_B+o_\epsilon(1).$$

\item [Case (C)]  This case is the simplest, since it reduces to the iso-perimetric inequality:
\begin{eqnarray*}
&&\int_{\R^3} \widetilde{\mathcal{W}}({s}_\epsilon, {n}_\epsilon,
\nabla {s}_\epsilon, \nabla {n}_\epsilon)\,dx
\ge \int_{\R^3}\big(|\nabla {s}_\epsilon|^2+\frac{1}{\epsilon^2}W({s}_\epsilon)\big) \,dx\nonumber\\
&&\ge\frac{2}{\epsilon}\int_{\R^3}\sqrt{\beta W(s_\epsilon)}|\nabla s_\epsilon|\,dx\nonumber\\
&&\ge \frac{2}{\epsilon}\int_{0}^{s_+} \sqrt{\beta W(\tau)}\mathcal{H}^2\big(\partial^*\big\{x\in\R^3: \ s_\epsilon(x)\ge\tau\big\}\big)\,d\tau.
\end{eqnarray*}
Since
$$\Big|\big\{x\in\R^3: \ s_\epsilon(x)\ge\tau\big\}\Big|\ge \Big|\big\{x\in\R^3: \ s_\epsilon(x)\ge s_+\big\}\Big|
=\big|B_1\big|, \ \forall 0<\tau<s_+,$$
it follows from the isoperimetric inequality in $\R^3$ that for all $0<\tau<s_+,$
\begin{eqnarray*}
\mathcal{H}^2\big(\partial^*\big\{x\in\R^3: \ s_\epsilon(x)\ge\tau\big\}\big)
&\ge& (36\pi)^\frac13 \Big|\big\{x\in\R^3: \ s_\epsilon(x)\ge\tau\big\}\Big|^\frac23\\
&\ge& (36\pi)^\frac13|B_1|^\frac23=\mathcal{H}^2(\partial B_1).
\end{eqnarray*}
Hence we obtain that
\begin{equation}
\int_{\R^3} \widetilde{\mathcal{W}}({s}_\epsilon, {n}_\epsilon,
\nabla {s}_\epsilon, \nabla {n}_\epsilon)\,dx\ge
\frac{\alpha_0}{\epsilon}\mathcal{H}^2(\partial B_1).
\end{equation}
\end{itemize}

\bigskip
\noindent{\it Completion of Proof of Theorem \ref{sharp1}}. It is readily seen that Theorem \ref{sharp1}
follows by combining the arguments from both section 4.1 and section 4.2. \qed

\bigskip
\bigskip
\bigskip
\noindent{\bf Acknowledgements}. The first author is partially supported by NSF grant DMS1955249. The second author is partially supported by NSF grant DMS1764417.

\end{document}